\documentclass[12pt]{amsart}

\usepackage{accents}
\usepackage{appendix}
\usepackage{amsfonts}
 \usepackage{amsaddr}
\usepackage{amsmath}
\usepackage{amssymb}	
\usepackage{amsthm,bm}
\usepackage{array,booktabs,multirow}
\usepackage{braket}
\usepackage{centernot}
\usepackage{cite}
\usepackage{comment}
\usepackage{dsfont}
\usepackage{mathalfa}
\usepackage[shortlabels]{enumitem}
\usepackage{etoolbox}
\usepackage{float}
\usepackage[hang, flushmargin]{footmisc}
\usepackage{latexsym}
\usepackage{lipsum}
\usepackage{needspace}
\usepackage{tikz}
\usepackage{hyperref}
\usetikzlibrary{matrix,arrows}
\usepackage{multicol}

\makeatletter
\def\namedlabel#1#2{\begingroup
   \def\@currentlabel{#2}%
   \label{#1}\endgroup
}
\makeatother

\newcommand{\floor}[1]{\left\lfloor #1 \right\rfloor}
\newcommand{\ceil}[1]{\left\lceil #1 \right\rceil}

\catcode`,\active

\catcode`\,12

\theoremstyle{plain}
\newtheorem{thm}{Theorem}[section]
\newtheorem{cor}[thm]{Corollary}
\newtheorem{lem}[thm]{Lemma}
\newtheorem{prop}[thm]{Proposition}

\theoremstyle{definition}

\newtheorem{defn}[thm]{Definition}
\newtheorem{exam}[thm]{Example}

\theoremstyle{remark}

\numberwithin{equation}{section}

\setlist[enumerate,1]{leftmargin=2em}

\def\A{\mathbf A}
\def\C{\mathbb C}

\def\H{\mathcal H_q}
\def\I{\mathbf I}

\def\L{\mathcal L}

\def\M{\mathbf M}
\def\N{\mathbb N}
\def\P{\mathbf P}
\def\T{\mathbf T}
\def\F{\mathbb F}
\def\Z{\mathbb Z}

\def\t{\tau_{x_0}^{x_1}}
\def\Px{\Phi_{x_0}^{x_1}}
\def\E{\mathcal E}

\def\U{U_q(\mathfrak{sl}_2)}
\def\W{\mathfrak{W}_q}
\def\e{\varepsilon}

\newcommand{\subsetdot}{\subset\mathrel{\mkern-7mu}\mathrel{\cdot}\,}

\title[The Clebsch--Gordan coefficients of $\U$ and Grassmann graphs]{An imperceptible connection between \\ the Clebsch--Gordan coefficients of $\U$ and \\ the Terwilliger algebras of Grassmann graphs}

\author{Hau-Wen Huang}
\address{
Department of Mathematics,
National Central University,
Chung-Li 32001, 
Taiwan
}
\email{hauwenh@math.ncu.edu.tw}

\begin{document}
\begin{abstract}
The Clebsch--Gordan coefficients of $U(\mathfrak{sl}_2)$ are  expressible in terms of Hahn polynomials. The phenomenon can be explained by an algebra homomorphism $\natural$ from the universal Hahn algebra $\mathcal H$ into 
$U(\mathfrak{sl}_2)\otimes U(\mathfrak{sl}_2)$. Let $\Omega$ denote a finite set of size $D$ and $2^\Omega$ denote the power set of $\Omega$.
It is generally known that $\C^{2^\Omega}$ supports a $U(\mathfrak{sl}_2)$-module. 
Let $k$ denote an integer with $0\leq k\leq D$ and fix a $k$-element subset $x_0$ of $\Omega$.
By identifying $\C^{2^\Omega}$ with $\C^{2^{\Omega\setminus x_0}}\otimes \C^{2^{x_0}}$ this induces a $U(\mathfrak{sl}_2)\otimes U(\mathfrak{sl}_2)$-module structure on $\C^{2^\Omega}$ denoted by $\C^{2^\Omega}(x_0)$. Pulling back via 
$\natural$ 
the $U(\mathfrak{sl}_2)\otimes U(\mathfrak{sl}_2)$-module $\C^{2^\Omega}(x_0)$ forms an $\mathcal H$-module. 
When $1\leq k\leq D-1$ the $\mathcal H$-module $\C^{2^\Omega}(x_0)$ enfolds the Terwilliger algebra of the Johnson graph $J(D,k)$ with respect to $x_0$. This result connects these two seemingly irrelevant topics: The Clebsch--Gordan coefficients of $U(\mathfrak{sl}_2)$ and the Terwilliger algebras of Johnson graphs. Unfortunately some steps break down in the $q$-analog case. By making detours, the imperceptible connection between the Clebsch--Gordan coefficients of $\U$ and the Terwilliger algebras of Grassmann graphs is successfully disclosed in this paper. 
\end{abstract}

\maketitle



\allowdisplaybreaks

\section{Introduction}\label{s:introduction}

The following conventions are in effect throughout the paper: 
Let $\N$ denote the set of nonnegative integers. 
Let $\Z$ denote the ring of integers. 
Let $\C$ denote the complex number field. 
The unadorned tensor product $\otimes$ is meant to be over $\C$. 
Suppose that $\Omega$ is a given vector space. The notation  $\L(\Omega)$ stands for the set of all subspaces of $\Omega$. 
Note that $(\L(\Omega),\subseteq)$ is a poset.
We adopt the symbol $\subsetdot$ to represent the covering relation of $(\L(\Omega),\subseteq)$.  For any $k\in \N$ let $\L_k(\Omega)$ denote the set of all $x\in \L(\Omega)$ with $\dim x=k$. 
For any set $X$ the notation $\C^X$ stands for the vector space over $\C$ that has a basis $X$. 
An {\it algebra} is meant to be a unital associative algebra. 
A subset of an algebra is called a {\it subalgebra} if the subset itself forms an algebra with the same operations and contains the unit of the parent algebra. 
An {\it algebra homomorphism} is meant to be a unital algebra homomorphism.  We always assume that $q$ is a nonzero complex number which is not a root of $1$. 
For any two elements $x,y$ of an algebra over $\C$, the bracket $[x,y]$ and the $q$-bracket $[x,y]_q$ are defined as 
\begin{align*}
[x,y]&=xy-yx,
\\
[x,y]_q&=qxy - q^{-1}yx.
\end{align*}
The $q$-analog $[n]_q$ of $n\in \Z$ is defined as 
$$
[n]_q=\frac{q^n-q^{-n}}{q-q^{-1}}.
$$ 
Note that the value of a vacuous summation is set to be $0$ and the value of a vacuous product is set to be $1$.

The subconstituent algebra of an association scheme was proposed in the pioneering papers \cite{TerAlgebraI, TerAlgebraII, TerAlgebraIII} of Terwilliger. 
The algebra is a combinatorial analog of the centralizer algebra of the stabilizer of a point $x_0$ in the automorphism group of an association scheme.
It has been renamed the Terwilliger algebra and studied for the hypercubes in \cite{hypercube2002, Lthypercube},  the Hamming graphs in \cite{Hamming2006, Huang:CG&Hamming, Hamming:2021}, the Johnson graphs in \cite{Huang:CG&Johnson, Johnson:2019, Johnson:2017, Johnson:2014, Johnson:2021}, the odd graphs in \cite{odd:2013, odd:2024}, the Grassmann graphs in \cite{Lee:2018, Watanabe2020}, the Doob graphs in \cite{Doob1997,Doob2023}, the dual polar graphs in \cite{Boyd:2013,Tanaka:2018,dualpolar:2022}, the connection with coding theory \cite{coding&Hamming, GTerAlg2005, semidefinite2006} and so on. 
In \cite{hypercube2002} it was proved that the Terwilliger algebra of a hypercube is a homomorphic image of $U(\mathfrak{sl}_2)$. 
Inspired by this result, the present author comes up with a series of papers on the connections between $U(\mathfrak{sl}_2)$ and the Terwilliger algebras of 
various $Q$-polynomial distance-regular graphs:  
In \cite{halved:2023,halved:2024} the Terwilliger algebras of halved cubes are studied through the even subalgebra of $U(\mathfrak{sl}_2)$. 
In \cite{odd:2024} the Terwilliger algebras of odd graphs are shown to be related to a skew group ring of $\Z/2\Z$ over $U(\mathfrak{sl}_2)$. In \cite{Huang:CG&Hamming} the Clebsch--Gordan rule for $U(\mathfrak{sl}_2)$ is applied to the study of Terwilliger algebras of Hamming graphs. In \cite{Huang:CG&Johnson} the Terwilliger algebras of Johnson graphs are revisited from the viewpoint of the algebra behind the Clebsch--Gordan coefficients of $U(\mathfrak{sl}_2)$.

To explain the motivation of this work, we give a quick sketch of \cite{Huang:CG&Johnson}: The Clebsch--Gordan coefficients of $U(\mathfrak{sl}_2)$ are expressible in terms of Hahn polynomials. The phenomenon can be explained by an algebra homomorphism $\natural$ from the universal Hahn algebra $\mathcal H$ into 
$U(\mathfrak{sl}_2)\otimes U(\mathfrak{sl}_2)$ 
in terms of the comultiplication $\Delta$ of $U(\mathfrak{sl}_2)$ \cite[Theorem 1.5]{Huang:CG&Johnson}. 
Let $\Omega$ denote a finite set of size $D$ and $2^\Omega$ denote the power set of $\Omega$. 
It is commonly known that $\C^{2^\Omega}$ has a  $U(\mathfrak{sl}_2)$-module structure \cite[Theorem 5.1]{Huang:CG&Johnson}. 
Let $k$ denote an integer with $0\leq k\leq D$ and fix a $k$-element subset $x_0$ of $\Omega$. 
There is a linear isomorphism $\C^{2^\Omega}\to \C^{2^{\Omega\setminus x_0}}\otimes \C^{2^{x_0}}$ given by $x\mapsto (x\setminus x_0)\otimes (x\cap x_0)$ for all $x\in 2^\Omega$. 
The $U(\mathfrak{sl}_2)$-module $\C^{2^\Omega}$ and the $U(\mathfrak{sl}_2)\otimes U(\mathfrak{sl}_2)$-module $\C^{2^{\Omega\setminus x_0}}\otimes \C^{2^{x_0}}$ are linked by the following property: The diagram 
\begin{table}[H]
\centering
\begin{tikzpicture}
\matrix(m)[matrix of math nodes,
row sep=4em, column sep=6em,
text height=1.5ex, text depth=0.25ex]
{
\C^{2^\Omega}
&\C^{2^{\Omega\setminus x_0}}\otimes \C^{2^{x_0}}
\\
\C^{2^\Omega}
&\C^{2^{\Omega\setminus x_0}}\otimes \C^{2^{x_0}}
\\
};
\path[->,font=\scriptsize,>=angle 90]
(m-1-1) edge node[left] {$X$} (m-2-1)
(m-1-1) edge node[above] {$\iota(x_0)$} (m-1-2)
(m-2-1) edge node[below] {$\iota(x_0)$} (m-2-2)
(m-1-2) edge node[right] {$\Delta(X)$} (m-2-2);
\end{tikzpicture}
\end{table}
\noindent is commutative for all $X\in U(\mathfrak{sl}_2)$ \cite[Lemma 5.4]{Huang:CG&Johnson}. 
By identifying $\C^{2^\Omega}$ with $\C^{2^{\Omega\setminus x_0}}\otimes \C^{2^{x_0}}$ via $\iota(x_0)$ this induces a $U(\mathfrak{sl}_2)\otimes U(\mathfrak{sl}_2)$-module $\C^{2^\Omega}$ denoted by $\C^{2^\Omega}(x_0)$. Pulling back via $\natural$ 
the $U(\mathfrak{sl}_2)\otimes U(\mathfrak{sl}_2)$-module $\C^{2^\Omega}(x_0)$ becomes an $\mathcal H$-module. 
In light of \cite[Lemma 5.4]{Huang:CG&Johnson}, it is easy to see that, when $1\leq k\leq D-1$, the $\mathcal H$-module $\C^{2^\Omega}(x_0)$ enfolds the Terwilliger algebra of the Johnson graph $J(D,k)$ with respect to $x_0$ \cite[Theorem 5.9]{Huang:CG&Johnson}. 
This result connects these two seemingly irrelevant topics: The Clebsch--Gordan coefficients of $U(\mathfrak{sl}_2)$ and the Terwilliger algebras of Johnson graphs. 
Regrettably, several steps malfunction in the $q$-analog case. By hurdling the challenges, the invisible connection between the Clebsch--Gordan coefficients of $\U$ and the Terwilliger algebras of Grassmann graphs is successfully unveiled in this paper.
We now begin this work by recalling the definition of $\U$:

\begin{defn}\label{defn:U}
The {\it quantum algebra $\U$ of $\mathfrak{sl}_2$} is an algebra over $\C$ generated by $E,F,K^{\pm 1}$ subject to the relations
\begin{gather}
KK^{-1}=K^{-1}K=1,
\label{U1}
\\
[E,K]_q=
[K,F]_q=0,
\label{U2}
\\
[E,F]=\frac{K-K^{-1}}{q-q^{-1}}.
\label{U3}
\end{gather}
\end{defn}

The  element 
\begin{gather}\label{Lambda}
\Lambda=(q-q^{-1})^2EF+q^{-1} K+ q K^{-1}
\end{gather}
is central in $\U$. Note that $\Lambda$ is a scalar multiple of the Casimir element of $\U$. In this paper we simply call $\Lambda$ the {\it Casimir element} of $\U$.

\begin{lem}
\label{lem:U_comultiplication}
There exists a unique algebra homomorphism $\Delta:\U\to \U\otimes \U$ given by 
\begin{eqnarray*}
\Delta(E) &= & E\otimes 1+K^{-1}\otimes E,
\\
\Delta(F) &= & F\otimes K+1\otimes F,
\\
\Delta(K^{\pm 1}) &= & K^{\pm 1}\otimes K^{\pm 1}.
\end{eqnarray*}
\end{lem}
\begin{proof}
Using Definition \ref{defn:U} yields that $E\otimes 1+K^{-1}\otimes E$, $F\otimes K+1\otimes F$, $K^{\pm 1}\otimes K^{\pm 1}$ satisfy the relations (\ref{U1})--(\ref{U3}). Hence the existence of $\Delta$ follows. Since the algebra $\U$ is generated by $E,F,K^{\pm 1}$ the uniqueness of $\Delta$ follows.
\end{proof}

The algebra $\U$ along with the algebra homomorphism $\Delta:\U\to \U\otimes \U$ and the algebra homomorphism $\U\to \C$ given by 
\begin{eqnarray*}
E &\mapsto & 0,
\qquad 
F \;\; \mapsto \;\; 0,
\qquad 
K^{\pm 1} \;\; \mapsto \;\; 1,
\end{eqnarray*}
is a {\it coalgebra}. 
The map $\Delta:\U\to \U\otimes \U$ is called the {\it comultiplication} of $\U$ and the map $\U\to \C$ is called the {\it counit} of $\U$.
To put a coalgebra structure on $\U$ 
the comultiplication has other choices. For example see \cite[Section 3.8]{jantzen}.

\begin{defn}
[Definition 2.1, \cite{Huang:CG}]
\label{defn:H}
The {\it universal $q$-Hahn algebra} $\H$ is an algebra over $\C$ generated by $A,B,C$ and the relations assert that each of 
\begin{align}\label{abc}
\frac{[B,C]_q}{q^2-q^{-2}}+A,
\qquad 
[C,A]_q,
\qquad 
\frac{[A,B]_q}{q^2-q^{-2}}+C
\end{align}
commutes with $A,B,C$. Let $\alpha,\beta,\gamma$ denote the central elements of $\H$ obtained from (\ref{abc}) by multiplying by $q+q^{-1},\frac{1}{q-q^{-1}},q+q^{-1}$ respectively. In other words 
\begin{align}
\alpha &= \frac{[B,C]_q}{q-q^{-1}}+(q+q^{-1})A,
\label{alpha}
\\
\beta &= \frac{[C,A]_q}{q-q^{-1}},
\label{beta}
\\
\gamma &= \frac{[A,B]_q}{q-q^{-1}}+(q+q^{-1})C.
\label{gamma}
\end{align}
\end{defn}

A $q$-analog of \cite[Theorem 1.5]{Huang:CG&Johnson} was appeared in \cite[Theorem 2.9]{Huang:CG} earlier than \cite[Theorem 1.5]{Huang:CG&Johnson}.
With respect to the current comultiplication $\Delta$ of $\U$ the algebraic treatment of the Clebsch--Gordan coefficients of $\U$ is as follows:

\begin{thm}
[Theorem 2.9, \cite{Huang:CG}]
\label{thm:H->U2}
There exists a unique algebra homomorphism $\natural:\H\to \U\otimes \U$ that sends 
\begin{eqnarray}
A &\mapsto & 1\otimes K^{-1},
\label{A:natural}
\\
B &\mapsto &\Delta(\Lambda),
\label{B:natural}
\\
C &\mapsto & 
K^{-1}\otimes 1
-
q(q-q^{-1})^2 
E\otimes F K^{-1}.
\label{C:natural}
\end{eqnarray}
The algebra homomorphism $\natural$ sends 
\begin{eqnarray}
\alpha  &\mapsto & 1\otimes \Lambda+(\Lambda\otimes 1) \Delta(K^{-1}),
\label{alpha:natural}
\\
\beta &\mapsto & \Delta(K^{-1}),
\label{beta:natural}
\\
\gamma &\mapsto & \Lambda\otimes 1+(1\otimes \Lambda) \Delta(K^{-1}).
\label{gamma:natural}
\end{eqnarray}
\end{thm}

To tie Theorem \ref{thm:H->U2} to the Terwilliger algebras of Grassmann graphs, we propose an algebra $\W$ 
and show that there is a surjective algebra homomorphism 
$$
\flat:\W\to \U\otimes \U.
$$ 
Please see Section \ref{s:W} for details. 
In Section \ref{s:lift1} we display a lift 
$$
\widetilde{\Delta}:\U\to \W
$$ 
of the comultiplication $\Delta:\U\to \U\otimes \U$ across $\flat$. 
In Section \ref{s:Wcentral} we reveal two central elements $\Lambda_1$ and $\Lambda_2$ of $\W$ which are the lifts of the central elements $\Lambda\otimes 1$ and  $1\otimes \Lambda$ of $\U\otimes \U$ across $\flat$.
In Section \ref{s:lift2} we apply the lift $\widetilde{\Delta}$ of $\Delta$ and the elements $\Lambda_1,\Lambda_2$ of $\W$ to depict a lift 
$$
\widetilde{\natural}:\H\to \W
$$ 
of the algebra homomorphism $\natural:\H\to \U\otimes \U$ across $\flat$. Therefore every $\W$-module is converted into an $\H$-module by pulling back via $\widetilde{\natural}$. 
In Section \ref{s:Wmodule} we classify the finite-dimensional irreducible $\W$-modules up to isomorphism. In Section \ref{s:Wreducible} we state the necessary and sufficient conditions for finite-dimensional $\W$-modules as completely reducible.
In Section \ref{s:Hmodule} we classify the finite-dimensional irreducible $\H$-modules up to isomorphism. In Section \ref{s:HWmodule} any finite-dimensional irreducible $\W$-module is decomposed into a direct sum of irreducible $\H$-modules.

Let $\Omega$ denote a vector space of finite dimension $D$ over a finite field $\F$. Fix a subspace $x_0$ of $\Omega$ and set the parameter $q$ to be $\sqrt{|\F|}$. 
In Section \ref{s:Dunkl} we recall the triple coordinate system $(\L(\Omega)_{x_0},\subseteq)$ for $(\L(\Omega),\subseteq)$ introduced by Dunkl in \cite[Section 4]{Dunkl77}.
In Section \ref{s:L(Omega)} we give a $\U$-module $\C^{\L(\Omega)}$ associated with any nonzero $\lambda\in \C$. This $\U$-module is denoted by $\C^{\L(\Omega)}(\lambda)$. 
Let $\iota(x_0):\C^{\L(\Omega)}\to \C^{\L(\Omega/x_0)}\otimes \C^{\L(x_0)}$ denote the linear map given by $x\mapsto (x+x_0)/x_0\otimes x\cap x_0$ for all $x\in \L(\Omega)$. 
Using the triple coordinate system yields a $q$-analog of \cite[Lemma 5.4]{Huang:CG&Johnson}: The diagram 
\begin{table}[H]
\centering
\begin{tikzpicture}
\matrix(m)[matrix of math nodes,
row sep=4em, column sep=4em,
text height=1.5ex, text depth=0.25ex]
{
\C^{\L(\Omega)}(q^{\dim x_0})
&\C^{\L(\Omega/ x_0)}(1)\otimes \C^{\L(x_0)}(q^{\dim x_0})\\
\C^{\L(\Omega)}(q^{\dim x_0})
&\C^{\L(\Omega/ x_0)}(1)\otimes \C^{\L(x_0)}(q^{\dim x_0})\\
};
\path[->,font=\scriptsize,>=angle 90]
(m-1-1) edge node[left] {$X$} (m-2-1)
(m-1-1) edge node[above] {$\iota(x_0)$} (m-1-2)
(m-2-1) edge node[below] {$\iota(x_0)$} (m-2-2)
(m-1-2) edge node[right] {$\Delta(X)$} (m-2-2);
\end{tikzpicture}
\end{table}
\noindent is commutative for all $X\in \U$. 
In Section \ref{s:bilinear} we discuss the inherent bilinear forms schemes in the triple coordinate system $(\L(\Omega)_{x_0},\subseteq)$. 
In Section \ref{s:WmoduleL(Omega)} we display a $\W$-module structure $\C^{\L(\Omega)}$ with respect to $x_0$. 
This $\W$-module is denoted by $\C^{\L(\Omega)}(x_0)$. 
To ensure the correctness of the $\W$-module $\C^{\L(\Omega)}(x_0)$, we manipulate the poset $(\L(\Omega)_{x_0},\subseteq)$ to carefully verify the equations from \cite{Watanabe:2017} again and build some new equations throughout Sections \ref{s:L(Omega)}--\ref{s:WmoduleL(Omega)}. 
It is noteworthy that the $\W$-module $\C^{\L(\Omega)}(x_0)$ is not only an extension of the $\U$-module $\C^{\L(\Omega)}(q^{\dim x_0})$ but also a lift of the $\U\otimes \U$-module $\C^{\L(\Omega/x_0)}(1)\otimes \C^{\L(x_0)}(q^{\dim x_0})$. 
In Section \ref{s:dec_WmoduleL(Omega)} the $\W$-module $\C^{\L(\Omega)}(x_0)$ is decomposed into a direct sum of irreducible $\W$-modules. 
Pulling back via $\widetilde{\natural}$ the $\W$-module $\C^{\L(\Omega)}(x_0)$ becomes an $\H$-module. 
Suppose that $k$ is an integer with $0\leq k\leq D$. The space $\C^{\L_k(\Omega)}$ is an $\H$-submodule of $\C^{\L(\Omega)}(x_0)$. This $\H$-module is denoted by $\C^{\L_k(\Omega)}(x_0)$. 
Assume that $x_0\in \L_k(\Omega)$. 
In Section \ref{s:dec_HmoduleL(Omega)} we decompose the $\H$-module $\C^{\L_k(\Omega)}(x_0)$ and study the image $\widetilde{\T}(x_0)$ of the corresponding representation $\H\to {\rm End}(\C^{\L_k(\Omega)})$. 
Suppose that $1\leq k\leq D-1$. 
Let $J_q(D,k)$ denote the Grassmann graph of $k$-dimensional subspaces of $\Omega$.
Let $\T(x_0)$ denote the Terwilliger algebra of $J_q(D,k)$ with respect to $x_0$. 
In Section \ref{s:H&Grassmann} we show that $\T(x_0)$ is a subalgebra of $\widetilde{\T}(x_0)$ and investigate $\T(x_0)$ from this perspective. 
The related research of $\T(x_0)$ can also be found in \cite{Dunkl77,Watanabe2020,
Watanabe:2017,Watanabe2015,TerAlgebraIII}.

\section{An algebraic covering of $\U\otimes \U$}\label{s:W}

\begin{defn}\label{defn:W}
The algebra $\W$ is an algebra over $\C$ defined by generators and relations. 
The generators are $E_1,E_2,F_1,F_2,K_1^{\pm 1},K_2^{\pm 1},I^{\pm 1}$. The relations assert that 
\begin{gather}
\hbox{$I$ is central in $\W$},
\label{W:Icenter}
\\
I I^{-1}=I^{-1}I=1,
\label{W:Ii}
\\
K_1K_1^{-1}=K_1^{-1}K_1=1,
\label{W:K1i}
\\
K_2K_2^{-1}=K_2^{-1}K_2=1,
\label{W:K2i}
\\
[K_1, E_2]=[K_1,F_2]
=[K_1,K_2]
=[K_2,E_1]
=[K_2,F_1]=0,
\label{W:K1E2}
\\
[E_1,K_1]_q=
[K_1,F_1]_q=
[E_2,K_2]_q=
[K_2,F_2]_q=0,
\label{W:E1K1q}
\\
[E_1,E_2]=[E_1,F_2]=[F_1,E_2]=[F_1,F_2]=0,
\label{W:E1E2F1F2}
\\
[E_1, F_1]=\frac{K_1-I K_1^{-1}}{q-q^{-1}},
\label{W:E1F1}
\\
[E_2,F_2]=\frac{I K_2-K_2^{-1}}{q-q^{-1}}.
\label{W:E2F2}
\end{gather}
\end{defn}

For the sake of simplicity the relations (\ref{W:Icenter})--(\ref{W:K2i}) will be used tacitly in the rest of this paper.

\begin{thm}\label{thm:W->U2}
There exists a unique algebra homomorphism $\flat:\W\to \U\otimes \U$ that sends
\begin{alignat*}{2}
E_1 \;\;&\mapsto \;\; E\otimes 1, 
\qquad 
\quad \;\: \,
E_2 &&\;\;\mapsto \;\; 1\otimes E,
\\
F_1 \;\;&\mapsto \;\; F\otimes 1,
\qquad 
\quad \;\: \:
F_2 &&\;\;\mapsto \;\; 1\otimes F,
\\
K_1^{\pm 1} \;\;&\mapsto \;\; K^{\pm 1}\otimes 1,
\qquad 
K_2^{\pm 1} &&\;\;\mapsto \;\; 1\otimes K^{\pm 1},
\\
I^{\pm 1} \;\;&\mapsto \;\; 1\otimes 1.
\end{alignat*}
Moreover $\flat$ is surjective.
\end{thm}
\begin{proof}
By Definition \ref{defn:U} the elements $E\otimes 1,1\otimes E,F\otimes 1,1\otimes F,K^{\pm 1}\otimes 1,1\otimes K^{\pm 1},1\otimes 1$ of $\U\otimes \U$ satisfy the relations (\ref{W:Icenter})--(\ref{W:E2F2}).
Hence the existence of $\flat$ follows. Since the algebra $\W$ is generated by $E_1,E_2,F_1,F_2,K_1^{\pm 1},K_2^{\pm 1},I^{\pm 1}$ the uniqueness of $\flat$ follows. Since the algebra $\U$ is generated by $E,F,K^{\pm 1}$ the algebra homomorphism $\flat$ is surjective.
\end{proof}

In light of Theorem \ref{thm:W->U2} we call the algebra $\W$ an {\it algebraic covering} of $\U\otimes \U$.

\begin{thm}
\label{thm:Wbasis}
The elements 
\begin{gather}\label{Wbasis}
E_1^{i_1} E_2^{i_2}
F_1^{j_1} F_2^{j_2}
K_1^{\ell_1} K_2^{\ell_2} 
I^n
\qquad 
\hbox{for all $i_1,i_2,j_1,j_2\in \N$ and all $\ell_1,\ell_2,n\in \Z$}
\end{gather}
are a basis for $\W$.
\end{thm}
\begin{proof}
Let $\W^{'}$ denote the subspace of $\W$ spanned by (\ref{Wbasis}). Clearly  
\begin{gather}
I^{\pm 1}\W^{'}\subseteq\W^{'},
\label{IW}
\\
E_1\W^{'}\subseteq\W^{'}. 
\label{E1W}
\end{gather}
Using (\ref{W:E1E2F1F2}) yields that 
\begin{gather}\label{E2W}
E_2\W^{'}\subseteq\W^{'}.
\end{gather}
Using (\ref{W:K1E2}) and (\ref{W:E1K1q}) yields that 
\begin{gather}
K_1^{\pm 1}\W^{'}\subseteq\W^{'},
\label{K1W}
\\
K_2^{\pm 1}\W^{'}\subseteq\W^{'}.
\label{K2W}
\end{gather}
We now show that 
\begin{gather}
F_1 \W^{'}\subseteq \W^{'},
\label{F1W}
\\
F_2 \W^{'}\subseteq \W^{'}.
\label{F2W}
\end{gather}
To see (\ref{F1W}) it suffices to show that 
\begin{gather}\label{F1W'}
F_1 E_1^{i_1} E_2^{i_2}
F_1^{j_1} F_2^{j_2}
K_1^{\ell_1} K_2^{\ell_2} 
I^n \in \W^{'}
\end{gather}
for all $i_1,i_2,j_2,j_2\in \N$ and $\ell_1,\ell_2,n\in \Z$. We proceed by induction on $i_1$. If $i_1=0$ then (\ref{F1W'}) holds by (\ref{W:E1E2F1F2}). Suppose that $i_1\geq 1$. Using (\ref{W:E1F1}) yields that the left-hand side of (\ref{F1W'}) is equal to 
\begin{gather}\label{F1W'_1}
E_1 F_1 E_1^{i_1-1} E_2^{i_2}
F_1^{j_1} F_2^{j_2}
K_1^{\ell_1} K_2^{\ell_2} I^{n}
-
\frac{K_1-I K_1^{-1}}{q-q^{-1}}
E_1^{i_1-1} E_2^{i_2}
F_1^{j_1} F_2^{j_2}
K_1^{\ell_1} K_2^{\ell_2} I^{n}.
\end{gather}
By the induction hypothesis and (\ref{E1W}) the first summand of (\ref{F1W'_1}) is in $\W^{'}$. By (\ref{IW}) and (\ref{K1W}) the second summand of (\ref{F1W'_1}) is in $\W^{'}$. Therefore (\ref{F1W}) follows. 
To see (\ref{F2W}), by (\ref{W:E1E2F1F2}) and (\ref{E1W}) it suffices to show that 
\begin{gather}\label{F2W'}
F_2 E_2^{i_2}
F_1^{j_1} F_2^{j_2}
K_1^{\ell_1} K_2^{\ell_2} I^{n}
\in \W^{'}
\end{gather}
for all $i_2,j_1,j_2\in \N$ and $\ell_1,\ell_2,n\in \Z$. We proceed by induction on $i_2$. If $i_2=0$ then (\ref{F2W'}) holds by (\ref{W:E1E2F1F2}). 
Suppose that $i_2\geq 1$. Using (\ref{W:E2F2}) yields that the left-hand side of (\ref{F2W'}) is equal to 
\begin{gather}\label{F2W'_1}
E_2 F_2 E_2^{i_2-1}
F_1^{j_1} F_2^{j_2}
K_1^{\ell_1} K_2^{\ell_2} I^{n}
-
\frac{I K_2-K_2^{-1}}{q-q^{-1}}
E_2^{i_2-1}
F_1^{j_1} F_2^{j_2}
K_1^{\ell_1} K_2^{\ell_2} I^{n}.
\end{gather}
By the induction hypothesis and (\ref{E2W}) the first summand of (\ref{F2W'_1}) is in $\W^{'}$.
By (\ref{IW}) and (\ref{K2W}) the second summand of (\ref{F2W'_1}) is in $\W^{'}$. 
Therefore (\ref{F2W}) follows. 
By (\ref{IW})--(\ref{F2W})  the space $\W^{'}$ is a left ideal of $\W$. Since $1\in \W^{'}$ it follows that $\W^{'}=\W$. Therefore $\W$ is spanned by (\ref{Wbasis}). 

It remains to show the linear independence of (\ref{Wbasis}). Consider the Laurent polynomial ring 
$$
\C[x_1,x_2,y_1,y_2,z_1^{\pm 1},z_2^{\pm 1},w^{\pm 1}]
$$ 
where $x_1,x_2,y_1,y_2,z_1,z_2,w$ are mutually distinct indeterminates over $\C$. Using Definition \ref{defn:W} it is straightforward to verify that  $\C[x_1,x_2,y_1,y_2,z_1^{\pm 1},z_2^{\pm 1},w^{\pm 1}]$ has a $\W$-module structure given by 
\begin{align*}
E_1 x_1^{i_1}x_2^{i_2} y_1^{j_1}y_2^{j_2}z_1^{\ell_1} z_2^{\ell_2} w^{n} 
&=x_1^{i_1+1}x_2^{i_2} y_1^{j_1}y_2^{j_2}z_1^{\ell_1} z_2^{\ell_2} w^{n}, 
\\
E_2 x_1^{i_1}x_2^{i_2} y_1^{j_1}y_2^{j_2}z_1^{\ell_1} z_2^{\ell_2} w^{n} 
&=x_1^{i_1}x_2^{i_2+1} y_1^{j_1}y_2^{j_2}z_1^{\ell_1} z_2^{\ell_2} w^{n}, 
\\
F_1 x_1^{i_1}x_2^{i_2} y_1^{j_1}y_2^{j_2}z_1^{\ell_1} z_2^{\ell_2} w^{n} 
&=
x_1^{i_1}x_2^{i_2} y_1^{j_1+1}y_2^{j_2}z_1^{\ell_1} z_2^{\ell_2} w^{n} 
\\
&
-\left\{
\begin{array}{ll}
[i_1]_q x_1^{i_1-1} x_2^{i_2} 
y_1^{j_1}
y_2^{j_2}
\left(
\frac{q^{i_1-2j_1-1} z_1-q^{2j_1-i_1+1} w z_1^{-1}}{q-q^{-1}}
\right)
z_1^{\ell_1} z_2^{\ell_2} w^{n} 
\qquad &\hbox{if $i_1\geq 1$},
\\
0 \qquad &\hbox{if $i_1=0$},
\end{array}
\right.
\\
F_2 x_1^{i_1}x_2^{i_2} y_1^{j_1}y_2^{j_2}z_1^{\ell_1} z_2^{\ell_2} w^{n} 
&=
x_1^{i_1}x_2^{i_2} y_1^{j_1}y_2^{j_2+1} z_1^{\ell_1} z_2^{\ell_2} w^{n} 
\\
&
-\left\{
\begin{array}{ll}
[i_2]_q x_1^{i_1} x_2^{i_2-1} y_1^{j_1} y_2^{j_2} 
\left(
\frac{q^{i_2-2j_2-1} w z_2-q^{2j_2-i_2+1} z_2^{-1}}{q-q^{-1}} 
\right)
z_1^{\ell_1} z_2^{\ell_2} w^n
\qquad &\hbox{if $i_2\geq 1$},
\\
0 \qquad &\hbox{if $i_2=0$},
\end{array}
\right.
\\
K_1^{\pm 1} x_1^{i_1}x_2^{i_2} y_1^{j_1}y_2^{j_2}z_1^{\ell_1} z_2^{\ell_2} w^{n} 
&=
q^{\pm2(i_1-j_1)} x_1^{i_1}x_2^{i_2} y_1^{j_1}y_2^{j_2}z_1^{\ell_1\pm 1} z_2^{\ell_2} w^{n} ,
\\
K_2^{\pm 1} x_1^{i_1}x_2^{i_2} y_1^{j_1}y_2^{j_2}z_1^{\ell_1} z_2^{\ell_2} w^{n} 
&=
q^{\pm2(i_2-j_2)}  x_1^{i_1}x_2^{i_2} y_1^{j_1}y_2^{j_2}z_1^{\ell_1} z_2^{\ell_2\pm 1} w^{n} ,
\\
I^{\pm 1} x_1^{i_1}x_2^{i_2} y_1^{j_1}y_2^{j_2}z_1^{\ell_1} z_2^{\ell_2} w^{n} 
&=
x_1^{i_1}x_2^{i_2} y_1^{j_1}y_2^{j_2}z_1^{\ell_1} z_2^{\ell_2} w^{n\pm 1}
\end{align*} 
for all $i_1,i_2,j_1,j_2\in \N$ and $\ell_1,\ell_2,n\in \Z$. Observe that $E_1^{i_1} E_2^{i_2}
F_1^{j_1} F_2^{j_2}
K_1^{\ell_1} K_2^{\ell_2} I^{n}
$ maps $1$ to 
\begin{gather}\label{e:Rbasis}
x_1^{i_1} x_2^{i_2}
y_1^{j_1} y_2^{j_2}
z_1^{\ell_1} z_2^{\ell_2} w^{n}
\qquad 
\hbox{for all $i_1,i_2,j_1,j_2\in \N$ and $\ell_1,\ell_2,n\in \Z$}.
\end{gather}
Since the elements (\ref{e:Rbasis}) are a basis for $\C[x_1,x_2,y_1,y_2,z_1^{\pm 1},z_2^{\pm 1},w^{\pm 1}]$ it follows that the elements (\ref{Wbasis}) are linearly independent. The result follows.
\end{proof}

\begin{prop}
The kernel of $\flat$ is equal to the two-sided ideal of $\W$ generated by $I-1$. In particular $\flat:\W\to \U\otimes \U$ is not injective.
\end{prop}
\begin{proof}
Let $\mathcal I$ denote the two-sided ideal of $\W$ generated by $I-1$. By Theorem \ref{thm:W->U2} the two-sided ideal $\mathcal I$ of $\W$ is contained in the kernel of $\flat$. 
Hence there exists a unique algebra homomorphism $\overline{\flat}:\W/\mathcal I\to \U\otimes \U$ given by $\overline{\flat}(X+\mathcal I)=\flat(X)$ for all $X\in \W$. By Theorem \ref{thm:Wbasis} the quotient algebra $\W/\mathcal I$ is spanned by 
\begin{gather}
\label{W/Kbasis}
E_1^{i_1} E_2^{i_2}
F_1^{j_1} F_2^{j_2}
K_1^{\ell_1} K_2^{\ell_2}
+\mathcal I 
\qquad 
\hbox{for all $i_1,i_2,j_1,j_2\in \N$ and all $\ell_1,\ell_2\in \Z$}.
\end{gather}
By Theorem \ref{thm:W->U2} the algebra homomorphism $\overline \flat$ sends (\ref{W/Kbasis}) to 
\begin{gather}
\label{U2basis}
E^{i_1}F^{j_1} K^{\ell_1}
\otimes 
E^{i_2}F^{j_2} K^{\ell_2}
\qquad 
\hbox{for all $i_1,i_2,j_1,j_2\in \N$ and all $\ell_1,\ell_2\in \Z$}.
\end{gather}
By \cite[Theorem 1.5]{jantzen} the elements (\ref{U2basis}) are a basis for $\U\otimes \U$. Therefore $\overline \flat$ is an algebra isomorphism. 
The proposition follows.
\end{proof}

\section{A lift of the algebra homomorphism $\Delta$ to $\W$}\label{s:lift1}

\begin{thm}\label{thm1:U->W}
There exists a unique algebra homomorphism $\widetilde{\Delta}:\U\to \W$ that sends
\begin{eqnarray}
E &\mapsto & E_1+K_1^{-1} E_2,
\label{E'}
\\
F &\mapsto & F_1 K_2+F_2,
\label{F'}
\\
K^{\pm 1} &\mapsto & K_1^{\pm 1} K_2^{\pm 1}.
\label{K'}
\end{eqnarray}
\end{thm}
\begin{proof}
Let $\widetilde{\Delta}(E),\widetilde{\Delta}(F),\widetilde{\Delta}(K^{\pm 1})$ denote the right-hand sides of (\ref{E'})--(\ref{K'}) respectively.  To see the existence of $\widetilde{\Delta}$, by Definition \ref{defn:U} it suffices to show that 
\begin{gather}
\label{wideDelta-1}
\widetilde{\Delta}(K)\widetilde{\Delta}(K^{-1})=\widetilde{\Delta}(K^{-1})\widetilde{\Delta}(K)=1,
\\
[\widetilde{\Delta}(E),\widetilde{\Delta}(K)]_q=0,
\label{wideDelta-2}
\\
[\widetilde{\Delta}(K),\widetilde{\Delta}(F)]_q=0,
\label{wideDelta-3}
\\
[\widetilde{\Delta}(E),\widetilde{\Delta}(F)]=\frac{\widetilde{\Delta}(K)-\widetilde{\Delta}(K^{-1})}{q-q^{-1}}.
\label{wideDelta-4}
\end{gather} 
Clearly (\ref{wideDelta-1}) holds. 
Using (\ref{W:K1E2}) and (\ref{W:E1K1q}) yields that 
\begin{align}
[E_1,K_1 K_2]_q&=
[K_1^{-1} E_2,K_1 K_2]_q =0,
\label{wideDelta-2'}
\\
[K_1K_2,F_1 K_2]_q&=
[K_1K_2,F_2]_q =0.
\label{wideDelta-3'}
\end{align}
The relations (\ref{wideDelta-2}) and (\ref{wideDelta-3}) are immediate from (\ref{wideDelta-2'}) and (\ref{wideDelta-3'}) respectively. 
Using (\ref{W:K1E2})--(\ref{W:E2F2}) yields that 
\begin{align*}
[E_1,F_1 K_2]&=\frac{K_1K_2-I K_1^{-1} K_2}{q-q^{-1}},
\\
[K_1^{-1} E_2, F_2]&=\frac{I K_1^{-1} K_2-K_1^{-1} K_2^{-1}}{q-q^{-1}},
\\
[K_1^{-1} E_2, F_1 K_2]&=[E_1,F_2]=0.
\end{align*}
The relation (\ref{wideDelta-4}) is immediate from the above equations. 
The existence of $\widetilde{\Delta}$ follows. 
Since the algebra $\U$ is generated by $E,F,K^{\pm 1}$ the uniqueness follows. 
\end{proof}

Recall the comultiplication $\Delta$ of $\U$ from Lemma \ref{lem:U_comultiplication}. 

\begin{thm}\label{thm2:U->W}
The following diagram commutes:
\begin{table}[H]
\centering
\begin{tikzpicture}
\matrix(m)[matrix of math nodes,
row sep=4em, column sep=4em,
text height=1.5ex, text depth=0.25ex]
{
\U
&\W\\
&\U\otimes \U\\
};
\path[->,font=\scriptsize,>=angle 90]
(m-1-1) edge node[above] {$\widetilde{\Delta}$} (m-1-2)
(m-1-2) edge node[right] {$\flat$} (m-2-2)
(m-1-1) edge[bend right] node[below] {$\Delta$} (m-2-2);
\end{tikzpicture}
\end{table}
\end{thm}
\begin{proof}
By Theorems \ref{thm:W->U2} and \ref{thm1:U->W} the algebra homomorphisms $\Delta$ and $\flat \circ \widetilde{\Delta}$ agree at $E,F,K^{\pm 1}$. 
Since the algebra $\U$ is generated by $E,F,K^{\pm 1}$ the result follows.
\end{proof}

By Theorem \ref{thm2:U->W} the algebra homomorphism $\widetilde{\Delta}$ is a lift of the algebra homomorphism $\Delta$ to $\W$ across $\flat$.

\section{Two central elements $\Lambda_1$ and $\Lambda_2$ of $\W$}\label{s:Wcentral}

Let $\Lambda_1$ and $\Lambda_2$ denote the elements of $\W$ given by  
\begin{align}
\Lambda_1 &=
(q-q^{-1})^2 E_1 F_1+q^{-1} K_1+q I K_1^{-1},
\label{Lambda1}
\\
\Lambda_2 &=
(q-q^{-1})^2 E_2 F_2+q^{-1} I K_2+q K_2^{-1}. 
\label{Lambda2}
\end{align}
Recall the Casimir element $\Lambda$ of $\U$ from (\ref{Lambda}).

\begin{lem}\label{lem:Lambda12&flat}
The following equations hold in $\U\otimes \U$:
\begin{enumerate}
\item $\flat(\Lambda_1)=\Lambda\otimes 1$.

\item $\flat(\Lambda_2)=1\otimes \Lambda$.
\end{enumerate}
\end{lem}
\begin{proof}
The lemma follows by employing Theorem \ref{thm:W->U2} to  evaluate the equations obtained by applying $\flat$ to both sides of (\ref{Lambda1}) and  (\ref{Lambda2}).
\end{proof}

\begin{lem}\label{lem:Lambda12}
The following equations hold in $\W$:
\begin{enumerate}
\item $\Lambda_1=(q-q^{-1})^2 F_1 E_1+q K_1+q^{-1} I K_1^{-1}$.

\item $\Lambda_2=(q-q^{-1})^2 F_2 E_2+q I K_2+q^{-1} K_2^{-1}$.
\end{enumerate}
\end{lem}
\begin{proof}
(i): 
By (\ref{W:E1F1}) we may replace $E_1F_1$ by $F_1 E_1+\frac{K_1-I K_1^{-1}}{q-q^{-1}}$ in (\ref{Lambda1}).

(ii): By (\ref{W:E2F2}) we may replace $E_2 F_2$ by $F_2 E_2+\frac{I K_2-K_2^{-1}}{q-q^{-1}}$ in (\ref{Lambda2}).
\end{proof}

\begin{lem}\label{lem:Lambda12_central}
\begin{enumerate}
\item The element $\Lambda_1$ is central in $\W$.

\item The element $\Lambda_2$ is central in $\W$.
\end{enumerate}
\end{lem}
\begin{proof}
(i): Using (\ref{W:K1E2})--(\ref{W:E1E2F1F2}) yields that 
\begin{align*}
[E_1F_1,E_2]&=
[K_1^{\pm 1},E_2]=0,
\\
[E_1F_1,F_2]&=
[K_1^{\pm 1},F_2]=0,
\\
[E_1F_1,K_1^{\pm 1}]&=[E_1F_1,K_2^{\pm 1}]=0.
\end{align*}
By the above equations the element (\ref{Lambda1})  commutes with $E_2,F_2,K_1^{\pm 1},K_2^{\pm 1}$.

Using (\ref{W:E1K1q}) and (\ref{W:E1F1})  yields that 
\begin{align*}
[E_1F_1,E_1]&=E_1\frac{I K_1^{-1}-K_1}{q-q^{-1}},
\\
[K_1,E_1]&=q(q-q^{-1}) E_1 K_1,
\\
[K_1^{-1},E_1]&=-q^{-1}(q-q^{-1}) E_1 K_1^{-1}.
\end{align*}
By the above equations the element (\ref{Lambda1})  commutes with $E_1$. Using (\ref{W:E1K1q}) and (\ref{W:E1F1}) yields that 
\begin{align*}
[E_1 F_1,F_1]&=\frac{K_1-I K_1^{-1}}{q-q^{-1}} F_1,
\\
[K_1,F_1]&=-q(q-q^{-1})  K_1 F_1,
\\
[K_1^{-1},F_1]&=q^{-1}(q-q^{-1})  K_1^{-1} F_1 .
\end{align*}
By the above equations the element (\ref{Lambda1})  commutes with $F_1$. Since the algebra $\W$ is generated by $E_1,E_2,F_1,F_2,K_1^{\pm 1},K_2^{\pm 1},I^{\pm 1}$ the statement (i) follows.

(ii): Using (\ref{W:K1E2})--(\ref{W:E1E2F1F2}) yields that 
\begin{align*}
[E_2F_2,E_1]&=
[K_2^{\pm 1},E_1]=0,
\\
[E_2F_2,F_1]&=
[K_2^{\pm 1},F_1]=0,
\\
[E_2F_2,K_1^{\pm 1}]&=[E_2F_2,K_2^{\pm 1}]=0.
\end{align*}
By the above equations the element (\ref{Lambda2}) commutes with $E_1,F_1,K_1^{\pm 1},K_2^{\pm 1}$.

Using (\ref{W:E1K1q}) and (\ref{W:E2F2}) yields that 
\begin{align*}
[E_2F_2,E_2]&=E_2\frac{K_2^{-1}-I K_2}{q-q^{-1}},
\\
[K_2,E_2]&=q(q-q^{-1}) E_2 K_2,
\\
[K_2^{-1},E_2]&=-q^{-1}(q-q^{-1}) E_2 K_2^{-1}.
\end{align*}
By the above equations the element (\ref{Lambda2})  commutes with $E_2$. Using (\ref{W:E1K1q}) and (\ref{W:E2F2}) yields that 
\begin{align*}
[E_2F_2,F_2]&=\frac{I K_2-K_2^{-1}}{q-q^{-1}} F_2,
\\
[K_2,F_2]&=-q(q-q^{-1}) K_2 F_2,
\\
[K_2^{-1},F_2]&=q^{-1}(q-q^{-1}) K_2^{-1} F_2.
\end{align*}
By the above equations the element (\ref{Lambda2})  commutes with $F_2$. Since the algebra $\W$ is generated by $E_1,E_2,F_1,F_2,K_1^{\pm 1},K_2^{\pm 1},I^{\pm 1}$ the statement (ii) follows.
\end{proof}

\section{A lift of the algebra homomorphism $\natural$ to $\W$}\label{s:lift2}

Recall the expression (\ref{Lambda}) for the Casimir element $\Lambda$ of $\U$. Recall the algebra homomorphism $\widetilde{\Delta}$ from Theorem \ref{thm1:U->W}.

\begin{lem}
\label{lem:wideLambda}
The image of $\Lambda$ under $\widetilde{\Delta}$ is equal to 
\begin{gather}\label{wideLambda1}
(q-q^{-1})^2 (K_1^{-1}E_2 F_1 K_2+ E_1F_2)+\Lambda_1 K_2+K_1^{-1} \Lambda_2 -(q+q^{-1})
I K_1^{-1} K_2.
\end{gather}
\end{lem}
\begin{proof}
Applying $\widetilde{\Delta}$ to either side of (\ref{Lambda}) we use Theorem \ref{thm1:U->W} to expand the resulting equation. It follows that $\widetilde{\Delta}(\Lambda)$ is equal to 
\begin{gather}\label{wideLambda2}
(q-q^{-1})^2(K_1^{-1} E_2 F_1 K_2+E_1F_2+E_1 F_1 K_2+K_1^{-1} E_2 F_2)+q^{-1} K_1 K_2 +q K_1^{-1} K_2^{-1}.
\end{gather}
The subtraction of (\ref{wideLambda2}) from (\ref{wideLambda1}) is equal to 
\begin{align*}
&(\Lambda_1-(q-q^{-1})^2 E_1 F_1-q^{-1} K_1-q I K_1^{-1}) K_2
\\
&\quad +\;
K_1^{-1} (\Lambda_2-(q-q^{-1})^2 E_2 F_2-q^{-1} I K_2-q K_2^{-1}).
\end{align*}
By (\ref{Lambda1}) and (\ref{Lambda2}) the above element is equal to zero.  The lemma follows. 
\end{proof}

\begin{thm}
\label{thm1:H->W}
There exists a unique algebra homomorphism $\widetilde{\natural}: \H\to \W$ that sends
\begin{eqnarray}
A &\mapsto & K_2^{-1},
\label{A:widenatural}
\\
B &\mapsto & 
\widetilde{\Delta}(\Lambda),
\label{B:widenatural}
\\
C &\mapsto & 
I K_1^{-1}
-q(q-q^{-1})^2 E_1 F_2 K_2^{-1}.
\label{C:widenatural}
\end{eqnarray}
The algebra homomorphism $\widetilde{\natural}$ sends 
\begin{eqnarray}
\alpha &\mapsto & \Lambda_2+\Lambda_1 
\widetilde{\Delta}(K^{-1}),
\label{alpha:widenatural}
\\
\beta &\mapsto &  
I \widetilde{\Delta}(K^{-1}),
\label{beta:widenatural}
\\
\gamma  &\mapsto & \Lambda_1+\Lambda_2 
\widetilde{\Delta}(K^{-1}).
\label{gamma:widenatural}
\end{eqnarray}
\end{thm}
\begin{proof}
Let $\widetilde{\natural}(A)$, $\widetilde{\natural}(B)$, $\widetilde{\natural}(C)$, $\widetilde{\natural}(\alpha)$, $\widetilde{\natural}(\beta)$, $\widetilde{\natural}(\gamma)$ denote the right-hand sides of (\ref{A:widenatural})--(\ref{gamma:widenatural}) respectively. Since $\Lambda$ is central in $\U$ the element $\widetilde{\Delta}(\Lambda)$ commutes with $\widetilde{\Delta}(K^{-1})$. Recall from (\ref{K'}) that $\widetilde{\Delta}(K^{-1})=K_1^{-1}K_2^{-1}$. 
Using (\ref{W:K1E2}) and (\ref{W:E1K1q}) yields that $E_1F_2$ commutes with $\widetilde{\Delta}(K^{-1})$. 
Combined with Lemma \ref{lem:Lambda12_central} each of $\widetilde{\natural}(\alpha)$, $\widetilde{\natural}(\beta)$, $\widetilde{\natural}(\gamma)$ commutes with $\widetilde{\natural}(A)$, $\widetilde{\natural}(B)$, $\widetilde{\natural}(C)$.
To see the existence of $\widetilde{\natural}$, by Definition \ref{defn:H} it remains to show that 
\begin{align}
\widetilde{\natural}(\alpha)
&=
\frac{[\widetilde{\natural}(B),\widetilde{\natural}(C)]_q}{q-q^{-1}}
+
(q+q^{-1})\widetilde{\natural}(A),
\label{widenatural-1}
\\
\widetilde{\natural}(\beta)
&=
\frac{[\widetilde{\natural}(C),\widetilde{\natural}(A)]_q}{q-q^{-1}},
\label{widenatural-2}
\\
\widetilde{\natural}(\gamma)
&=\frac{[\widetilde{\natural}(A),\widetilde{\natural}(B)]_q}{q-q^{-1}}
+(q+q^{-1})\widetilde{\natural}(C).
\label{widenatural-3}
\end{align}

Using (\ref{W:K1E2}) and (\ref{W:E1K1q}) yields that 
\begin{align}
[K_1^{-1} E_2 F_1 K_2, K_1^{-1}]_q &=0,
\label{widenatural-1-1}
\\
[E_1 F_2, K_1^{-1}]_q &=q^{-1}(q^2-q^{-2}) E_1 F_2  K_1^{-1},
\label{widenatural-1-2}
\\
[E_1F_2,E_1 F_2 K_2^{-1}]_q&=0,
\label{widenatural-1-3}
\\
[K_2, E_1 F_2 K_2^{-1}]_q &=0,
\label{widenatural-1-4}
\\
[K_1^{-1}, E_1 F_2 K_2^{-1}]_q&=0,
\label{widenatural-1-5}
\\
[K_1^{-1} K_2, E_1 F_2 K_2^{-1}]_q &=-q^{-2}(q-q^{-1}) E_1 F_2 K_1^{-1}.
\label{widenatural-1-6}
\end{align}
By (\ref{W:K1E2}) and (\ref{W:E1K1q}) the product of $K_1^{-1} E_2 F_1 K_2$ and $E_1 F_2 K_2^{-1}$ is equal to 
$$
q^{-2} K_1^{-1} E_2 F_1 E_1 F_2.
$$
By (\ref{Lambda2}) and Lemma \ref{lem:Lambda12}(i) along with Lemma \ref{lem:Lambda12_central}(i), the above element is equal to 
$$
\frac{K_1^{-1}(\Lambda_1-q K_1-q^{-1} I K_1^{-1})(\Lambda_2-q^{-1} I K_2-q K_2^{-1})}{q^2(q-q^{-1})^4}.
$$
By (\ref{W:K1E2}) and (\ref{W:E1K1q}) the product of $E_1 F_2 K_2^{-1}$ and $K_1^{-1} E_2 F_1 K_2$ is equal to 
$$
K_1^{-1} E_1 F_2 E_2 F_1.
$$
By (\ref{Lambda1}) and Lemma \ref{lem:Lambda12}(ii) along with Lemma \ref{lem:Lambda12_central}(ii), the above element is equal to 
$$
\frac{K_1^{-1}(\Lambda_1-q^{-1} K_1-q I K_1^{-1})(\Lambda_2-q I K_2-q^{-1} K_2^{-1})}{(q-q^{-1})^4}.
$$
Hence $[K_1^{-1} E_2 F_1 K_2, E_1 F_2 K_2^{-1}]_q$ is equal to 
\begin{align}
\label{widenatural-1-7}
\frac{\Lambda_1 K_1^{-1} (I K_2-K_2^{-1})+\Lambda_2(I K_1^{-2}-1)+(q+q^{-1})(K_2^{-1}-I^2 K_1^{-2} K_2)}{q(q-q^{-1})^3}.
\end{align}
Recall an expression for $\widetilde{\Delta}(\Lambda)$ from Lemma \ref{lem:wideLambda}. 
Applying Lemma \ref{lem:Lambda12_central} and (\ref{widenatural-1-1})--(\ref{widenatural-1-7}) a routine calculation shows that 
$$
[\widetilde{\natural}(B),\widetilde{\natural}(C)]_q
=(q-q^{-1})(\Lambda_1 K_1^{-1} K_2^{-1}+\Lambda_2)-(q^2-q^{-2}) K_2^{-1}.
$$
Hence the relation (\ref{widenatural-1}) holds.

Using (\ref{W:K1E2}) and (\ref{W:E1K1q}) yields that 
$[E_1 F_2 K_2^{-1}, K_2^{-1}]_q=0$. 
The relation (\ref{widenatural-2}) follows by applying the above equation to evaluate $[\widetilde{\natural}(C),\widetilde{\natural}(A)]_q$. 
Using (\ref{W:K1E2}) and (\ref{W:E1K1q}) yields that  
\begin{align*}
[K_2^{-1}, K_1^{-1} E_2 F_1 K_2]_q&=0,
\\
[K_2^{-1}, E_1 F_2]_q &=q(q^2-q^{-2}) E_1 F_2 K_2^{-1}.
\end{align*}
The relation (\ref{widenatural-3}) follows by applying the above two equations and Lemma \ref{lem:Lambda12_central} to evaluate $[\widetilde{\natural}(A),\widetilde{\natural}(B)]_q$. 
The existence of $\widetilde{\natural}$ follows. 
Since the algebra $\H$ is generated by $A,B,C$ the uniqueness of $\widetilde{\natural}$ follows. 
\end{proof}

Recall the algebra homomorphism $\natural:\H\to \U\otimes \U$ from Theorem \ref{thm:H->U2}.

\begin{thm}
\label{thm2:H->W}
The following diagram commutes:
\begin{table}[H]
\centering
\begin{tikzpicture}
\matrix(m)[matrix of math nodes,
row sep=4em, column sep=4em,
text height=1.5ex, text depth=0.25ex]
{
\H
&\W\\
&\U\otimes \U\\
};
\path[->,font=\scriptsize,>=angle 90]
(m-1-1) edge node[above] {$\widetilde{\natural}$} (m-1-2)
(m-1-2) edge node[right] {$\flat$} (m-2-2)
(m-1-1) edge[bend right] node[below] {$\natural$} (m-2-2);
\end{tikzpicture}
\end{table}
\end{thm}
\begin{proof}
Recall the images of $A,B,C$ under $\widetilde{\natural}$ from Theorem \ref{thm1:H->W}.
By Theorem \ref{thm:W->U2} the images of $A$ and $C$ under $\flat\circ \widetilde{\natural}$ are equal to the right-hand sides of (\ref{A:natural}) and (\ref{C:natural}) respectively. By Theorem \ref{thm2:U->W} the image of $B$ under $\flat\circ \widetilde{\natural}$ is equal to the right-hand side of (\ref{B:natural}). Since the algebra $\H$ is generated by $A,B,C$ the result follows.
\end{proof}

By Theorem \ref{thm2:H->W} the algebra homomorphism $\widetilde{\natural}$ is a lift of the algebra homomorphism $\natural$  to $\W$ across $\flat$.

\section{Finite-dimensional irreducible $\W$-modules}\label{s:Wmodule}

In this section we lay the groundwork for finite-dimensional irreducible $\W$-modules. Since $K_1,K_2,I$ are invertible in $\W$ any eigenvalues of $K_1,K_2,I$ on a $\W$-module are nonzero.

\begin{defn}
\label{defn:weight}
Assume that $V$ is a $\W$-module. Let $\lambda_1,\lambda_2,\mu$ denote three nonzero complex numbers. 
We define the following terms:
\begin{enumerate}
\item A vector $v\in V$ is called a {\it weight vector with weight $(\lambda_1,\lambda_2,\mu)$} if $v\not=0$ and 
\begin{gather*}
K_1 v=\lambda_1 v,
\qquad 
K_2 v=\lambda_2 v,
\qquad 
I v=\mu v.
\end{gather*}

\item A weight vector $v$ of $V$ is said to be {\it highest} if 
\begin{gather*}
E_1 v=0,
\qquad 
E_2 v=0.
\end{gather*}
\end{enumerate}
\end{defn}

\begin{defn}\label{defn:VX(theta)}
Assume that $V$ is a $\W$-module. Let $X\in \W$ and $\theta\in \C$ be given. For any integer $n\geq 1$ we let 
$$
V_X^{(n)}(\theta)=\{v\in V\,|\, (X-\theta)^n v=0\}.
$$
For convenience we simply write $V_X(\theta)=V_X^{(1)}(\theta)$.
\end{defn}

\begin{lem}
\label{lem:Wmodule_highest}
Every nonzero finite-dimensional $\W$-module contains a highest weight vector.
\end{lem}
\begin{proof}
Let $V$ denote a nonzero finite-dimensional $\W$-module. 
Since $\C$ is algebraically closed there exists a nonzero scalar $\theta\in \C$ such that $V_{K_1}(\theta)\not=\{0\}$.
Since $q$ is not a root of unity the numbers $\{q^{2i}\theta\}_{i\in \N}$ are mutually distinct. Hence there is a nonzero scalar $\lambda_1=q^{2i}\theta$ for some $i\in \N$ such that $V_{K_1}(\lambda_1)\not=\{0\}$ and 
\begin{gather}\label{V1q2}
V_{K_1}(q^{2}\lambda_1)=\{0\}.
\end{gather}
Since $[K_1,K_2]=0$ the nonzero space $V_{K_1}(\lambda_1)$ is $K_2$-invariant. 
By a similar argument there exists a nonzero $\lambda_2\in \C$ such that  
$V_{K_1}(\lambda_1)\cap V_{K_2}(\lambda_2)\not=\{0\}$ and 
\begin{gather}\label{V12q2}
V_{K_1}(\lambda_1)\cap V_{K_2}(q^2\lambda_2)=\{0\}. 
\end{gather} 
Since $I$ commutes with $K_1$ and $K_2$ there is a nonzero $\mu\in \C$ such that $V_{K_1}(\lambda_1)\cap V_{K_2}(\lambda_2)\cap V_I(\mu)\not=\{0\}$.

Pick any nonzero $v\in V_{K_1}(\lambda_1)\cap V_{K_2}(\lambda_2)\cap V_I(\mu)$. Then $v$ is a weight vector of $V$ with weight $(\lambda_1,\lambda_2,\mu)$ by Definition \ref{defn:weight}(i).
Using (\ref{W:E1K1q}) yields that $E_1v\in V_{K_1}(q^2\lambda_1)$ and $E_2v\in V_{K_2}(q^2\lambda_2)$. It follows from (\ref{V1q2}) that $E_1v=0$. 
Using (\ref{W:K1E2}) yields that $E_2 v\in V_{K_1}(\lambda_1)$.
It follows from (\ref{V12q2}) that $E_2v=0$.
Thus $v$ is a highest weight vector of $V$ by Definition \ref{defn:weight}(ii). 
The lemma follows.
\end{proof}

We are now going to construct finite-dimensional $\W$-modules explicitly.

\begin{lem}
\label{lem:Ln}
For any $n\in \N$ there exists an $(n+1)$-dimensional $\U$-module $L_n$ satisfying the following conditions: 
\begin{enumerate}
\item There exists a basis $\{v_i^{(n)}\}_{i=0}^n$ such that 
\begin{align*}
E v_i^{(n)} &=[i]_q [n-i+1]_q v_{i-1}^{(n)} 
\qquad (1\leq i\leq n),
\qquad 
E v_0^{(n)}=0,
\\
F v_i^{(n)} &= v_{i+1}^{(n)} 
\qquad (0\leq i\leq n-1),
\qquad 
Fv_n^{(n)}=0,
\\
K v_i^{(n)}&=q^{n-2i} v_i^{(n)}
\qquad (0\leq i\leq n).
\end{align*}

\item  The Casimir element $\Lambda$ of $\U$ acts on $L_n$ as scalar multiplication by 
$q^{n+1}+q^{-n-1}$.
\end{enumerate}
\end{lem}
\begin{proof}
For example see \cite[Theorem 2.6]{jantzen}.
\end{proof}

Note that the $\U$-module $L_n$ ($n\in \N$) is irreducible. Recall the algebra homomorphism $\flat:\W\to \U\otimes \U$ from Theorem \ref{thm:W->U2}.
Every $\U\otimes \U$-module can be considered as a $\W$-module by pulling back via $\flat$.  Thus we may view $L_m\otimes L_n$ ($m,n\in \N$) as a $\W$-module. For the sake of convenience, we give a description of the $\W$-module $L_m\otimes L_n$:

\begin{prop}
\label{prop:Lmn}
Suppose that $m,n\in \N$. Then the actions of $E_1,E_2,F_1,F_2,K_1,K_2,I$ on the $\W$-module $L_m\otimes L_n$ are as follows:
\begin{align*}
E_1 v_{i}^{(m)}\otimes  v_{j}^{(n)} &= 
\left\{
\begin{array}{ll}
0 \qquad &\hbox{if $i=0$},
\\
\hbox{$[i]_q[m-i+1]_q v_{i-1}^{(m)}\otimes  v_{j}^{(n)}$}  \qquad &\hbox{if $i\not=0$},
\end{array}
\right.
\\
E_2 v_{i}^{(m)}\otimes  v_{j}^{(n)} &=
\left\{
\begin{array}{ll}
0 \qquad &\hbox{if $j=0$},
\\
\hbox{$[j]_q[n-j+1]_q v_{i}^{(m)}\otimes  v_{j-1}^{(n)}$}  \qquad &\hbox{if $j\not=0$},
\end{array}
\right.
\\
F_1 v_{i}^{(m)}\otimes  v_{j}^{(n)} &= 
\left\{
\begin{array}{ll}
v_{i+1}^{(m)}\otimes  v_{j}^{(n)}  \qquad &\hbox{if $i\not=m$},
\\
0  \qquad &\hbox{if $i=m$},
\end{array}
\right.
\\
F_2 v_{i}^{(m)}\otimes  v_{j}^{(n)} &=
\left\{
\begin{array}{ll}
v_{i}^{(m)}\otimes  v_{j+1}^{(n)} \qquad &\hbox{if $j\not=n$},
\\
0  \qquad &\hbox{if $j=n$},
\end{array}
\right.
\\
K_1 v_{i}^{(m)}\otimes  v_{j}^{(n)} &= q^{m-2i} v_{i}^{(m)}\otimes  v_{j}^{(n)}, 
\\
K_2 v_{i}^{(m)}\otimes  v_{j}^{(n)} &=q^{n-2j} v_{i}^{(m)}\otimes  v_{j}^{(n)}, 
\\
I v_{i}^{(m)}\otimes  v_{j}^{(n)} &=v_{i}^{(m)}\otimes  v_{j}^{(n)}
\end{align*}
for all integers $i$ and $j$ with $0\leq i\leq m$ and $0\leq j\leq n$.
\end{prop}
\begin{proof}
Immediate from Theorem \ref{thm:W->U2} and Lemma \ref{lem:Ln}.
\end{proof}

\begin{lem}
\label{lem:Lmn_weight}
Suppose that $m,n\in \N$. For any nonzero $v\in L_m\otimes L_n$ the following conditions are equivalent:
\begin{enumerate}
\item $v$ is a weight vector of the $\W$-module $L_m\otimes L_n$.

\item $v$ is a scalar multiple of $v_i^{(m)}\otimes v_j^{(n)}$ for some integers $i$ and $j$ with $0\leq i\leq m$ and $0\leq j\leq n$.
\end{enumerate}
\end{lem}
\begin{proof}
(ii) $\Rightarrow$ (i): Immediate from Definition \ref{defn:weight}(i) and Proposition \ref{prop:Lmn}. 

(i) $\Rightarrow$ (ii): By Definition \ref{defn:weight}(i) the vector $v$ is a simultaneous eigenvector of $K_1$ and $K_2$. Since $q$ is not a root of unity and by  Proposition \ref{prop:Lmn}, the statement (ii) follows.
\end{proof}

\begin{lem}
\label{lem:Lmn_highest}
Suppose that $m,n\in \N$. For any nonzero $v\in L_m\otimes L_n$ the following conditions are equivalent:
\begin{enumerate}
\item $v$ is a highest weight vector of the $\W$-module $L_m\otimes L_n$.

\item $v$ is a scalar multiple of $v_0^{(m)}\otimes v_0^{(n)}$.
\end{enumerate}
\end{lem}
\begin{proof}
(ii) $\Rightarrow$ (i): Immediate from Definition \ref{defn:weight} and Proposition \ref{prop:Lmn}. 

(i) $\Rightarrow$ (ii): 
By Lemma \ref{lem:Lmn_weight} there are integers $i$ and $j$ with $0\leq i\leq m$ and $0\leq j\leq n$ such that $v$ is a scalar multiple of $v_i^{(m)}\otimes v_j^{(n)}$.
By Definition \ref{defn:weight}(ii) and Proposition \ref{prop:Lmn} this implies that $i=j=0$. 
\end{proof}

\begin{lem}
\label{lem:Lmn_irr}
Suppose that $m,n\in \N$. Then the following statements hold:
\begin{enumerate}
\item The $\W$-module $L_m\otimes L_n$ is generated by $v_0^{(m)}\otimes v_0^{(n)}$.

\item The $\W$-module $L_m\otimes L_n$ is irreducible.
\end{enumerate}
\end{lem}
\begin{proof}
(i): Using Proposition \ref{prop:Lmn} yields that 
$$
v_i^{(m)}\otimes v_j^{(n)}=F_1^i F_2^j\, v_0^{(m)}\otimes v_0^{(n)}
$$
for all integers $i$ and $j$ with $0\leq i\leq m$ and $0\leq j\leq n$.

(ii): Let $V$ denote a nonzero $\W$-submodule of $L_m\otimes L_n$. It follows from Lemma \ref{lem:Wmodule_highest} that $V$ contains a highest weight vector $v$. Since $v$ is also a highest weight vector of $L_m\otimes L_n$ it follows from Lemma \ref{lem:Lmn_highest} that $v$ is a scalar multiple of $v_0^{(m)}\otimes v_0^{(n)}$. By Lemma \ref{lem:Lmn_irr}(i) it follows that $V= L_m\otimes L_n$. The statement (ii) follows.
\end{proof}

\begin{lem}\label{lem:delta&lambda}
Suppose that $\delta,\lambda \in \C$ with $\delta\in\{\pm 1\}$ and $\lambda\not=0$. Then there exists a unique algebra homomorphism $\sigma_{\delta,\lambda}:\W\to \W$ that sends 
\begin{alignat*}{2}
&E_1 \;\; \mapsto\;\; \lambda E_1, 
\qquad 
&& E_2 \;\; \mapsto\;\; \delta \lambda E_2,
\\
&F_1 \;\; \mapsto\;\; F_1,
\qquad 
&&F_2 \;\; \mapsto\;\; F_2,
\\
&K_1^{\pm 1} \;\; \mapsto\;\;  \lambda^{\pm 1} K_1^{\pm 1},
\qquad 
&&K_2^{\pm 1} \;\; \mapsto\;\; \delta  \lambda^{\mp 1} K_2^{\pm 1},
\\
& I^{\pm 1}\;\; \mapsto\;\; \lambda^{\pm 2} I^{\pm 1}.
\end{alignat*}
Moreover $\sigma_{\delta,\lambda}$ is an algebra automorphism with the inverse $\sigma_{\delta,\lambda^{-1}}$.
\end{lem}
\begin{proof}
It is routine  to verify the existence and uniqueness of $\sigma_{\delta,\lambda}$ by using Definition \ref{defn:W}.  
\end{proof}

Let $V$ denote a $\W$-module. For any $\delta,\lambda\in \C$ with $\delta\in\{\pm 1\}$ and $\lambda\not=0$, we adopt the notation 
$
V^{\delta,\lambda}
$ 
to stand for the $\W$-module obtained by twisting the $\W$-module $V$ via $\sigma_{\delta,\lambda}$. 
For notational convenience let 
$
\I(\W)
$ 
denote the set consisting of all quadruples $(m,n,\delta,\lambda)$ where $m,n\in \N$  and $\delta,\lambda\in \C$ with $\delta\in\{\pm 1\}$ and $\lambda\not=0$.

\begin{lem}
\label{lem1:irrWmodule}
For any $(m,n,\delta,\lambda)\in \I(\W)$ every highest weight vector of the $\W$-module $(L_m\otimes L_n)^{\delta,\lambda}$ has the weight
$
(\lambda q^m,\delta \lambda^{-1} q^n,\lambda^2)$. 
\end{lem}
\begin{proof}
Immediate from Proposition \ref{prop:Lmn} and Lemmas \ref{lem:Lmn_highest} and \ref{lem:delta&lambda}.
\end{proof}

\begin{lem}
\label{lem2:irrWmodule}
For any $(m,n,\delta,\lambda),(m',n',\delta',\lambda')\in \I(\W)$ the following conditions are equivalent:
\begin{enumerate}
\item $(\lambda q^m,\delta \lambda^{-1} q^n,\lambda^2)=(\lambda' q^{m'},\delta' \lambda'^{-1} q^{n'},\lambda'^2)$.

\item $(m,n,\delta,\lambda)=(m',n',\delta',\lambda')$.

\item The $\W$-module $(L_m\otimes L_n)^{\delta,\lambda}$ is isomorphic to $(L_{m'}\otimes L_{n'})^{\delta',\lambda'}$.

\end{enumerate}
\end{lem}
\begin{proof}
(i) $\Rightarrow$ (ii): Let $(\lambda_1,\lambda_2,\mu)$ denote the triple given in (i). Then 
\begin{align*}
q^{2m} &= \lambda_1^2\mu^{-1}=q^{2m'},
\\
q^{2n} &= \lambda_2^2\mu=q^{2n'}.
\end{align*}
Since $q$ is not a root of unity the pair $(m,n)=(m',n')$. Since $m=m'$ and $\lambda q^m=\lambda' q^{m'}$ it follows that $\lambda=\lambda'$.  Since $(n,\lambda)=(n',\lambda')$ and $\delta \lambda^{-1} q^n=\delta' \lambda'^{-1} q^{n'}$ it follows that $\delta=\delta'$. The condition (ii) holds.

(ii) $\Rightarrow$ (iii): Trivial.

(iii)  $\Rightarrow$ (i): Immediate from Lemma \ref{lem1:irrWmodule}.
\end{proof}

We are in a position to prove the main result of this section.

\begin{thm}
\label{thm:irrWmodule} 
Suppose that $V$ is a finite-dimensional $\W$-module and $v$ is a highest weight vector of $V$. Then there exists a unique quadruple $(m,n,\delta,\lambda)\in \I(\W)$ such that 
the $\W$-module $(L_m\otimes L_n)^{\delta,\lambda}$ is isomorphic to the $\W$-submodule of $V$ generated by $v$. 
\end{thm}
\begin{proof}
Let $(\lambda_1,\lambda_2,\mu)$ denote the weight of $v$ in the $\W$-module $V$.
Set 
\begin{align}\label{vij}
v_{ij}=F_1^i F_2^j v
\qquad 
\hbox{for all $i,j\in \N$}.
\end{align}
It is immediate from (\ref{vij}) that 
\begin{align}\label{F1vij}
F_1 v_{ij} = v_{i+1,j}
\qquad 
\hbox{for all $i,j\in \N$}.
\end{align}
Using (\ref{W:E1E2F1F2}) and (\ref{vij}) yields that 
\begin{align}\label{F2vij}
F_2 v_{ij} &= v_{i,j+1}
\qquad 
\hbox{for all $i,j\in \N$}.
\end{align}
Note that $v=v_{00}$.
Recall from Definition \ref{defn:weight}(i) that $I v_{00}=\mu v_{00}$. Applying $I$ to either side of (\ref{vij}) it follows that 
\begin{align}\label{Ivij}
I v_{ij} &=\mu  v_{ij}
\qquad 
\hbox{for all $i,j\in \N$}.
\end{align}
Recall from Definition \ref{defn:weight}(i) that $K_1 v_{00}=\lambda_1 v_{00}$. 
Using (\ref{W:K1E2}) and (\ref{F2vij}) yields that $K_1 v_{0j}=\lambda_1 v_{0j}$ for all $j\in \N$. 
Applying (\ref{W:E1K1q})
and  (\ref{F1vij}) a routine induction implies that 
\begin{align}\label{K1vij}
K_1 v_{ij} &=\lambda_1 q^{-2i}   v_{ij}
\qquad 
\hbox{for all $i,j\in \N$}.
\end{align}
Recall from Definition \ref{defn:weight}(i) that $K_2 v_{00}=\lambda_2 v_{00}$.
Using (\ref{W:K1E2}) and (\ref{F1vij}) yields that $K_2 v_{i0}=\lambda_2 v_{i0}$ for all $i\in \N$. 
Applying (\ref{W:E1K1q}) and (\ref{F2vij}) a routine induction implies that 
\begin{align}\label{K2vij}
K_2 v_{ij}&=\lambda_2 q^{-2j}  v_{ij}
\qquad 
\hbox{for all $i,j\in \N$}.
\end{align}
Recall from Definition \ref{defn:weight}(ii) that $E_1 v_{00}=0$ and $E_2v_{00}=0$. 
Using (\ref{W:E1E2F1F2}) and (\ref{F2vij}) yields that $E_1 v_{0j}=0$ for all $j\in \N$. 
Applying (\ref{W:E1F1}) along with (\ref{F1vij}), (\ref{Ivij}) and (\ref{K1vij}) a routine induction implies that 
\begin{align}\label{E1vij}
E_1 v_{ij} &=
\left\{
\begin{array}{ll}
0 \qquad &\hbox{if $i=0$},
\\
\displaystyle 
[i]_q 
\frac{\lambda_1 q^{1-i}-\lambda_1^{-1}\mu q^{i-1}}{q-q^{-1}}
v_{i-1,j}
\qquad 
&\hbox{if $i\not=0$}
\end{array}
\right.
\qquad 
\hbox{for all $i,j\in \N$}.
\end{align}
Using (\ref{W:E1E2F1F2}) and (\ref{F1vij}) yields that $E_2 v_{i0}=0$ for all $i\in \N$. 
Applying (\ref{W:E2F2}) along with (\ref{F2vij}), (\ref{Ivij}) and (\ref{K2vij}) a routine induction implies that 
\begin{align}\label{E2vij}
E_2 v_{ij} &=
\left\{
\begin{array}{ll}
0 \qquad &\hbox{if $j=0$},
\\
\displaystyle 
[j]_q
\frac{\lambda_2 \mu q^{1-j}-\lambda_2^{-1} q^{j-1}}{q-q^{-1}}
v_{i,j-1}
\qquad 
&\hbox{if $j\not=0$}
\end{array}
\right.
\qquad 
\hbox{for all $i,j\in \N$}.
\end{align}

Note that $v_{00}$ is nonzero. 
By (\ref{K1vij}) there exists an integer $i\geq 1$ such that $v_{i0}=0$; otherwise the vectors $\{v_{i0}\}_{i\in \N}$ are linearly independent, a contradiction to the finite-dimensionality of $V$. By (\ref{K2vij}) there exists an integer $j\geq 1$ such that $v_{0j}=0$; otherwise the vectors $\{v_{0j}\}_{j\in \N}$ are linearly independent, a contradiction to the finite-dimensionality of $V$. Let
\begin{align}
m &= \min_{i\geq 1}\{i\,|\, v_{i0}=0\}-1,
\label{m}
\\
n &= \min_{j\geq 1}\{j\,|\, v_{0j}=0\}-1.
\label{n}
\end{align}
In addition we let $\delta$ and $\lambda$ denote the nonzero complex scalars given by
\begin{align}
\delta &=\lambda_1 \lambda_2 q^{-m-n},
\label{delta}
\\
\lambda&=\lambda_1 q^{-m}.
\label{lambda}
\end{align}
By (\ref{m}) we have $v_{m+1,0}=0$ and $v_{m0}\not=0$. 
Applying (\ref{E1vij}) with $(i,j)=(m+1,0)$ yields that 
\begin{gather}\label{delta'}
\lambda_1 q^{-m}=\lambda_1^{-1}\mu q^m.
\end{gather}
By (\ref{n}) we have $v_{0,n+1}=0$ and $v_{0n}\not=0$. 
Applying (\ref{E2vij}) with $(i,j)=(0,n+1)$ yields that 
\begin{gather}\label{delta''}
\lambda_2^{-1} q^n=\lambda_2 \mu q^{-n}.
\end{gather}
Dividing (\ref{delta'}) by (\ref{delta''}) yields that $\lambda_1^2 \lambda_2^2 =q^{2(m+n)}$. Hence $\delta\in\{\pm 1\}$ by the setting (\ref{delta}). This shows that $(m,n,\delta,\lambda)\in \I(\W)$. 
Rewriting (\ref{delta})--(\ref{delta'}) yields that 
\begin{gather}\label{lambda123}
(\lambda_1,\lambda_2,\mu)=(\lambda q^m, \delta \lambda^{-1} q^n,\lambda^2)
\end{gather}

Since $v_{m+1,0}=0$ and $v_{0,n+1}=0$ it follows from (\ref{F1vij}) and (\ref{F2vij}) that 
\begin{align}
F_1 v_{ij}&=
\left\{
\begin{array}{ll}
0 \qquad &\hbox{if $i=m$},
\\
v_{i+1,j}
\qquad 
&\hbox{if $i\not= m$},
\end{array}
\right. 
\label{F1vij_new}
\\
F_2 v_{ij}&=
\left\{
\begin{array}{ll}
0 \qquad &\hbox{if $j=n$},
\\
v_{i,j+1}
\qquad 
&\hbox{if $j\not= n$}
\end{array}
\right.
\label{F2vij_new}
\end{align}
for all integers $i$ and $j$ with $0\leq i\leq m$ and $0\leq j\leq n$.
Substituting (\ref{lambda123}) into (\ref{Ivij})--(\ref{E2vij}) yields that 
\begin{align}
I v_{ij}&= \lambda^2 v_{ij},
\label{Ivij_new}
\\
K_1 v_{ij} &=\lambda q^{m-2i} v_{ij},
\label{K1vij_new}
\\
K_2 v_{ij}&= \delta \lambda^{-1} q^{n-2j} v_{ij},
\label{K2vij_new}
\\
E_1 v_{ij}&=
\left\{
\begin{array}{ll}
0 \qquad &\hbox{if $i=0$},
\\
\lambda [i]_q [m-i+1]_q
v_{i-1,j}
\qquad 
&\hbox{if $i\not=0$},
\end{array}
\right.
\label{E1vij_new}
\\
E_2 v_{ij}&=
\left\{
\begin{array}{ll}
0 \qquad &\hbox{if $j=0$},
\\
\delta 
\lambda
[j]_q
[n-j+1]_q
v_{i,j-1}
\qquad 
&\hbox{if $j\not=0$}
\end{array}
\right.
\label{E2vij_new}
\end{align}
for all $i,j\in \N$.
By Proposition \ref{prop:Lmn} and Lemma \ref{lem:delta&lambda} along with (\ref{F1vij_new})--(\ref{E2vij_new}) there exists a unique $\W$-module homomorphism 
\begin{gather}\label{Lmn->V}
(L_m\otimes L_n)^{\delta,\lambda}\to V
\end{gather}
that sends $v_i^{(m)}\otimes v_j^{(n)}$ to $v_{ij}$ for all integers $i$ and $j$ with $0\leq i\leq m$ and $0\leq j\leq n$. 
By Lemma \ref{lem:Lmn_irr}(ii) the $\W$-module $(L_m\otimes L_n)^{\delta,\lambda}$ is irreducible.
Hence the map (\ref{Lmn->V}) is injective. The existence follows. The uniqueness is immediate from Lemma \ref{lem2:irrWmodule}.
\end{proof}

\begin{thm}
\label{thm:irrWmodule'} 
Suppose that $V$ is a finite-dimensional irreducible $\W$-module. Then there exists a unique quadruple $(m,n,\delta,\lambda)\in \I(\W)$ such that 
the $\W$-module $(L_m\otimes L_n)^{\delta,\lambda}$ is isomorphic to $V$. 
\end{thm}
\begin{proof}
Immediate from Lemma \ref{lem:Wmodule_highest} and Theorem \ref{thm:irrWmodule}. 
\end{proof}

\begin{cor}
\label{cor:irrWmodule}
Suppose that $V$ is a finite-dimensional $\W$-module generated by a highest weight vector $v$. For any $(m,n,\delta,\lambda)\in \I(\W)$ the following conditions are equivalent:
\begin{enumerate}
\item The $\W$-module $V$ is isomorphic to $(L_m\otimes L_n)^{\delta,\lambda}$.

\item The highest weight vector $v$ of $V$ has the weight
$(\lambda q^m,\delta \lambda^{-1} q^n,\lambda^2)$.
\end{enumerate}
\end{cor}
\begin{proof}
(i) $\Rightarrow$ (ii): Immediate from Lemma \ref{lem1:irrWmodule}.

(ii) $\Rightarrow$ (i): Immediate from Lemma \ref{lem2:irrWmodule} and Theorem \ref{thm:irrWmodule}.
\end{proof}

As a consequence we have a classification of finite-dimensional irreducible $\W$-modules as follows:
Let $\M(\W)$ denote the set of all isomorphism classes of finite-dimensional irreducible $\W$-modules.
By Lemma \ref{lem:Lmn_irr}(ii) there exists a map $\I(\W)\to \M(\W)$ given by 
\begin{eqnarray*}
(m,n,\delta,\lambda) &\mapsto & 
\hbox{the isomorphism class of the $\W$-module $(L_m\otimes L_n)^{\delta,\lambda}$}
\end{eqnarray*}
for all $(m,n,\delta,\lambda)\in \I(\W)$. By Theorem \ref{thm:irrWmodule'} the map $\I(\W)\to \M(\W)$ is a bijection.

\section{The necessary and sufficient conditions for finite-dimensional $\W$-modules as completely reducible}\label{s:Wreducible}

Recall the elements $\Lambda_1$ and $\Lambda_2$ of $\W$ from (\ref{Lambda1}) and (\ref{Lambda2}).

\begin{lem}
\label{lem:Lmn&Lambda12}
Suppose that $(m,n,\delta,\lambda)\in \I(\W)$. 
Then the following statements hold: 
\begin{enumerate}
\item $\Lambda_1$ acts on $(L_m\otimes L_n)^{\delta,\lambda}$ as scalar multiplication by $\lambda(q^{m+1}+q^{-m-1})$.

\item $\Lambda_2$ acts on $(L_m\otimes L_n)^{\delta,\lambda}$ as scalar multiplication by $\delta\lambda(q^{n+1}+q^{-n-1})$.

\item $I$ acts on $(L_m\otimes L_n)^{\delta,\lambda}$ as scalar multiplication by $\lambda^2$.
\end{enumerate}
\end{lem}
\begin{proof}
(i), (ii): By Lemma \ref{lem:Lambda12&flat} the actions of $\Lambda_1$ and $\Lambda_2$ on the $\W$-module $L_m\otimes L_n$ are identical to the actions of $\Lambda\otimes 1$ and $1\otimes \Lambda$ on the $\U\otimes \U$-module $L_m\otimes L_n$ respectively. 
Combined with Lemma \ref{lem:Ln}(ii) this implies that $\Lambda_1$ and $\Lambda_2$ act on $L_m\otimes L_n$ as scalar multiplication by $q^{m+1}+q^{-m-1}$ and $q^{n+1}+q^{-n-1}$ respectively. By Lemma \ref{lem:delta&lambda} the automorphism $\sigma_{\delta,\lambda}$ sends $\Lambda_1$ to $\lambda \Lambda_1$ and $\Lambda_2$ to $\delta\lambda \Lambda_2$. Therefore (i) and (ii) follow.

(iii): Immediate from Proposition \ref{prop:Lmn} and Lemma \ref{lem:delta&lambda}.
\end{proof}

\begin{lem}
\label{lem2:irrWmodule'}
For any $(m,n,\delta,\lambda),(m',n',\delta',\lambda')\in \I(\W)$ the following conditions are equivalent:
\begin{enumerate}
\item $(m,n,\delta,\lambda)=(m',n',\delta',\lambda')$.

\item $(
\lambda(q^{m+1}+q^{-m-1}),
\delta\lambda(q^{n+1}+q^{-n-1})
,\lambda^2)
=
(
\lambda'(q^{m'+1}+q^{-m'-1}),
\delta'\lambda'(q^{n'+1}+q^{-n'-1})
,\lambda'^2)$.
\end{enumerate}
\end{lem}
\begin{proof}
(i) $\Rightarrow$ (ii): It is obvious. 

(ii) $\Rightarrow$ (i): Let $(c_1,c_2,\mu)$ denote the triple given in (ii). Then 
\begin{align*}
q^{2(m+1)}+q^{-2(m+1)} &= c_1^{2} \mu^{-1}-2=q^{2(m'+1)}+q^{-2(m'+1)},
\\
q^{2(n+1)}+q^{-2(n+1)}&= c_2^{2} \mu^{-1}-2=q^{2(n'+1)}+q^{-2(n'+1)}.
\end{align*}
Since $q$ is not a root of unity and $m,n\in \N$ it follows that $(m,n)=(m',n')$. 
Note that $q^{m+1}+q^{-m-1}$ and $q^{n+1}+q^{-n-1}$ are  nonzero. 
Since $m=m'$ and $\lambda(q^{m+1}+q^{-m-1})=\lambda'(q^{m'+1}+q^{-m'-1})$ it follows that $\lambda=\lambda'$. Since $(n,\lambda)=(n',\lambda')$ and $\delta\lambda(q^{n+1}+q^{-n-1})=\delta'\lambda'(q^{n'+1}+q^{-n'-1})$ it follows that $\delta=\delta'$. The condition (i) holds.
\end{proof}

\begin{cor}
\label{cor:irrWmodule'}
Suppose that $V$ is a finite-dimensional irreducible $\W$-module. For any $(m,n,\delta,\lambda)\in \I(\W)$ the $\W$-module $V$ is isomorphic to $(L_m\otimes L_n)^{\delta,\lambda}$ if and only if the following conditions hold:
\begin{enumerate}
\item $\Lambda_1$ acts on $V$ as scalar multiplication by $\lambda(q^{m+1}+q^{-m-1})$.

\item $\Lambda_2$ acts on $V$ as scalar multiplication by $\delta\lambda(q^{n+1}+q^{-n-1})$.

\item $I$ acts on $V$ as scalar multiplication by $\lambda^2$.
\end{enumerate}
\end{cor}
\begin{proof}
($\Rightarrow$): Immediate from Lemma \ref{lem:Lmn&Lambda12}.

($\Leftarrow$): Immediate from Theorem \ref{thm:irrWmodule'} and Lemmas \ref{lem:Lmn&Lambda12} and \ref{lem2:irrWmodule'}.
\end{proof}

\begin{lem}
\label{lem:F1F2V=0}
Suppose that $V$ is a finite-dimensional $\W$-module. Then there is an integer $n\geq 1$ such that $F_1^n$ and $F_2^n$ vanish on $V$. 
\end{lem}
\begin{proof}
Using (\ref{W:E1K1q}) two routine inductions show that 
\begin{align}
F_1^i V_{K_1}^{(m)}(\theta)\subseteq V_{K_1}^{(m)}(q^{-2i}\theta)
\qquad 
\hbox{for all $\theta\in \C$ and all $i,m\in \N$ with $m\geq 1$},
\label{F1V1n}
\\
F_2^i V_{K_2}^{(m)}(\theta)\subseteq V_{K_2}^{(m)}(q^{-2i}\theta)
\qquad 
\hbox{for all $\theta\in \C$ and all $i,m\in \N$ with $m\geq 1$}.
\label{F2V2n}
\end{align}
Let $\{\theta_i^{(1)}\}_{i=0}^{d_1}$ and $\{\theta_i^{(2)}\}_{i=0}^{d_2}$ denote the eigenvalues of $K_1$ and $K_2$ on $V$ respectively. Since $V$ is finite-dimensional there is an integer $m\geq 1$ such that 
\begin{align}
\label{V}
V=\bigoplus_{i=0}^{d_1} V_{K_1}^{(m)}(\theta_i^{(1)})=\bigoplus_{i=0}^{d_2} V_{K_2}^{(m)}(\theta_i^{(2)}).
\end{align}
Since $K_1$ and $K_2$ are invertible in $\W$ the scalars $\{\theta_i^{(1)}\}_{i=0}^{d_1}$ and $\{\theta_i^{(2)}\}_{i=0}^{d_2}$ are nonzero. 
Since $q$ is not a root of unity there is an integer $n\geq 1$ such that none of $\{q^{-2n}\theta_i^{(1)}\}_{i=0}^{d_1}$ is an eigenvalue of $K_1$ on $V$ and none of $\{q^{-2n}\theta_i^{(2)}\}_{i=0}^{d_2}$ is an eigenvalue of $K_2$ on $V$. 
To see that $F_1^n$ and $F_2^n$ vanish on $V$, one may 
apply (\ref{F1V1n}) and (\ref{F2V2n}) to evaluate the equations obtained by left multiplying (\ref{V}) by $F_1^n$ and $F_2^n$ respectively.
\end{proof}

\begin{lem}
\label{lem:[m]q}
For all $i,m\in \Z$ the following equalities hold:
\begin{align*}
[m]_q&=q^{-i}[m-i]_q+q^{m-i} [i]_q\\
&=q^{i}[m-i]_q+q^{i-m}[i]_q.
\end{align*}
\end{lem}
\begin{proof}
It is straightforward to verify the lemma.
\end{proof}

\begin{lem}
\label{lem:E1F1n}
For all integers $n\geq 1$ the following equations hold in $\W$:
\begin{enumerate}
\item $(q-q^{-1})E_1 F_1^n=(q-q^{-1}) F_1^n E_1+[n]_q (q^{n-1}K_1-q^{1-n}I K_1^{-1})F_1^{n-1}$.

\item $(q-q^{-1})E_2 F_2^n=(q-q^{-1})F_2^n E_2+[n]_q (q^{n-1}I K_2-q^{1-n}K_2^{-1})F_2^{n-1}$.
\end{enumerate}
\end{lem}
\begin{proof}
(i): We proceed by induction on $n$. For $n=1$ it is immediate from (\ref{W:E1F1}). Suppose that $n\geq 2$. By the induction hypothesis the element $(q-q^{-1})E_1 F_1^n$ is equal to 
\begin{gather*}
(q-q^{-1})F_1^{n-1}E_1F_1+[n-1]_q (q^{n-2} K_1-q^{2-n} I K_1^{-1}) F_1^{n-1}.
\end{gather*}
By (\ref{W:E1F1}) the above element is equal to 
\begin{gather}\label{E1F1n}
(q-q^{-1})F_1^n E_1
+
F_1^{n-1}(K_1-I K_1^{-1})
+[n-1]_q (q^{n-2} K_1-q^{2-n} I K_1^{-1}) F_1^{n-1}.
\end{gather}
By (\ref{W:E1K1q}) the sum of the last two terms in  (\ref{E1F1n}) is equal to 
\begin{gather*}
(q^{2(n-1)}+ q^{n-2}[n-1]_q) K_1 F_1^{n-1}-(q^{2(1-n)}+ q^{2-n} [n-1]_q) I K_1^{-1} F_1^{n-1}.
\end{gather*}
Applying Lemma \ref{lem:[m]q} with $(i,m)=(1,n)$ the coefficients of $K_1 F_1^{n-1}$ and $I K_1^{-1} F_1^{n-1}$ are equal to $q^{n-1}[n]_q$ and $-q^{1-n}[n]_q$ respectively.
The equation (i) follows.

(ii): We proceed by induction on $n$. For $n=1$ it is immediate from (\ref{W:E2F2}). Suppose that $n\geq 2$. By the induction hypothesis the element $(q-q^{-1}) E_2 F_2^n$ is equal to 
\begin{gather*}
(q-q^{-1})F_2^{n-1}E_2F_2+[n-1]_q(q^{n-2} I K_2-q^{2-n} K_2^{-1}) F_2^{n-1}.
\end{gather*}
By (\ref{W:E2F2}) the above element is equal to 
\begin{gather}\label{E2F2n}
(q-q^{-1})F_2^{n}E_2
+
F_2^{n-1}( I K_2-K_2^{-1})
+
[n-1]_q(q^{n-2}  I K_2-q^{2-n} K_2^{-1}) F_2^{n-1}.
\end{gather}
By (\ref{W:E1K1q}) the sum of the last two terms in (\ref{E2F2n}) is equal to  
\begin{gather*}
(q^{2n-2}+q^{n-2}[n-1]_q) I K_2 F_2^{n-1}-(q^{2-2n}+q^{2-n}[n-1]_q) K_2^{-1} F_2^{n-1}.
\end{gather*}
Applying Lemma \ref{lem:[m]q} with $(i,m)=(1,n)$ the coefficients of $ I K_2 F_2^{n-1}$ and $K_2^{-1} F_2^{n-1}$ are equal to $q^{n-1}[n]_q$ and $-q^{1-n}[n]_q$ respectively. 
The equation (ii) follows.
\end{proof}

Recall that the $q$-binomial coefficients are defined as 
\begin{gather}\label{qbino}
{n\brack \ell}_q=\prod_{i=1}^\ell
\frac{[n-i+1]_q}{[\ell-i+1]_q}
\qquad 
\hbox{for all $n\in \Z$ and $\ell\in \N$}.
\end{gather}

\begin{lem}
\label{lem:E1mF1n}
For all integers $m$ and $n$ with $0\leq m\leq n$ the following equations hold in $\W$:
\begin{enumerate}
\item $E_1^m F_1^n
=
\displaystyle
\sum_{i=0}^m 
{m\brack i}_q 
\left(
\prod\limits_{j=0}^{i-1}
[n-j]_q
\frac{q^{n-m-j}K_1-q^{m-n+j}I K_1^{-1}}{q-q^{-1}}
\right)
F_1^{n-i}E_1^{m-i}$.

\item $E_2^m F_2^n
=
\displaystyle
\sum_{i=0}^m 
{m\brack i}_q 
\left(
\prod\limits_{j=0}^{i-1}
[n-j]_q
\frac{q^{n-m-j}I K_2-q^{m-n+j}K_2^{-1}}{q-q^{-1}}
\right)
F_2^{n-i}
E_2^{m-i}$.
\end{enumerate}
\end{lem}
\begin{proof}
(i): We proceed by induction on $m$. For $m=0$ there is nothing to prove. Suppose that $m\geq 1$. 
For convenience let 
\begin{gather}\label{Ci}
C_i^{(\ell)}
=\prod\limits_{j=0}^{i-1}
\frac{q^{\ell-j}K_1-q^{j-\ell}I K_1^{-1}}{q-q^{-1}}
\qquad 
\hbox{for all $i,\ell\in \Z$}.
\end{gather}
By the induction hypothesis the element $E_1^m F_1^n$ is equal to 
$$
\sum_{i=0}^{m-1} 
{m-1\brack i}_q 
\left(
\prod\limits_{j=0}^{i-1}
[n-j]_q
\right)
E_1
C_i^{(n-m+1)}
F_1^{n-i}E_1^{m-i-1}.
$$
By (\ref{W:E1K1q}) the above element is equal to 
$$
\sum_{i=0}^{m-1} 
{m-1\brack i}_q 
\left(
\prod\limits_{j=0}^{i-1}
[n-j]_q
\right)
C_i^{(n-m-1)}
E_1 F_1^{n-i}E_1^{m-i-1}.
$$
Applying Lemma \ref{lem:E1F1n}(i) the above element is equal to 
\begin{align*}
&\sum_{i=0}^{m-1} 
{m-1\brack i}_q 
\left(
\prod\limits_{j=0}^{i-1}
[n-j]_q
\right)
C_i^{(n-m-1)}
F_1^{n-i}E_1^{m-i}
\\
&\quad +\;
\sum_{i=0}^{m-1} 
{m-1\brack i}_q 
\left(
\prod\limits_{j=0}^i
[n-j]_q
\right)
C_1^{(n-i-1)}
C_i^{(n-m-1)}
F_1^{n-i-1}E_1^{m-i-1}.
\end{align*}
Rewriting the index $i$ of the second summation yields the above element is equal to 
\begin{align*}
&\sum_{i=0}^{m-1} 
{m-1\brack i}_q 
\left(
\prod\limits_{j=0}^{i-1}
[n-j]_q
\right)
C_i^{(n-m-1)}
F_1^{n-i}E_1^{m-i}
\\
&\quad +\;
\sum_{i=1}^{m} 
{m-1\brack i-1}_q 
\left(
\prod\limits_{j=0}^{i-1}
[n-j]_q
\right)
C_1^{(n-i)}
C_{i-1}^{(n-m-1)}
F_1^{n-i}E_1^{m-i}.
\end{align*}
Combining the above two summations yields that $E_1^m F_1^n$ is equal to $F_1^n E_1^m$ plus
$$
\sum_{i=1}^{m}
{m\brack i}_q
\left(
\prod\limits_{j=0}^{i-1}
[n-j]_q
\right)
C_{i-1}^{(n-m-1)}
\frac{
[m-i]_q
C_1^{(n-m-i)}
+
[i]_q
C_1^{(n-i)}
}{[m]_q}
F_1^{n-i}E_1^{m-i}.
$$
Applying Lemma \ref{lem:[m]q} yields that 
$
[m-i]_q
C_1^{(n-m-i)}
+
[i]_q
C_1^{(n-i)}
=[m]_qC_1^{(n-m)}$. The equation (i) follows.

(ii): We proceed by induction on $m$. For $m=0$ there is nothing to do. Suppose that $m\geq 1$. For convenience let 
\begin{gather}\label{Di}
D_i^{(\ell)}
=
\prod_{j=0}^{i-1}
\frac{q^{\ell-j} I K_2-q^{j-\ell} K_2^{-1}}{q-q^{-1}}
\qquad 
\hbox{for all $i,\ell\in \Z$}.
\end{gather}
By the induction hypothesis the element $E_2^m F_2^n$ is equal to 
$$
\sum_{i=0}^{m-1}
{m-1\brack i}_q
\left(
\prod_{j=0}^{i-1}
[n-j]_q
\right)
E_2 D_i^{(n-m+1)} F_2^{n-i} E_2^{m-i-1}.
$$
By (\ref{W:E1K1q}) the above element is equal to 
$$
\sum_{i=0}^{m-1}
{m-1\brack i}_q
\left(
\prod_{j=0}^{i-1}
[n-j]_q
\right)
D_i^{(n-m-1)} E_2 F_2^{n-i} E_2^{m-i-1}.
$$
Applying Lemma \ref{lem:E1F1n}(ii) the above element is equal to 
\begin{align*}
&\sum_{i=0}^{m-1} 
{m-1\brack i}_q 
\left(
\prod\limits_{j=0}^{i-1}
[n-j]_q
\right)
D_i^{(n-m-1)}
F_2^{n-i}E_2^{m-i}
\\
&\quad +\;
\sum_{i=0}^{m-1} 
{m-1\brack i}_q 
\left(
\prod\limits_{j=0}^i
[n-j]_q
\right)
D_1^{(n-i-1)}
D_i^{(n-m-1)}
F_2^{n-i-1}E_2^{m-i-1}.
\end{align*}
Rewriting the index $i$ of the second summation yields that $E_2^m F_2^n$ is equal to $F_2^n E_2^m$ plus
$$
\sum_{i=1}^m
{m\brack i}_q
\left(
\prod_{j=0}^{i-1}
[n-j]_q
\right)
D_{i-1}^{(n-m-1)}
\frac{[m-i]_q D_1^{(n-m-i)}+[i]_q D_1^{(n-i)}}{[m]_q}
F_2^{n-i} E_2^{m-i}.
$$
Applying Lemma \ref{lem:[m]q} yields that 
$[m-i]_q
D_1^{(n-m-i)}
+
[i]_q
D_1^{(n-i)}
=[m]_qD_1^{(n-m)}$. The equation (ii) follows.
\end{proof}

\begin{lem}
\label{lem:K123}
Suppose that $V$ is a finite-dimensional $\W$-module. Then there exists an integer $n\geq 1$ satisfying the following conditions:
\begin{enumerate}
\item $
\prod\limits_{i=1}^{2n-1}\left(
I^{-1} K_1^2-q^{2(i-n)}
\right)$ vanishes on $V$.

\item $
\prod\limits_{i=1}^{2n-1}\left(
I K_2^2-q^{2(i-n)}
\right)$ vanishes on $V$.
\end{enumerate}
\end{lem}
\begin{proof}
By Lemma \ref{lem:F1F2V=0} there exists an integer $n\geq 1$ such that $F_1^n$ and $F_2^n$ vanish on $V$.

(i): 
Recall from (\ref{Ci}) the elements $C_i^{(\ell)}$ of $\W$ for all $i,\ell\in \Z$. 
It suffices to show that $C_{2n-1}^{(n-1)}$ vanishes on $V$.  To see this we proceed by induction to prove that 
$
F_1^{n-m} C^{(2m-n-1)}_{2m-1}
$ vanishes on $V$ for any integer $m$ with $0\leq m\leq n$. 
For $m=0$ the claim holds by the choice of $n$. Suppose that $m\geq 1$. Consider the element 
\begin{gather}\label{X}
E_1^m F_1^n
C_{m-1}^{(2m-n-1)}.
\end{gather} 
By the choice of $n$ the element (\ref{X}) vanishes on $V$.
By Lemma \ref{lem:E1mF1n}(i) the element (\ref{X}) is equal to 
$$
\sum_{i=0}^m
{m\brack i}_q
\left(
\prod_{j=0}^{i-1}
[n-j]_q
\right)
C_i^{(n-m)} F_1^{n-i} E_1^{m-i}
C_{m-1}^{(2m-n-1)}.
$$
Using (\ref{W:E1K1q}) yields that the above element is equal to 
$$
\sum_{i=0}^m
{m\brack i}_q
\left(
\prod_{j=0}^{i-1}
[n-j]_q
\right)
F_1^{n-i} 
C_i^{(2i-m-n)} 
C_{m-1}^{(2i-n-1)}
E_1^{m-i}.
$$
Since $C_i^{(2i-m-n)} 
C_{m-1}^{(2i-n-1)}=C_{2i-1}^{(2i-n-1)}C_{m-i}^{(-n)}$ for any integer $i$ with $1\leq i\leq m$, the element (\ref{X}) is equal to 
\begin{gather}\label{X1}
F_1^n C_{m-1}^{(-n-1)} E_1^m
+
\sum_{i=1}^m
{m\brack i}_q
\left(
\prod_{j=0}^{i-1}
[n-j]_q
\right)
F_1^{n-i} 
C_{2i-1}^{(2i-n-1)}
C_{m-i}^{(-n)}
E_1^{m-i}.
\end{gather}
Applying (\ref{X1}) to $V$ it follows from the induction hypothesis that 
$F_1^{n-m} C_{2m-1}^{(2m-n-1)}$ vanishes on $V$. The statement (i) follows.

(ii): Recall from (\ref{Di}) the elements $D_i^{(\ell)}$ of $\W$ for all $i,\ell\in \Z$. 
It suffices to show that $D_{2n-1}^{(n-1)}$ vanishes on $V$.  To see this we proceed by induction to prove that 
$
F_2^{n-m} D^{(2m-n-1)}_{2m-1}
$ vanishes on $V$ for any integer $m$ with $0\leq m\leq n$. 
For $m=0$ the claim holds by the choice of $n$. Suppose that $m\geq 1$. Consider the element 
\begin{gather}\label{Y}
E_2^m F_2^n
D_{m-1}^{(2m-n-1)}.
\end{gather} 
By the choice of $n$ the element (\ref{Y}) vanishes on $V$.
By Lemma \ref{lem:E1mF1n}(ii) the element (\ref{Y}) is equal to 
$$
\sum_{i=0}^m
{m\brack i}_q
\left(
\prod_{j=0}^{i-1}
[n-j]_q
\right)
D_i^{(n-m)} F_2^{n-i} E_2^{m-i}
D_{m-1}^{(2m-n-1)}.
$$
Using (\ref{W:E1K1q}) yields that the above element is equal to 
$$
\sum_{i=0}^m
{m\brack i}_q
\left(
\prod_{j=0}^{i-1}
[n-j]_q
\right)
F_2^{n-i} 
D_i^{(2i-m-n)} 
D_{m-1}^{(2i-n-1)}
E_2^{m-i}.
$$
Since $D_i^{(2i-m-n)} 
D_{m-1}^{(2i-n-1)}=D_{2i-1}^{(2i-n-1)}D_{m-i}^{(-n)}$ for any integer $i$ with $1\leq i\leq m$, the element (\ref{Y}) is equal to 
\begin{gather}\label{Y1}
F_2^n D_{m-1}^{(-n-1)} E_2^m
+
\sum_{i=1}^m
{m\brack i}_q
\left(
\prod_{j=0}^{i-1}
[n-j]_q
\right)
F_2^{n-i} 
D_{2i-1}^{(2i-n-1)}
D_{m-i}^{(-n)}
E_2^{m-i}.
\end{gather}
Applying (\ref{Y1}) to $V$ it follows from the induction hypothesis that 
$F_2^{n-m} D_{2m-1}^{(2m-n-1)}$ vanishes on $V$. The statement (ii) follows.
\end{proof}

\begin{lem}
\label{lem:K3V}
Suppose that $V$ is a finite-dimensional $\W$-module. Then $I^{-1} K_1^2,I K_2^2,K_1K_2$ are simultaneously diagonalizable on $V$.
\end{lem}
\begin{proof}
Since $K_1,K_2,I^{\pm 1}$ commute with each other it suffices to show that each of $I^{-1} K_1^2$, $I K_2^2$, $K_1K_2$ is diagonalizable on $V$. 
By Lemma \ref{lem:K123}(i) and since $q$ is not a root of unity, the minimal polynomial of $I^{-1} K_1^2$
on $V$ splits into distinct linear factors. Therefore $I^{-1} K_1^2$ is diagonalizable on $V$.
By Lemma \ref{lem:K123}(ii) the element $I K_2^2$ is diagonalizable on $V$. 

Since $(K_1K_2)^2=(I^{-1} K_1^2)(I K_2^2)$ the element $(K_1K_2)^2$ is diagonalizable on $V$. 
Since $K_1K_2$ is invertible in $\W$, all roots of the minimal polynomial of $(K_1K_2)^2$ on $V$ are nonzero. 
It follows that every root of the minimal polynomial of $K_1K_2$ on $V$ is simple. Therefore $K_1K_2$ is diagonalizable on $V$.
\end{proof}

\begin{thm}
\label{thm:W_reducible}
Suppose that $V$ is a finite-dimensional $\W$-module. Then the following conditions are equivalent:
\begin{enumerate}
\item The $\W$-module $V$ is completely reducible.

\item $K_1$ is diagonalizable on $V$.

\item $K_2$ is diagonalizable on $V$.

\item $I$ is diagonalizable on $V$.
\end{enumerate}
\end{thm}
\begin{proof}
(ii) $\Leftrightarrow$ (iii) $\Leftrightarrow$ (iv): 
Since $K_1,K_2,I$ are invertible in $\W$ the equivalence of (ii)--(iv) is immediate from Lemma \ref{lem:K3V}.

(i) $\Rightarrow$ (ii)--(iv): Immediate from Proposition \ref{prop:Lmn}, Lemma \ref{lem:delta&lambda} and Theorem \ref{thm:irrWmodule'}.

(ii)--(iv) $\Rightarrow$ (i): 
Observe that the generalized eigenspace of any central element of $\W$ in $V$ is a $\W$-submodule of $V$.
Since $\C$ is algebraically closed every finite-dimensional $\W$-module $V$ is a direct sum of the generalized eigenspaces of any specific element of $\W$ in $V$.

In view of the above facts we may assume that the finite-dimensional $\W$-module $V$ is the generalized eigenspaces of $\Lambda_1,\Lambda_2, I$ simultaneously. Let $c_1,c_2,\mu$ denote the eigenvalues of $\Lambda_1,\Lambda_2, I$ on $V$ respectively.  
By Schur's lemma the central elements $\Lambda_1,\Lambda_2, I$ act on each composition factor of $V$ as scalar multiplication by $c_1,c_2,\mu$ respectively. 
Let $\ell$ denote the length of any composition series of $V$.
Let a composition factor of $V$ be given. By Theorem \ref{thm:irrWmodule'} there exists a unique quadruple $(m,n,\delta,\lambda)\in \I(\W)$ such that the given composition factor of $V$ is isomorphic to the $\W$-module $(L_m\otimes L_n)^{\delta,\lambda}$. 
By Lemma \ref{lem:Lmn&Lambda12} the triple $(c_1,c_2,\mu)$ is equal to 
$$
(\lambda(q^{m+1}+q^{-m-1}),\delta\lambda(q^{n+1}+q^{-n-1}),\lambda^2).
$$  
Furthermore, combined with Corollary \ref{cor:irrWmodule'} each composition factor of $V$ is isomorphic to the $\W$-module $(L_m\otimes L_n)^{\delta,\lambda}$. 
Hence 
\begin{gather}\label{dimV}
\dim V=\ell\cdot \dim (L_m\otimes L_n)^{\delta,\lambda}.
\end{gather}

Set  
$(\lambda_1,\lambda_2)=(\lambda q^m,\delta \lambda^{-1} q^n)$ and let 
$$
P=\prod\limits_{i=1}^m(K_1-\lambda_1 q^{-2i})
\prod\limits_{j=1}^n (K_2-\lambda_2 q^{-2j}).
$$ 
By Proposition \ref{prop:Lmn} and Lemma \ref{lem:delta&lambda} the subspace
$
P(L_m\otimes L_n)^{\delta,\lambda}
$
of $(L_m\otimes L_n)^{\delta,\lambda}$ is spanned by $v_0^{(m)}\otimes v_0^{(n)}$. 
Hence the dimension of
$P^\ell V$ is equal to $\ell$. 
Since $P$ commutes with $K_1,K_2,I$ it follows that $P^\ell V$ is invariant under $K_1,K_2,I$. Combined with the conditions (ii)--(iv) the elements $K_1,K_2,I$ are simultaneously diagonalizable on $P^\ell V$. 
Since $\lambda_1$ is the only eigenvalue of $K_1$ on $P^\ell V$ and $\lambda_2$ is the only eigenvalue of $K_2$ on $P^\ell V$, there exists a basis $\{v_i\}_{i=0}^{\ell-1}$ for $P^\ell V$ such that $v_i$ is a weight vector of $V$ with weight $(\lambda_1,\lambda_2,\mu)$ for every integer $i$ with $0\leq i\leq \ell-1$. Note that $q^2\lambda_1$ and $q^2\lambda_2$ are not the eigenvalues of $K_1$ and $K_2$ on $V$, respectively. Using (\ref{W:K1E2}) and (\ref{W:E1K1q}) yields that 
\begin{align*}
&E_1 P^\ell V\subseteq V_{K_1}(q^2\lambda_1)=\{0\},
\\
&E_2 P^\ell V\subseteq V_{K_2}(q^2\lambda_2)=\{0\}.
\end{align*}
Hence each of the vectors $\{v_i\}_{i=0}^{\ell-1}$ is a highest weight vector of $V$. For any integer $i$ with $0\leq i\leq \ell-1$
let $V_i$ denote the $\W$-submodule of $V$ generated by $v_i$. 
By Corollary \ref{cor:irrWmodule} the $\W$-module $V_i$ is isomorphic to the irreducible $\W$-module $(L_m\otimes L_n)^{\delta,\lambda}$ for any integer $i$ with $0\leq i\leq \ell-1$. Since the vectors $\{v_i\}_{i=0}^{\ell-1}$ are linearly independent and by (\ref{dimV}) it follows that
$$
V=\bigoplus_{i=0}^{\ell-1} V_i.
$$
The result follows.
\end{proof}

We end this section with a simple example of a finite-dimensional $\W$-module which is not completely reducible.

\begin{exam}
The matrices 
\begin{align*}
&E_1=E_2=F_1=F_2
=
\begin{pmatrix}
0 &0
\\
0 &0
\end{pmatrix},
\\
&K_1^{\pm 1}= 
\begin{pmatrix}
1  &\pm 1
\\
0 &1
\end{pmatrix},
\quad 
K_2^{\pm 1}
=
\begin{pmatrix}
1  &\mp 1
\\
0 &1
\end{pmatrix},
\quad 
I^{\pm 1}
=
\begin{pmatrix}
1  &\pm 2
\\
0 &1
\end{pmatrix}
\end{align*}
satisfy the relations (\ref{W:Icenter})--(\ref{W:E2F2}). By Theorem \ref{thm:W_reducible} the corresponding two-dimensional $\W$-module is not completely reducible.
\end{exam}

\section{Finite-dimensional irreducible $\H$-modules}\label{s:Hmodule}

In this section we establish some background on finite-dimensional irreducible $\H$-modules.

\begin{prop}
[Section 4.1, \cite{Huang:CG}]
\label{prop:Mnu(abc)}
For any $a,c\in \C$ and nonzero $b,\nu\in \C$ there exists an $\H$-module $M_\nu(a,b,c)$ satisfying the following conditions:
\begin{enumerate}
\item There exists a basis $\{m_i\}_{i\in\N}$ for $M_\nu(a,b,c)$ such that 
\begin{align*}
(A-\theta_i) m_i&=\varphi_i m_{i-1}
\qquad \hbox{for all integers $i\geq 1$},
\qquad 
(A-\theta_0) m_0=0,
\\
(B-\theta_i^*) m_i &= m_{i+1}
\qquad \hbox{for all $i\in \N$},
\\
(C-\theta_i^\e) m_i &=\varphi_i^\e m_{i-1}
\qquad \hbox{for all integers $i\geq 1$},
\qquad 
(C-\theta_0^\e) m_0=0,
\end{align*}
where 
\begin{align*}
\theta_i &= a \nu q^{-2i} 
\qquad 
\hbox{for all $i\in \N$},
\\
\theta_i^* &= b \nu q^{-2i}+b^{-1} \nu^{-1} q^{2i}
\qquad 
\hbox{for all $i\in \N$},
\\
\theta_i^\e &= c \nu^{-1} q^{2i} 
\qquad 
\hbox{for all $i\in \N$},
\\
\varphi_i &=(q^i-q^{-i})(\nu q^{1-i}-\nu^{-1} q^{i-1})(c-ab \nu q^{1-2i})
\qquad \hbox{for all integers $i\geq 1$},
\\
\varphi_i^\e &=(q^i-q^{-i})(\nu q^{1-i}-\nu^{-1} q^{i-1})(a-b^{-1}c \nu^{-1} q^{2i-1})
\qquad \hbox{for all integers $i\geq 1$}.
\end{align*}

\item The elements $\alpha,\beta,\gamma$ act on $M_\nu(a,b,c)$ as scalar multiplication by 
\begin{align*}
\omega&=c(b+b^{-1})+a(\nu q+\nu^{-1} q^{-1}),
\\
\omega^*&=ac,
\\
\omega^\e &=a(b+b^{-1})+c(\nu q+\nu^{-1} q^{-1}),
\end{align*}
respectively.
\end{enumerate}
\end{prop}

\begin{prop}
[Proposition 4.3, \cite{Huang:CG}]
\label{prop:Mnu(abc)_homo}
Let the notations be as in Proposition \ref{prop:Mnu(abc)}. Suppose that $V$ is an $\H$-module which has a vector $v\in V$ satisfying 
\begin{gather*}
Av=\theta_0 v, 
\qquad 
Cv=\theta_0^\e v,
\\
\alpha v=\omega v,
\qquad 
\beta v=\omega^* v,
\qquad 
\gamma v=\omega^\e v.
\end{gather*}
Then there exists a unique $\H$-module homomorphism $M_\nu (a,b,c)\to V$ that sends $m_0$ to $v$. 
\end{prop}

\begin{defn}
[Section 4.1, \cite{Huang:CG}]
\label{defn:Nnu(a,b,c)&Vd(a,b,c)}
Let the notations be as in Proposition \ref{prop:Mnu(abc)}. Assume that $d\in \N$ and set $\nu=q^d$.
\begin{enumerate}
\item Define $N_\nu(a,b,c)$ as the subspace of $M_\nu(a,b,c)$ spanned by $
\{m_i\}_{i=d+1}^\infty$. Equivalently 
$$
N_\nu(a,b,c)=
\prod_{i=0}^d (B-\theta_i^*)\, 
M_\nu(a,b,c).
$$
Since $\varphi_{d+1}$ and $\varphi^\e_{d+1}$ are zero $N_\nu(a,b,c)$ is an $\H$-submodule of $M_\nu(a,b,c)$.

\item Define $V_d(a,b,c)$ as the quotient $\H$-module 
$$
M_\nu(a,b,c)/N_\nu(a,b,c).
$$
\end{enumerate}
\end{defn}

\begin{lem}
\label{lem:Vd(a,b,c)}
For any $a,b,c\in \C$ with $b\not=0$ and $d\in \N$ the following statements hold:
\begin{enumerate}
\item $a[d+1]_q$ is equal to the trace of $A$ on $V_d(a,b,c)$.

\item $(b+b^{-1})[d+1]_q$ is equal to the trace of $B$ on $V_d(a,b,c)$.

\item $c[d+1]_q$ is equal to the trace of $C$ on $V_d(a,b,c)$.

\item $d+1$ is equal to the dimension of $V_d(a,b,c)$.
\end{enumerate}
\end{lem}
\begin{proof}
It is routine to verify the lemma by using Proposition \ref{prop:Mnu(abc)} and Definition \ref{defn:Nnu(a,b,c)&Vd(a,b,c)}.
\end{proof}

\begin{lem}
[Proposition 4.7, \cite{Huang:CG}]
\label{lem:Vd(abc)_irr}
For any $a,b,c\in \C$ with $b\not=0$ and $d\in \N$
the $\H$-module $V_d(a,b,c)$ is irreducible if and only if the following conditions hold:
\begin{enumerate}
\item $a\not=bc q^{d-2i+1}$ for all integers $i$ with $1\leq i\leq d$.

\item $c\not=ab q^{d-2i+1}$ for all integers $i$ with $1\leq i\leq d$.
\end{enumerate}
\end{lem}

\begin{lem}
\label{lem:irrHmodule}
Suppose that $A$ is nilpotent on a nonzero finite-dimensional $\H$-module $V$. Then $\beta$ vanishes at some nonzero vector of $V$.
\end{lem}
\begin{proof}
Since $V$ is nonzero and finite-dimensional over $\C$ there exists an eigenvalue $c$ of $C$ on $V$. Let $w$ denote an eigenvector of $C$ in $V$ corresponding to the eigenvalue $c$. Set 
$w_i=A^i w$ for all $i\in \N$. 
Using (\ref{beta}) a routine induction yields that 
\begin{gather*}
(C-c q^{-2i-2}) w_{i+1}=q^{-i-1}(q^{i+1}-q^{-i-1}) \beta w_i
\qquad \hbox{for all $i\in \N$}.
\end{gather*}
Since $A$ is nilpotent on $V$ there is a $d\in \N$ such that $w_{d+1}=0$ and $w_d\not=0$. By the above equation with $i=d$ and  
since $q$ is not a root of unity it follows that $\beta w_d=0$. The lemma follows.
\end{proof}

We are ready to determine finite-dimensional irreducible $\H$-modules up to isomorphism.

\begin{thm}
\label{thm:irrHmodule}
Suppose that $V$ is a finite-dimensional irreducible $\H$-module. Then there are $a,b,c\in \C$ with $b\not=0$ and $d\in \N$ such that the $\H$-module $V$ is isomorphic to $V_d(a,b,c)$. 
\end{thm}
\begin{proof}
The proof idea is similar to that of the proofs in \cite[Theorem 4.7]{Huang:2015}, \cite[Theorem 6.3]{Huang:BImodule}, \cite[Theorems 6.1 and 8.1]{Huang:DAHAmodule} and \cite[Theorem 6.3]{SH:2019-1}. Since some details are significantly different we emphasize the different parts in this proof. 

We choose $d\in \N$ such that $d+1=\dim V$ and set $\nu=q^d$. 
If $A$ is nilpotent on $V$ then we choose 
$$
a=0.
$$ 
Suppose that $A$ is not nilpotent on $V$. Then there exists a nonzero eigenvalue $\lambda$ of $A$ on $V$. Since $q$ is not a root of unity the numbers $\{\lambda q^{2i}\}_{i\in \N}$ are mutually distinct. Since $V$ is finite-dimensional there exists an integer $h\geq 0$ such that $\lambda q^{2h}$ is an eigenvalue of $A$ on $V$ but $\lambda q^{2h+2}$ is not an eigenvalue of $A$ on $V$. In this case we set 
$$
a=\lambda \nu^{-1} q^{2h}. 
$$
Since $\C$ is algebraically closed there exists a nonzero $\mu\in \C$ such that $\mu+\mu^{-1}$ is an eigenvalue of $B$ on $V$. For any distinct $i,j\in \N$ the scalars $\mu q^{2i}+\mu^{-1} q^{-2i}$ and $\mu q^{2j}+\mu^{-1} q^{-2j}$ are equal if and only if $\mu^2=q^{-2(i+j)}$. Hence there are infinitely many values among the numbers $\{\mu q^{2i}+\mu^{-1} q^{-2i}\}_{i\in \N}$. Since $V$ is finite-dimensional there exists an integer $k\geq 0$ such that $\mu q^{2k}+\mu^{-1} q^{-2k}$ is an eigenvalue of $B$ on $V$ but $\mu q^{2k+2}+\mu^{-1} q^{-2k-2}$ is not an eigenvalue of $B$ on $V$. We choose
$$
b=\mu \nu^{-1}  q^{2k}.
$$
By Schur's lemma and since $q$ is not a root of unity there exists a scalar $c\in \C$ such that $\gamma$ acts on $V$ as scalar multiplication by $a(b+b^{-1})+c(\nu q+\nu^{-1} q^{-1})$. 
Let $\{\theta_i\}_{i=0}^d$, $\{\theta_i^*\}_{i=0}^d$, $\theta_0^\e$, $\{\varphi_i\}_{i=1}^d$ and $\omega,\omega^*,\omega^\e$ denote the parameters described in Proposition \ref{prop:Mnu(abc)} associated with the current values of $a,b,c,\nu$.

We claim that the $\theta_0$-eigenspace of $A$ in $V$ is invariant under $(A-\theta_1)B$. 
Let $v\in V$ with $Av=\theta_0 v$. 
To see the claim it suffices to prove that 
\begin{gather}\label{claim1}
(A-\theta_0)(A-\theta_1)B v=0.
\end{gather}
Applying $v$ to either side of (\ref{gamma}) yields that  
\begin{gather*}
(q^2-q^{-2})Cv=(q-q^{-1})\omega^\e v-qAB v+q^{-1}\theta_0 Bv.
\end{gather*}
Applying $v$ to either side of (\ref{beta}) and using the above equation to eliminate $Cv$, we obtain that  
\begin{gather}
\label{beta_v}
(q-q^{-1})(q^2-q^{-2})\beta v
=(A-a q^{d+2})(A-\theta_1)B v+(q-q^{-1})^2\theta_0 \omega^\e v.
\end{gather}
Suppose that $A$ is nilpotent on $V$. Since $a=0$ the equation (\ref{beta_v}) becomes  
$$
(q-q^{-1})(q^2-q^{-2})\beta v=A^2 B v.
$$
By Lemma \ref{lem:irrHmodule} and Schur's lemma the left-hand side of the above equation is equal to zero. Therefore (\ref{claim1}) holds in this case. Suppose that $A$ is not nilpotent on $V$. Left multiplying either side of (\ref{beta_v}) by $(A-\theta_0)$ it follows that 
$$
(A-a q^{d+2})(A-\theta_0)(A-\theta_1)B v=0.
$$
By the setting of $a$ the scalar $a q^{d+2}$ is not an eigenvalue of $A$ on $V$ and (\ref{claim1}) holds. In light of the claim we may apply the approach of \cite{Huang:2015, SH:2019-1, Huang:BImodule, Huang:DAHAmodule} to obtain the following results: There exists a basis $\{v_i\}_{i=0}^d$ for $V$ such that 
\begin{align*}
(A-\theta_0)v_0&=0,
\\
(A-\theta_i)v_i&=\varphi_i v_{i-1}+\hbox{a linear combination of $\{v_0,v_1,\ldots, v_{i-2}\}$}
\qquad (1\leq i\leq d),
\\
(B-\theta_i^*) v_i&=v_{i+1} 
\qquad (0\leq i\leq d-1),
\qquad 
(B-\theta_d^*) v_d=0,
\\
(\alpha-\omega) v_0&=0.
\end{align*}
Also, there exists a basis $\{w_i\}_{i=0}^d$ for $V$ such that 
\begin{align*}
(B-\theta_0^*)w_0&=0,
\\
(B-\theta_i^*)w_i&=\varphi_i w_{i-1}+\hbox{a linear combination of $\{w_0,w_1,\ldots, w_{i-2}\}$}
\qquad (1\leq i\leq d),
\\
(A-\theta_i) w_i&=w_{i+1} 
\qquad (0\leq i\leq d-1),
\qquad 
(A-\theta_d) w_d=0,
\\
(\beta-\omega^*) w_0 &=0.
\end{align*}
It follows from Schur's lemma that $(\beta-\omega^*) v_0 =0$. Applying $v_0$ to either side of (\ref{gamma}) yields that 
$(C -\theta_0^\e)v_0=0$. 
By Proposition \ref{prop:Mnu(abc)_homo} there exists a unique $\H$-module homomorphism 
\begin{gather}\label{Mnu(a,b,c)->V}
M_\nu(a,b,c)\to V
\end{gather} that sends $m_0$ to $v_0$. Since the matrix representing $B$ with respect to the basis $\{w_i\}_{i=0}^d$ for $V$ is upper triangular with diagonal entries $\{\theta_i^*\}_{i=0}^d$, the element $\prod_{i=0}^d(B-\theta_i^*) $ vanishes on $V$. 
By Definition \ref{defn:Nnu(a,b,c)&Vd(a,b,c)} the map (\ref{Mnu(a,b,c)->V}) induces an $\H$-module homomorphism 
\begin{gather}\label{Vd(a,b,c)->V}
V_d(a,b,c)\to V
\end{gather} 
that sends $m_0+N_\nu(a,b,c)$ to $v_0$. Since $V$ is a $(d+1)$-dimensional irreducible $\H$-module this implies that (\ref{Vd(a,b,c)->V}) is an isomorphism. The result follows.
\end{proof}

\begin{cor}
\label{cor:irrHmodule}
Suppose that $V$ is a finite-dimensional irreducible $\H$-module.  For any $a,b,c\in \C$ with $b\not=0$ and $d\in \N$ the $\H$-module $V$ is isomorphic to $V_d(a,b,c)$ if and only if the following conditions hold:
\begin{enumerate}
\item $a[d+1]_q$ is equal to the trace of $A$ on $V$.

\item $(b+b^{-1})[d+1]_q$ is equal to the trace of $B$ on $V$.

\item $c[d+1]_q$ is equal to the trace of $C$ on $V$.

\item $d+1$ is equal to the dimension of $V$.
\end{enumerate}
\end{cor}
\begin{proof}
($\Rightarrow$): Immediate from Lemma \ref{lem:Vd(a,b,c)}. 

($\Leftarrow$):  By Theorem \ref{thm:irrHmodule} there are $a',b',c'\in \C$ and $d'\in \N$ with $b'\not=0$ such that the $\H$-module $V$ is isomorphic to $V_{d'}(a',b',c')$. Since $q$ is not a root of unity and by Lemma \ref{lem:Vd(a,b,c)} it follows that $(a',c',d')=(a,c,d)$ and $b'\in \{b^{\pm 1}\}$. Hence the $\H$-module $V$ is isomorphic to $V_d(a,b,c)$ or $V_d(a,b^{-1},c)$. By \cite[Lemma 4.6]{Huang:CG} the $\H$-module $V_d(a,b,c)$ is isomorphic to $V_d(a,b^{-1},c)$. The corollary follows.
\end{proof}

\begin{cor}
\label{cor':irrHmodule}
Suppose that $V$ and $V'$ are finite-dimensional irreducible $\H$-modules with dimensions at least two. 
Then the following conditions are equivalent:
\begin{enumerate}
\item The $\H$-module $V$ is isomorphic to $V'$.

\item There exists a linear isomorphism 
$V\to V'$
which commutes with $A$ and $B$.
\end{enumerate}
\end{cor}
\begin{proof}
(i) $\Rightarrow$ (ii): Trivial.

(ii) $\Rightarrow$ (i): 
By Theorem \ref{thm:irrHmodule} there are $a,b,c\in \C$ with $b\not=0$ and an integer $d\geq 1$ such that the $\H$-module $V_d(a,b,c)$ is isomorphic to $V$. 
By (ii) there exists a linear isomorphism 
\begin{gather}
\label{dim>=2}
V_d(a,b,c)\to V'
\end{gather}
which commutes with $A$ and $B$. Combined with Corollary \ref{cor:irrHmodule} there is a scalar $c'\in \C$ such that the $\H$-module $V'$ is isomorphic to $V_d(a,b,c')$.  
To see (i) it remains to prove that $c=c'$.

Let $\{\theta_i\}_{i=0}^d$, $\{\theta_i^*\}_{i=0}^d$, $\{\varphi_i\}_{i=1}^d$ denote the parameters described in Proposition \ref{prop:Mnu(abc)}(i) associated with the current parameters $a,b,c$ and $\nu=q^d$. 
By Definition \ref{defn:Nnu(a,b,c)&Vd(a,b,c)} there exists a basis $\{v_i\}_{i=0}^d$ for $V_d(a,b,c)$ such that 
\begin{align*}
(A-\theta_i )v_i&=\varphi_i v_{i-1}
\qquad (1\leq i\leq d), 
\qquad 
(A-\theta_0) v_0=0,
\\
(B-\theta_i^*) v_i&=v_{i+1}
\qquad 
(0\leq i\leq d-1),
\qquad 
(B-\theta_d^*) v_d=0.
\end{align*}
In addition there exists a basis $\{w_i\}_{i=0}^d$ for $V'$ such that 
\begin{align*}
(A-\theta_i )w_i&=\varphi_i' w_{i-1}
\qquad (1\leq i\leq d), 
\qquad 
(A-\theta_0) w_0=0,
\\
(B-\theta_i^*) w_i&=w_{i+1}
\qquad 
(0\leq i\leq d-1),
\qquad 
(B-\theta_d^*) w_d=0,
\end{align*} 
where the scalar $\varphi_i'$ is obtained from the above setting of $\varphi_i$ by replacing $c$ with $c'$ for each integer $i$ with $1\leq i\leq d$. 
By Lemma \ref{lem:Vd(abc)_irr} the scalars $\varphi_i\not=0$ $(1\leq i\leq d)$ and the scalars $\varphi_i'\not=0$ $(1\leq i\leq d)$. 
Hence the $\theta_0$-eigenspace of $A$ on $V_d(a,b,c)$
is spanned by $v_0$ and 
that on $V'$ is spanned by $w_0$. 
Since (\ref{dim>=2}) commutes with $A$ there exists a nonzero scalar $\lambda\in \C$ such that (\ref{dim>=2}) sends $v_0$ to $\lambda w_0$.
Since (\ref{dim>=2}) commutes with $B$ the map (\ref{dim>=2}) sends $v_i$ to $\lambda w_i$ for all integers $i$ with $0\leq i\leq d$. 
Hence $\varphi_i=\varphi_i'$ for all integers $i$ with $1\leq i\leq d$. Since $d \geq 1$ it follows that $c=c'$. The corollary follows.
\end{proof}

We summarize this section as follows: 
Let $\M(\H)$ denote the set of all isomorphism classes of finite-dimensional irreducible $\H$-modules. Let $\I(\H)$ denote the set consisting of all quadruples $(a,b,c,d)$ where $a,b,c\in \C$ with $b\not=0$ and $d\in \N$ satisfy Lemma \ref{lem:Vd(abc)_irr}(i), (ii).  By Lemma \ref{lem:Vd(abc)_irr} there exists a map $\I(\H)\to \M(\H)$ given by 
\begin{eqnarray*}
(a,b,c,d) &\mapsto &
\hbox{the isomorphism class of the $\H$-module $V_d(a,b,c)$}
\end{eqnarray*}
for all $(a,b,c,d)\in \I(\H)$. By Theorem \ref{thm:irrHmodule} the map $\I(\H)\to \M(\H)$ is onto. Observe that there is an action of the multiplicative group $\{\pm 1\}$ on $\I(\H)$ given by
\begin{gather*}
(a,b,c,d)^{-1}=(a,b^{-1},c,d) 
\qquad \hbox{for all $(a,b,c,d)\in \I(\H)$}. 
\end{gather*}
Let $\I(\H)/\{\pm 1\}$ denote the set of all $\{\pm1\}$-orbits of $\I(\H)$. Combined with Corollary \ref{cor:irrHmodule} the map $\I(\H)\to \M(\H)$ induces a bijection $\I(\H)/\{\pm 1\}\to \M(\H)$.

\section{The decomposition of the $\H$-module $(L_m\otimes L_n)^{\delta,\lambda}$}
\label{s:HWmodule}

Recall the algebra homomorphism $\widetilde{\Delta}:\U\to \W$ from Theorem \ref{thm1:U->W}.  
Recall the algebra homomorphism $\widetilde{\natural}:\H\to \W$ from Theorem \ref{thm1:H->W}.

\begin{lem}
\label{lem:centralizer}
$\widetilde{\Delta}(K)$ is in the centralizer of $\widetilde{\natural}(\H)$ in $\W$. 
\end{lem}
\begin{proof}
Apparently $\widetilde{\Delta}(K)$ commutes with $\widetilde{\natural}(A)$. Since $[K,\Lambda]=0$ the element $\widetilde{\Delta}(K)$  commutes with $\widetilde{\natural}(B)$. By Lemma \ref{lem:Lambda12_central} the element $\widetilde{\Delta}(K)$  commutes with $\widetilde{\natural}(\gamma)$.
Since the algebra $\H$ is generated by $A,B,\gamma$ the lemma follows.
\end{proof}

\begin{defn}
\label{defn:V(theta)}
Assume that $V$ is a $\W$-module. 
Recall from Definition \ref{defn:VX(theta)} the notation $V_X(\theta)$  where $X\in \W$ and $\theta\in \C$. For convenience we write 
$$
V(\theta)=V_{\widetilde{\Delta}(K)}(\theta)
\qquad 
\hbox{for any $\theta\in \C$}.
$$ 
\end{defn}

Recall that every $\W$-module can be viewed as an $\H$-module by pulling back via $\widetilde{\natural}$.

\begin{lem}
\label{lem:widetildeK&Hmodule}
Suppose that $V$ is a $\W$-module. Then $V(\theta)$ is an $\H$-submodule of $V$ for any $\theta\in \C$.
\end{lem}
\begin{proof}
Immediate from Lemma \ref{lem:centralizer}.
\end{proof}

Recall from Section \ref{s:Wmodule} the $\W$-module $(L_m\otimes L_n)^{\delta,\lambda}$ for any $(m,n,\delta,\lambda)\in \I(\W)$.

\begin{thm}
\label{thm:Lmn_Hmodule}
Suppose that $(m,n,\delta,\lambda)\in \I(\W)$. Then the following statements hold: 
\begin{enumerate}
\item The $\H$-module $(L_m\otimes L_n)^{\delta,\lambda}$ is equal to 
\begin{gather*}
\bigoplus_{\ell=0}^{m+n}
(L_m\otimes L_n)^{\delta,\lambda}(\delta q^{m+n-2\ell}).
\end{gather*}

\item Suppose that $\ell$ is an integer with $0\leq \ell \leq m+n$. Then the $\H$-module $(L_m\otimes L_n)^{\delta,\lambda}(\delta q^{m+n-2\ell})$ is isomorphic to the irreducible $\H$-module $V_d(a,b,c)$ where 
\begin{align}
a &=\delta \lambda q^{\ell-n-\min\{m,\ell\}+\min\{n,\ell\}},
\label{Lmn_a}
\\
b&=\delta q^{m+n+\ell-\min\{m,\ell\}-\min\{n,\ell\}+1},
\label{Lmn_b}
\\
c&=\lambda q^{\ell-m-\min\{n,\ell\}+\min\{m,\ell\}},
\label{Lmn_c}
\\
d&=\min\{m,\ell\}+\min\{n,\ell\}-\ell.
\label{Lmn_d}
\end{align}
\end{enumerate}
\end{thm}
\begin{proof}
(i): Immediate from Proposition \ref{prop:Lmn} and Lemma \ref{lem:delta&lambda}.

(ii): For convenience we write $V=(L_m\otimes L_n)^{\delta,\lambda}(\delta q^{m+n-2\ell})$. 
Let $a,b,c,d$ denote the values (\ref{Lmn_a})--(\ref{Lmn_d}) and we set $\nu=q^d$. Let $\{\theta_i\}_{i\in \N}$, $\{\theta_i^*\}_{i\in \N}$, $\theta_0^\e$, $\varphi_{d+1}$
and $\omega,\omega^*,\omega^\e$ denote the parameters given in Proposition \ref{prop:Mnu(abc)} associated with the current values of $a,b,c,\nu$.

By Proposition \ref{prop:Lmn} and Lemma \ref{lem:delta&lambda} the vectors
\begin{gather*}
v^{(m)}_{\ell-k}
\otimes 
v^{(n)}_k
\qquad 
(\max\{m,\ell\}-m\leq k\leq \min\{n,\ell\})
\end{gather*}
form a basis for $V$. Replacing the index $k$ by $\min\{n,\ell\}-i$ the expression for the above vectors becomes 
\begin{gather*}
v_i=v^{(m)}_{\ell-\min\{n,\ell\}+i}
\otimes 
v^{(n)}_{\min\{n,\ell\}-i}
\qquad 
(0\leq i\leq d).
\end{gather*}
By (\ref{A:widenatural}) the action of $A$ on $V$ is as follows:
\begin{gather}
\label{Avi}
Av_i=\theta_iv_i
\qquad 
(0\leq i\leq d).
\end{gather}
Using (\ref{C:widenatural}) yields that 
$Cv_0=I K_1^{-1} v_0=\theta_0^\e v_0$. Using (\ref{beta:widenatural}) yields that $\beta v_0=\omega^* v_0$. By Lemma \ref{lem:Lmn&Lambda12}(i), (ii) the elements $\alpha$ and $\gamma$ act on $V$ as scalar multiplication by 
\begin{gather}
\delta\lambda
(
q^{n+1}+q^{-n-1}+q^{2\ell-m-n} (q^{m+1}+q^{-m-1})
),
\label{alpha_Lmn}
\\
\lambda 
(q^{m+1}+q^{-m-1}+q^{2\ell-m-n} (q^{n+1}+q^{-n-1})
),
\label{gamma_Lmn}
\end{gather}
respectively.
It is routine to verify that (\ref{alpha_Lmn}) is equal to $\omega$ and  (\ref{gamma_Lmn}) is equal to $\omega^\e$ by dividing the argument into the cases:
(1) $\ell>\max\{m,n\}$; (2) $n\leq \ell\leq m$; (3) $m\leq \ell\leq n$; (4) $\ell<\min\{m,n\}$.

By Proposition \ref{prop:Mnu(abc)_homo} there exists a unique $\H$-module homomorphism 
\begin{gather}\label{Mnu(a,b,c)->Lmn(delta,lambda)}
M_\nu(a,b,c) \to V
\end{gather}
that sends $m_0$ to $v_0$. Since $\varphi_{d+1}=0$ it follows from Proposition \ref{prop:Mnu(abc)}(i) that 
$
(A-\theta_{d+1}) m_{d+1}=0$. 
By (\ref{Avi}) the scalar $\theta_{d+1}$ is not an eigenvalue of $A$ on $V$. Hence $m_{d+1}$ lies in the kernel of (\ref{Mnu(a,b,c)->Lmn(delta,lambda)}). By Proposition \ref{prop:Mnu(abc)}(i) we have 
$$
m_i=
\prod_{h=d+1}^{i-1} (B-\theta_h^*) 
\,
m_{d+1}
\qquad 
\hbox{for all integers $i\geq d+1$}.
$$
Hence the vectors $\{m_i\}_{i=d+1}^\infty$ are in the kernel of (\ref{Mnu(a,b,c)->Lmn(delta,lambda)}). By Definition \ref{defn:Nnu(a,b,c)&Vd(a,b,c)} the map (\ref{Mnu(a,b,c)->Lmn(delta,lambda)}) induces an $\H$-module homomorphism 
\begin{gather}\label{Vd(a,b,c)->Lmn(delta,lambda)}
V_d(a,b,c) \to V
\end{gather}
that sends $m_0+N_\nu(a,b,c)$ to $v_0$. In each of the cases (1)--(4) it is straightforward to check that $(a,b,c,d)\in \I(\H)$. By Lemma \ref{lem:Vd(abc)_irr}
the $\H$-module $V_d(a,b,c)$ is irreducible. Hence (\ref{Vd(a,b,c)->Lmn(delta,lambda)}) is an isomorphism. The result follows. 
\end{proof}

\begin{lem}
\label{lem:Lmn_Hmodule_iso}
Suppose that $(m,n,\delta,\lambda), (m',n',\delta',\lambda')\in \I(\W)$ and $\ell,\ell'$ are two integers with $0\leq \ell \leq m+n$ and $0\leq \ell'\leq m'+n'$. Then the following conditions are equivalent: 
\begin{enumerate}
\item The $\H$-module $(L_m\otimes L_n)^{\delta,\lambda}(\delta q^{m+n-2\ell})$ is isomorphic to $(L_{m'}\otimes L_{n'})^{\delta',\lambda'}(\delta' q^{m'+n'-2\ell'})$.  

\item The following equations hold:
\begin{align}
&\min\{m,\ell\}-\min\{n,\ell\}-m=\min\{m',\ell'\}-\min\{n',\ell'\}-m',
\label{Lmn_Hmodule_iso1}
\\
&\min\{m,\ell\}+\min\{n,\ell\}-\ell=\min\{m',\ell'\}+\min\{n',\ell'\}-\ell',
\label{Lmn_Hmodule_iso2}
\\
&m+n=m'+n',
\label{Lmn_Hmodule_iso3}
\\
&\lambda q^{\ell}=\lambda'q^{\ell'},
\label{Lmn_Hmodule_iso4}
\\
&\delta=\delta'.
\label{Lmn_Hmodule_iso5}
\end{align}
\end{enumerate}
\end{lem}
\begin{proof}
Let $a,b,c,d$ denote the parameters (\ref{Lmn_a})--(\ref{Lmn_d}). Let $a',b',c',d'$ denote the parameters (\ref{Lmn_a})--(\ref{Lmn_d}) with the substitution $(m,n,\delta,\lambda,\ell)=(m',n',\delta',\lambda',\ell')$. 
By Theorem \ref{thm:Lmn_Hmodule}(ii) the $\H$-modules $(L_m\otimes L_n)^{\delta,\lambda}(\delta q^{m+n-2\ell})$ and $(L_{m'}\otimes L_{n'})^{\delta',\lambda'}(\delta' q^{m'+n'-2\ell'})$ are isomorphic to the irreducible $\H$-modules $V_d(a,b,c)$ and $V_{d'}(a',b',c')$ respectively. 
Using (\ref{Lmn_b}) and (\ref{Lmn_d}) yields that $b=\delta q^{m+n-d+1}$ and $b'=\delta' q^{m'+n'-d'+1}$. 
By (\ref{Lmn_d}) the inequalities $d\leq m+n$ and $d'\leq m'+n'$ hold. 
Since $q$ is not a root of unity it follows that $b^{-1}\not=b'$.
Combined with Corollary \ref{cor:irrHmodule} the condition (i) is equivalent to 
\begin{gather*}
(a,b,c,d)=(a',b',c',d').
\end{gather*}

(ii) $\Rightarrow$ (i):  By (\ref{Lmn_Hmodule_iso1}) and (\ref{Lmn_Hmodule_iso3})--(\ref{Lmn_Hmodule_iso5}) the equality $a=a'$ holds. 
By (\ref{Lmn_Hmodule_iso2}), (\ref{Lmn_Hmodule_iso3}) and (\ref{Lmn_Hmodule_iso5}) the equality $b=b'$ holds. 
By (\ref{Lmn_Hmodule_iso1}) and (\ref{Lmn_Hmodule_iso4}) the equality $c=c'$ holds. 
By (\ref{Lmn_Hmodule_iso2}) the equality $d=d'$ holds. 

(i) $\Rightarrow$ (ii):
Since $d=d'$ the equation (\ref{Lmn_Hmodule_iso2}) follows. 
The equation (\ref{Lmn_Hmodule_iso3}) follows by using (\ref{Lmn_Hmodule_iso2}) to simplify $b^2=b'^2$. 
The equation (\ref{Lmn_Hmodule_iso5}) follows by using (\ref{Lmn_Hmodule_iso2}) and (\ref{Lmn_Hmodule_iso3}) to simplify $b=b'$. 
The equation (\ref{Lmn_Hmodule_iso1}) follows by using (\ref{Lmn_Hmodule_iso3}) to simplify $\delta a/c=\delta' a'/c'$. 
The equation (\ref{Lmn_Hmodule_iso4}) follows by using (\ref{Lmn_Hmodule_iso1}) to simplify $c=c'$.
\end{proof}

\begin{thm}
\label{thm:Lmn_Hmodule_iso}
Suppose that $(m,n,\delta,\lambda), (m',n',\delta',\lambda')\in \I(\W)$ and $\ell,\ell'$ are two integers with $0\leq \ell \leq m+n$ and $0\leq \ell'\leq m'+n'$. Then the following conditions are equivalent:
\begin{enumerate}
\item The $\H$-module $(L_m\otimes L_n)^{\delta,\lambda}(\delta q^{m+n-2\ell})$ is isomorphic to $(L_{m'}\otimes L_{n'})^{\delta',\lambda'}(\delta' q^{m'+n'-2\ell'})$.

\item  $\delta=\delta'$ and 
\begin{gather}
\label{Lmn_Hmodule_iso}
(m',n',\ell',\lambda')
\in 
\left\{
\begin{split}
&(m,n,\ell,\lambda),
(\ell,m+n-\ell,m,\lambda q^{\ell-m}),
\\
&(m+n-\ell,\ell,n,\lambda q^{\ell-n}),
(n,m,m+n-\ell,\lambda q^{2\ell-m-n})
\end{split}
\right\}.
\end{gather}
\end{enumerate}
\end{thm}
\begin{proof}
By symmetry it is enough to deal with the following  cases: 
\begin{multicols}{2}
\begin{enumerate}
\item[(1)] $\ell>\max\{m,n\}$ and $\ell'>\max\{m',n'\}$;
\item[(2)] $\ell>\max\{m,n\}$ and $n'\leq \ell'\leq m'$;
\item[(3)] $\ell>\max\{m,n\}$ and $m'\leq \ell'\leq n'$;
\item[(4)] $\ell>\max\{m,n\}$ and $\ell'<\min\{m', n'\}$;
\item[(5)] $n\leq \ell\leq m$ and $n'\leq \ell'\leq m'$;
\item[(6)] $n\leq \ell\leq m$ and $m'\leq \ell'\leq n'$;
\item[(7)] $n\leq \ell\leq m$ and $\ell'<\min\{m', n'\}$;
\item[(8)] $m\leq \ell\leq n$ and $m'\leq \ell'\leq n'$;
\item[(9)] $m\leq \ell\leq n$ and $\ell'<\min\{m',n'\}$;
\item[(10)] $\ell<\min\{m, n\}$ and $\ell'<\min\{m',n'\}$. 
\end{enumerate}
\end{multicols}
\noindent In each case the condition (\ref{Lmn_Hmodule_iso}) is equivalent to the equations (\ref{Lmn_Hmodule_iso1})--(\ref{Lmn_Hmodule_iso4}). More precisely both are equivalent to

\begin{gather*}
(m',n',\ell',\lambda')=
\left\{
\begin{array}{ll}
(m,n,\ell,\lambda)
\qquad 
&\hbox{in the cases (1), (5), (8) and (10)};
\\
(\ell,m+n-\ell,m,\lambda q^{\ell-m})
\qquad 
&\hbox{in the cases (2) and (9)};
\\
(m+n-\ell,\ell,n,\lambda q^{\ell-n})
\qquad 
&\hbox{in the cases (3) and (7)};
\\
(n,m,m+n-\ell,\lambda q^{2\ell-m-n})
\qquad 
&\hbox{in the cases (4) and (6)}.
\end{array}
\right.
\end{gather*}
By Lemma \ref{lem:Lmn_Hmodule_iso} the conditions (i) and (ii) are equivalent. 
\end{proof}

\section{The triple coordinate system for $(\L(\Omega),\subseteq)$}\label{s:Dunkl}

From now on we will adopt the following conventions:
Let $\Omega$ denote a vector space and assume that there are two subspaces $x_0$ and $x_1$ of $\Omega$ with $\Omega=x_0\oplus x_1$. 
For any $u\in \Omega$ we write $u_0$ and $u_1$ for the unique vectors $u_0\in x_0$ and $u_1\in x_1$ with $u=u_0+u_1$.

Recall that $\L(\Omega)$ denotes the set of all subspaces of $\Omega$. 
In this section we recall the triple coordinate system for the poset $(\L(\Omega),\subseteq)$ introduced by Dunkl in \cite[Section 4]{Dunkl77}.

\begin{defn}
[Section 4, \cite{Dunkl77}]
\label{defn:Lx0}
Define $\L(\Omega)_{x_0}$ to be the set consisting of all triples $(y,z,\tau)$ where $y\in \L(\Omega/x_0)$, $z\in \L(x_0)$ and $\tau$ is a linear map from $y$ into $x_0/z$. 
For any $(y,z,\tau), (y',z',\tau')\in \L(\Omega)_{x_0}$ we write $(y,z,\tau)\subseteq (y',z',\tau')$ whenever the following conditions hold:
\begin{enumerate}
\item $y\subseteq y'$.

\item $z\subseteq z'$.

\item $\tau(u)\subseteq \tau'(u)$ for all $u\in y$.
\end{enumerate}
Note that $(\L(\Omega)_{x_0},\subseteq)$ is a poset.
\end{defn}

\begin{defn}
\label{defn:taux}
For any $x\in \L(\Omega)$ the linear map $\t(x):x+x_0/x_0\to x_0/x\cap x_0$ is defined by 
\begin{eqnarray*}
u+x_0
&\mapsto &
u_0+(x\cap x_0)
\qquad 
\hbox{for all $u\in x$}.
\end{eqnarray*}
\end{defn}

\begin{thm}
[Theorem 4.1 and Proposition 4.2, \cite{Dunkl77}]
\label{thm:L(Omega)->Lx0(Omega)}
The map 
$\Px:\L(\Omega)
\to 
\L(\Omega)_{x_0}$ given by 
\begin{eqnarray*}
x &\mapsto & (x+x_0/x_0, x\cap x_0, \t(x))
\qquad 
\hbox{for all $x\in \L(\Omega)$}
\end{eqnarray*}
is a bijection. Moreover, for any $x,x'\in \L(\Omega)$ the following conditions are equivalent:
\begin{enumerate}
\item $x\subseteq x'$.

\item $\Px(x)\subseteq \Px(x')$.
\end{enumerate}
\end{thm}
\begin{proof} 
We associate a triple $(y,z,\tau)\in \L(\Omega)_{x_0}$ to the subspace of $\Omega$ consisting of all $u\in \Omega$ for which $u_1+x_0\in y$ and $u_0\in\tau(u_1+x_0)$. This gives the inverse of $\Px$. 

(i) $\Rightarrow$ (ii): Immediate from Definitions \ref{defn:Lx0} and  \ref{defn:taux}.

(ii) $\Rightarrow$ (i): Let $u\in x$ be given. By Definitions \ref{defn:Lx0} and  \ref{defn:taux} there exists a vector $u'\in x'$ such that $u+x_0=u'+x_0$ and $u_0+(x\cap x_0)\subseteq u_0'+(x'\cap x_0)$. 
Since $u-u'\in x_0$ it follows that $u_1=u_1'$.  
Since $u_0+(x'\cap x_0)= u_0'+(x'\cap x_0)$ 
it follows that 
$$
u-u'=u_0-u_0'\in x'\cap x_0.
$$ 
Hence $u\in x'$. The result follows.
\end{proof}

In view of Theorem \ref{thm:L(Omega)->Lx0(Omega)} we may identify the poset $(\L(\Omega),\subseteq)$ with the poset $(\L(\Omega)_{x_0},\subseteq)$. 
We end this section with the linear endomorphisms of $\C^{\L(\Omega)}$ mentioned in \cite[Definition 4.16]{Dunkl77}.

\begin{defn}
[Definition 4.16, \cite{Dunkl77}]
\label{defn:L12&R12}
\begin{enumerate}
\item Define $L_1(x_0)$ to be the linear map $\C^{\L(\Omega)}\to \C^{\L(\Omega)}$ given by 
\begin{eqnarray*}
x &\mapsto &
\sum_{\substack{x'\subsetdot x
\\
x\cap x_0=x'\cap x_0}} x'
=
\sum_{(y,x\cap x_0,\tau)\subsetdot \Px(x)}
(y,x\cap x_0,\tau)
\qquad 
\hbox{for all $x\in \L(\Omega)$}.
\end{eqnarray*}

\item Define $L_2(x_0)$ to be the linear map $\C^{\L(\Omega)}\to \C^{\L(\Omega)}$ given by 
\begin{eqnarray*}
x &\mapsto &
\sum_{\substack{x'\subsetdot x
\\
x+x_0/x_0=x'+x_0/x_0}} x'
=
\sum_{(x+x_0/x_0,z,\tau)\subsetdot \Px(x)}
(x+x_0/x_0,z,\tau)
\qquad 
\hbox{for all $x\in \L(\Omega)$}.
\end{eqnarray*}

\item Define $R_1(x_0)$ to be the linear map $\C^{\L(\Omega)}\to \C^{\L(\Omega)}$ given by 
\begin{eqnarray*}
x &\mapsto &
\sum_{\substack{x\subsetdot x'
\\
x\cap x_0=x'\cap x_0}} x'
=
\sum_{\Px(x)\subsetdot (y,x\cap x_0,\tau)}
(y,x\cap x_0,\tau)
\qquad 
\hbox{for all $x\in \L(\Omega)$}.
\end{eqnarray*}

\item Define $R_2(x_0)$ to be the linear map $\C^{\L(\Omega)}\to \C^{\L(\Omega)}$ given by 
\begin{eqnarray*}
x &\mapsto &
\sum_{\substack{x\subsetdot x'
\\
x+x_0/x_0=x'+x_0/x_0}} x'
=
\sum_{\Px(x)\subsetdot (x+x_0/x_0,z,\tau)}
(x+x_0/x_0,z,\tau)
\qquad 
\hbox{for all $x\in \L(\Omega)$}.
\end{eqnarray*}
\end{enumerate}
\end{defn}

\section{The $\U$-module $\C^{\L(\Omega)}(\lambda)$}\label{s:L(Omega)}

In the rest of this paper we always assume that the vector space $\Omega$ is of finite dimension $D$ over a finite field $\F$. In addition we set the parameter 
$$
q=\sqrt{|\F|}.
$$

Recall the following formula for the number of all $k$-dimensional subspaces of an $n$-dimensional vector space over $\F$:

\begin{lem}
\label{lem:|Lk|}
For any integers $k$ and $n$ with $0\leq k\leq n$ the following equation holds:
\begin{gather*}
|\L_k(\F^n)|=q^{k(n-k)}
{n\brack k}_q.
\end{gather*}
\end{lem}

\begin{prop}\label{prop:L(Omega)_Umodule}
Suppose that $\lambda$ is a nonzero scalar in $\C$. Then there exists a unique $\U$-module $\C^{\L(\Omega)}$ such that  
\begin{align*}
E x&=\lambda q^{-D}\sum_{x'\subsetdot x} x'
\qquad 
\hbox{for all $x\in \L(\Omega)$},
\\
F x&=\lambda^{-1} q \sum_{x\subsetdot x'} x'
\qquad 
\hbox{for all $x\in \L(\Omega)$},
\\
K^{\pm 1} x&= q^{\pm(D-2 \dim x)}x
\qquad 
\hbox{for all $x\in \L(\Omega)$}.
\end{align*}
\end{prop}
\begin{proof} 
For the moment set $E,F,K^{\pm 1}$ to be the linear endomorphisms of $\C^{\L(\Omega)}$ stated in Proposition \ref{prop:L(Omega)_Umodule}. 
By construction the relations (\ref{U1}) and (\ref{U2}) hold. 
Let $x\in \L(\Omega)$ be given. 
Observe that 
$
\{x'\in \L_{\dim x}(\Omega)\,|\,x\subsetdot x+x'\}=\{x'\in \L_{\dim x}(\Omega)\,|\,x\cap x' \subsetdot x\}.
$
For convenience let $S$ denote the above set. 
A direct calculation yields that 
\begin{align*}
EFx&=
q^{1-D}
\left(
|\L_1(\Omega/x)|x+
\sum_{x'\in S} x'
\right),
\\
FEx&
=q^{1-D}\left(
|\L_{\dim x-1}(x)|x+
\sum_{x'\in S} x'
\right).
\end{align*}
Hence
$
(EF-FE)x
=q^{1-D}\left(
|\L_1(\Omega/x)|
-
|\L_{\dim x-1}(x)|
\right)
x$.
Applying Lemma \ref{lem:|Lk|} to evaluate the above equation yields that 
$$
(EF-FE)x=[D-2\dim x]_q x.
$$
Hence the relation (\ref{U3}) holds. 
By Definition \ref{defn:U} the existence of the $\U$-module $\C^{\L(\Omega)}$ follows. 
Since the algebra $\U$ is generated by $E,F,K^{\pm 1}$ the uniqueness follows. 
\end{proof}

\begin{defn}
\label{defn:L(Omega)_Umodule}
We denote the $\U$-module shown in Proposition \ref{prop:L(Omega)_Umodule} by $\C^{\L(\Omega)}(\lambda)$.
We remark that the $\U$-module $\C^{\L(\Omega)}(q)$ was mentioned in \cite[Section 33]{Askeyscheme}. 
\end{defn}

Recall the Casimir element $\Lambda$ of $\U$ from (\ref{Lambda}).

\begin{lem}
\label{lem:Lambda_L(Omega)}
Suppose that $\lambda$ is a nonzero scalar in $\C$. Then the action of $\Lambda$ on the $\U$-module $\C^{\L(\Omega)}(\lambda)$ is as follows: 
\begin{align*}
\Lambda x
&=
(q^{D-2\dim x+1}+q^{2\dim x-D+1}+q^{-1-D}-q^{1-D}) x
\\
&\qquad +\;
q^{1-D}(q-q^{-1})^2\sum_{\substack{x'\in \L_{\dim x}(\Omega)
\\
x\cap x' \subsetdot x}}
x'
\qquad 
\hbox{for all $x\in \L(\Omega)$}.
\end{align*}
\end{lem}
\begin{proof}
It is straightforward to evaluate the action of $\Lambda$ on the $\U$-module $\C^{\L(\Omega)}(\lambda)$ by applying Proposition \ref{prop:L(Omega)_Umodule}. 
\end{proof}

We are now devoted to developing a $q$-analog of \cite[Lemma 5.4]{Huang:CG&Johnson}.

\begin{defn}
\label{defn:D12}
\begin{enumerate}
\item Define $D_1(x_0)$ to be the linear map $\C^{\L(\Omega)}\to \C^{\L(\Omega)}$ given by 
\begin{eqnarray*}
x &\mapsto &
q^{\dim{\Omega/x_0}-2\dim (x+x_0/x_0)} 
x
\qquad 
\hbox{for all $x\in \L(\Omega)$}.
\end{eqnarray*}

\item Define $D_2(x_0)$ to be the linear map $\C^{\L(\Omega)}\to \C^{\L(\Omega)}$ given by 
\begin{eqnarray*}
x &\mapsto &
q^{\dim x_0-2 \dim x\cap x_0} x
\qquad 
\hbox{for all $x\in \L(\Omega)$}.
\end{eqnarray*}
\end{enumerate}
Note that $D_1(x_0)$ and $D_2(x_0)$ are invertible.
\end{defn}

\begin{lem}
\label{lem:x+x0/x0}
$\dim(x/x\cap x_0)=\dim(x+x_0/x_0)$ for all $x\in \L(\Omega)$.
\end{lem}
\begin{proof}
Let $x\in \L(\Omega)$ be given.
Consider the linear map $x\to x+x_0/x_0$ given by $u\mapsto u+x_0$ for all $u\in x$. Since the linear map is onto and its kernel is equal to $x\cap x_0$ the lemma follows. 
\end{proof}

\begin{lem}
\label{lem:L1+L2}
Suppose that $\lambda$ is a nonzero scalar in $\C$. Then the following equations hold on the $\U$-module $\C^{\L(\Omega)}(\lambda)$:
\begin{align*}
E&=\lambda q^{-D}(L_1(x_0)+L_2(x_0)),
\\
F&=\lambda^{-1} q(R_1(x_0)+R_2(x_0)),
\\
K^{\pm 1}&=D_1(x_0)^{\pm 1}\circ D_2(x_0)^{\pm 1}=D_2(x_0)^{\pm 1}\circ D_1(x_0)^{\pm 1}.
\end{align*}
\end{lem}
\begin{proof}
By Definition \ref{defn:L12&R12}(i), (ii) and Lemma \ref{lem:x+x0/x0} the linear map $L_1(x_0)+L_2(x_0)$ sends 
\begin{eqnarray*}
x &\mapsto & \sum_{x'\subsetdot x} x'
\qquad 
\hbox{for all $x\in \L(\Omega)$}.
\end{eqnarray*}
By Definition \ref{defn:L12&R12}(iii), (iv) and Lemma \ref{lem:x+x0/x0} the linear map $R_1(x_0)+R_2(x_0)$ sends 
\begin{eqnarray*}
x &\mapsto & \sum_{x\subsetdot x'} x'
\qquad 
\hbox{for all $x\in \L(\Omega)$}.
\end{eqnarray*}
By Definition \ref{defn:D12} and Lemma \ref{lem:x+x0/x0} the linear map $D_1(x_0)\circ D_2(x_0)=D_2(x_0)\circ D_1(x_0)$ sends 
\begin{eqnarray*}
x &\mapsto & q^{D-2\dim x} x
\qquad 
\hbox{for all $x\in \L(\Omega)$}.
\end{eqnarray*}
The lemma follows by comparing the above results with Proposition \ref{prop:L(Omega)_Umodule}.
\end{proof}

\begin{defn}
\label{defn:iota}
Define $\iota(x_0)$ to be the linear map $\C^{\L(\Omega)}\to \C^{\L(\Omega/x_0)}\otimes \C^{\L(x_0)}$ given by 
\begin{eqnarray*}
x &\mapsto & (x+x_0)/x_0 \otimes x\cap x_0
\qquad 
\hbox{for all $x\in \L(\Omega)$}.
\end{eqnarray*}
\end{defn}

\begin{lem}
\label{lem:E1diagram}
Suppose that $\lambda$ is a nonzero scalar in $\C$. Then the following diagram commutes:
\begin{table}[H]
\centering
\begin{tikzpicture}
\matrix(m)[matrix of math nodes,
row sep=4em, column sep=4em,
text height=1.5ex, text depth=0.25ex]
{
\C^{\L(\Omega)}
&\C^{\L(\Omega/ x_0)}(1)\otimes \C^{\L(x_0)}(\lambda)\\
\C^{\L(\Omega)}
&\C^{\L(\Omega/ x_0)}(1)\otimes \C^{\L(x_0)}(\lambda)\\
};
\path[->,font=\scriptsize,>=angle 90]
(m-1-1) edge node[left] {$q^{\dim x_0-D}L_1(x_0)$} (m-2-1)
(m-1-1) edge node[above] {$\iota(x_0)$} (m-1-2)
(m-2-1) edge node[below] {$\iota(x_0)$} (m-2-2)
(m-1-2) edge node[right] {$E\otimes 1$} (m-2-2);
\end{tikzpicture}
\end{table}
\end{lem}
\begin{proof}
Let $x\in \L(\Omega)$ be given. Given any $y\in \L(\Omega/x_0)$ with $y\subsetdot x+x_0/x_0$ there is exactly one linear map $\tau:y\to x_0/x\cap x_0$ such that $\tau(u)\subseteq \t(x)(u)$ for all $u\in y$. It follows from Definitions \ref{defn:L12&R12}(i) and \ref{defn:iota} that $\iota(x_0)\circ L_1(x_0)$ sends
\begin{eqnarray*}
x 
&\mapsto &
\left(
\sum_{y\subsetdot (x+x_0)/x_0}  y
\right) 
\otimes x\cap x_0.
\end{eqnarray*}
By Proposition \ref{prop:L(Omega)_Umodule} the element $E$ sends 
\begin{eqnarray*}
x+x_0/x_0
&\mapsto &
q^{\dim x_0-D}
\sum_{y\subsetdot (x+x_0)/x_0}  y
\end{eqnarray*}
on the $\U$-module $\C^{\L(\Omega/ x_0)}(1)$. The lemma follows.
\end{proof}

\begin{lem}
\label{lem:E2diagram}
Suppose that $\lambda$ is a nonzero scalar in $\C$. Then the following diagram commutes:
\begin{table}[H]
\centering
\begin{tikzpicture}
\matrix(m)[matrix of math nodes,
row sep=4em, column sep=4em,
text height=1.5ex, text depth=0.25ex]
{
\C^{\L(\Omega)}
&\C^{\L(\Omega/ x_0)}(\lambda)\otimes \C^{\L(x_0)}(q^{\dim x_0})\\
\C^{\L(\Omega)}
&\C^{\L(\Omega/ x_0)}(\lambda)\otimes \C^{\L(x_0)}(q^{\dim x_0})\\
};
\path[->,font=\scriptsize,>=angle 90]
(m-1-1) edge node[left] {$q^{\dim x_0-D} D_1(x_0)\circ L_2(x_0)$} (m-2-1)
(m-1-1) edge node[above] {$\iota(x_0)$} (m-1-2)
(m-2-1) edge node[below] {$\iota(x_0)$} (m-2-2)
(m-1-2) edge node[right] {$1\otimes E$} (m-2-2);
\end{tikzpicture}
\end{table}
\end{lem}
\begin{proof}
Let $x\in \L(\Omega)$ be given. Given any $z\in \L(x_0)$ with $z\subsetdot x\cap x_0$  there are exactly 
$$
|\F|^{\dim (x+x_0/x_0)}
=q^{2\dim (x+x_0/x_0)}
$$ 
linear maps $\tau:x+x_0/x_0\to x_0/z$ such that $\tau(u)\subseteq \t(x)(u)$ for all $u\in x+x_0/x_0$. It follows from Definitions \ref{defn:L12&R12}(ii), \ref{defn:D12}(i) and \ref{defn:iota} that $\iota(x_0)\circ D_1(x_0)\circ L_2(x_0)$ sends
\begin{eqnarray*}
x 
&\mapsto &
q^{D-\dim x_0} (x+x_0)/x_0\otimes  
\left(
\sum_{z\subsetdot x\cap x_0} z 
\right).
\end{eqnarray*}
By Proposition \ref{prop:L(Omega)_Umodule} the element $E$ sends 
\begin{eqnarray*}
x\cap x_0
&\mapsto &
\sum_{z\subsetdot x\cap x_0} z
\end{eqnarray*}
on the $\U$-module $\C^{\L(x_0)}(q^{\dim x_0})$.
The lemma follows.
\end{proof}

\begin{lem}
\label{lem:F1diagram}
Suppose that $\lambda$ is a nonzero scalar in $\C$. Then the following diagram commutes:
\begin{table}[H]
\centering
\begin{tikzpicture}
\matrix(m)[matrix of math nodes,
row sep=4em, column sep=4em,
text height=1.5ex, text depth=0.25ex]
{
\C^{\L(\Omega)}
&\C^{\L(\Omega/ x_0)}(1)\otimes \C^{\L(x_0)}(\lambda)\\
\C^{\L(\Omega)}
&\C^{\L(\Omega/ x_0)}(1)\otimes \C^{\L(x_0)}(\lambda)\\
};
\path[->,font=\scriptsize,>=angle 90]
(m-1-1) edge node[left] {$q^{1-\dim x_0} R_1(x_0)\circ D_2(x_0)^{-1}$} (m-2-1)
(m-1-1) edge node[above] {$\iota(x_0)$} (m-1-2)
(m-2-1) edge node[below] {$\iota(x_0)$} (m-2-2)
(m-1-2) edge node[right] {$F\otimes 1$} (m-2-2);
\end{tikzpicture}
\end{table}
\end{lem}
\begin{proof}
Let $x\in \L(\Omega)$ be given. Given any $y\in \L(\Omega/x_0)$ with $x+x_0/x_0\subsetdot y$ there are exactly 
$$
|x_0/x\cap x_0|
=|\F|^{\dim (x_0/x\cap x_0)}
=q^{2(\dim x_0-\dim x\cap x_0)}
$$ 
linear maps $\tau:y\to x_0/ x\cap x_0$ such that $\t(x)(u)\subseteq \tau(u)$ for all $u\in x+x_0/x_0$. It follows from Definitions \ref{defn:L12&R12}(iii), \ref{defn:D12}(ii) and \ref{defn:iota} that $\iota(x_0)\circ R_1(x_0)\circ D_2(x_0)^{-1}$ sends
\begin{eqnarray*}
x 
&\mapsto &
q^{\dim x_0} 
\left(
\sum_{(x+x_0)/x_0\subsetdot y} y
\right)
\otimes x\cap x_0.
\end{eqnarray*}
By Proposition \ref{prop:L(Omega)_Umodule} the element $F$ sends 
\begin{eqnarray*}
x+x_0/x_0
&\mapsto &
q \sum_{(x+x_0)/x_0\subsetdot y} y
\end{eqnarray*}
on the $\U$-module $\C^{L(\Omega/x_0)}(1)$. The lemma follows.
\end{proof}

\begin{lem}
\label{lem:F2diagram}
Suppose that $\lambda$ is a nonzero scalar in $\C$. Then the following diagram commutes:
\begin{table}[H]
\centering
\begin{tikzpicture}
\matrix(m)[matrix of math nodes,
row sep=4em, column sep=4em,
text height=1.5ex, text depth=0.25ex]
{
\C^{\L(\Omega)}
&\C^{\L(\Omega/ x_0)}(\lambda)\otimes \C^{\L(x_0)}(q^{\dim x_0})\\
\C^{\L(\Omega)}
&\C^{\L(\Omega/ x_0)}(\lambda)\otimes \C^{\L(x_0)}(q^{\dim x_0})\\
};
\path[->,font=\scriptsize,>=angle 90]
(m-1-1) edge node[left] {$q^{1-\dim x_0}  R_2(x_0)$} (m-2-1)
(m-1-1) edge node[above] {$\iota(x_0)$} (m-1-2)
(m-2-1) edge node[below] {$\iota(x_0)$} (m-2-2)
(m-1-2) edge node[right] {$1\otimes F$} (m-2-2);
\end{tikzpicture}
\end{table}
\end{lem}
\begin{proof}
Let $x\in \L(\Omega)$ be given. Given any $z\in \L(x_0)$ with $x\cap x_0 \subsetdot z$ there is exactly one linear map $\tau:x+x_0/x_0\to x_0/z$ satisfying $\t(x)(u)\subseteq \tau(u)$ for all $u\in x+x_0/x_0$. It follows from Definitions \ref{defn:L12&R12}(iv) and \ref{defn:iota} that $\iota(x_0)\circ R_2(x_0)$ sends 
\begin{eqnarray*}
x
&\mapsto &
(x+x_0)/x_0\otimes 
\left( 
\sum_{x\cap x_0\subsetdot z}  z
\right).
\end{eqnarray*}
By Proposition \ref{prop:L(Omega)_Umodule} the element $F$ sends  
\begin{eqnarray*}
x\cap x_0
&\mapsto &
q^{1-\dim x_0} 
\sum_{x\cap x_0\subsetdot z}  z
\end{eqnarray*}
on the $\U$-module $\C^{\L(x_0)}(q^{\dim x_0})$. 
The lemma follows.
\end{proof}

\begin{lem}
\label{lem:K1diagram}
Suppose that $\lambda$ and $\mu$ are two nonzero scalars in $\C$. Then the following diagram commutes:
\begin{table}[H]
\centering
\begin{tikzpicture}
\matrix(m)[matrix of math nodes,
row sep=4em, column sep=4em,
text height=1.5ex, text depth=0.25ex]
{
\C^{\L(\Omega)}
&\C^{\L(\Omega/ x_0)}(\lambda)\otimes \C^{\L(x_0)}(\mu)\\
\C^{\L(\Omega)}
&\C^{\L(\Omega/ x_0)}(\lambda)\otimes \C^{\L(x_0)}(\mu)\\
};
\path[->,font=\scriptsize,>=angle 90]
(m-1-1) edge node[left] {$D_1(x_0)$} (m-2-1)
(m-1-1) edge node[above] {$\iota(x_0)$} (m-1-2)
(m-2-1) edge node[below] {$\iota(x_0)$} (m-2-2)
(m-1-2) edge node[right] {$K\otimes 1$} (m-2-2);
\end{tikzpicture}
\end{table}
\end{lem}
\begin{proof}
Immediate from Proposition \ref{prop:L(Omega)_Umodule} and Definition \ref{defn:D12}(i).
\end{proof}

\begin{lem}
\label{lem:K2diagram}
Suppose that $\lambda$ and $\mu$ are two nonzero scalars in $\C$. Then the following diagram commutes:
\begin{table}[H]
\centering
\begin{tikzpicture}
\matrix(m)[matrix of math nodes,
row sep=4em, column sep=4em,
text height=1.5ex, text depth=0.25ex]
{
\C^{\L(\Omega)}
&\C^{\L(\Omega/ x_0)}(\lambda)\otimes \C^{\L(x_0)}(\mu)\\
\C^{\L(\Omega)}
&\C^{\L(\Omega/ x_0)}(\lambda)\otimes \C^{\L(x_0)}(\mu)\\
};
\path[->,font=\scriptsize,>=angle 90]
(m-1-1) edge node[left] {$D_2(x_0)$} (m-2-1)
(m-1-1) edge node[above] {$\iota(x_0)$} (m-1-2)
(m-2-1) edge node[below] {$\iota(x_0)$} (m-2-2)
(m-1-2) edge node[right] {$1\otimes K$} (m-2-2);
\end{tikzpicture}
\end{table}
\end{lem}
\begin{proof}
Immediate from Proposition \ref{prop:L(Omega)_Umodule} and Definition \ref{defn:D12}(ii).
\end{proof}

\begin{thm}
\label{thm:X&DeltaXq}
The following diagram commutes for each $X\in \U$:
\begin{table}[H]
\centering
\begin{tikzpicture}
\matrix(m)[matrix of math nodes,
row sep=4em, column sep=4em,
text height=1.5ex, text depth=0.25ex]
{
\C^{\L(\Omega)}(q^{\dim x_0})
&\C^{\L(\Omega/ x_0)}(1)\otimes \C^{\L(x_0)}(q^{\dim x_0})\\
\C^{\L(\Omega)}(q^{\dim x_0})
&\C^{\L(\Omega/ x_0)}(1)\otimes \C^{\L(x_0)}(q^{\dim x_0})\\
};
\path[->,font=\scriptsize,>=angle 90]
(m-1-1) edge node[left] {$X$} (m-2-1)
(m-1-1) edge node[above] {$\iota(x_0)$} (m-1-2)
(m-2-1) edge node[below] {$\iota(x_0)$} (m-2-2)
(m-1-2) edge node[right] {$\Delta(X)$} (m-2-2);
\end{tikzpicture}
\end{table}
\end{thm}
\begin{proof}
Recall from Lemma \ref{lem:U_comultiplication} that $\Delta(E)=E\otimes 1+K^{-1}\otimes E$. 
By Lemma \ref{lem:L1+L2} the equation $E=q^{\dim x_0-D}(L_1(x_0)+L_2(x_0))$ holds on the $\U$-module $\C^{\L(\Omega)}(q^{\dim x_0})$. 
Combined with Lemmas \ref{lem:E1diagram},  \ref{lem:E2diagram} and \ref{lem:K1diagram} the above diagram commutes for $X=E$. 
Recall from Lemma \ref{lem:U_comultiplication} that $\Delta(F)=F\otimes K+1\otimes F$. 
By Lemma \ref{lem:L1+L2} the equation $F=q^{1-\dim x_0}(R_1(x_0)+R_2(x_0))$ holds on the $\U$-module $\C^{\L(\Omega)}(q^{\dim x_0})$. 
By Lemmas \ref{lem:F1diagram} , \ref{lem:F2diagram} and \ref{lem:K2diagram} the above diagram commutes for $X=F$. 
Recall from Lemma \ref{lem:U_comultiplication} that $\Delta(K^{\pm 1})=K^{\pm 1}\otimes K^{\pm 1}$. By Lemma \ref{lem:L1+L2} the equation $K^{\pm 1}=D_1(x_0)^{\pm 1}\circ D_2(x_0)^{\pm 1}$ holds on the $\U$-module $\C^{\L(\Omega)}(q^{\dim x_0})$. By Lemmas \ref{lem:K1diagram} and \ref{lem:K2diagram} the above diagram commutes for $X=K^{\pm 1}$. Since the algebra $\U$ is generated by $E,F,K^{\pm 1}$ the result follows.
\end{proof}

Theorem \ref{thm:X&DeltaXq} is a $q$-analog of \cite[Lemma 5.4]{Huang:CG&Johnson}. It is disappointing that $\iota(x_0)$ is not a linear isomorphism in general by Theorem \ref{thm:L(Omega)->Lx0(Omega)}.
We end this section with some trivial relations among the linear maps $L_1(x_0), L_2(x_0), R_1(x_0), R_2(x_0), D_1(x_0), D_2(x_0)$.

\begin{lem}
[Lemma 7.4, \cite{Watanabe:2017}]
\label{lem:L12R12D12}
The following equations hold:
\begin{multicols}{2}
\begin{enumerate}
\item $[L_1(x_0),D_2(x_0)]=0$.  

\item $[L_2(x_0),D_1(x_0)]=0$.

\item $[R_1(x_0),D_2(x_0)]=0$.

\item $[R_2(x_0),D_1(x_0)]=0$.

\item $[L_1(x_0),D_1(x_0)]_q=0$.

\item $[L_2(x_0),D_2(x_0)]_q=0$.

\item $[D_1(x_0),R_1(x_0)]_q=0$.

\item $[D_2(x_0),R_2(x_0)]_q=0$.
\end{enumerate}
\end{multicols}
\end{lem}
\begin{proof}
It is straightforward to verify these equations by using Definitions \ref{defn:L12&R12} and \ref{defn:D12}.
\end{proof}

\begin{lem}
[Lemma 7.5, \cite{Watanabe:2017}]
\label{lem:L1L2R1R2}
The following equations hold:
\begin{multicols}{2}
\begin{enumerate}
\item $[L_2(x_0),L_1(x_0)]_q=0$.

\item $[R_1(x_0),R_2(x_0)]_q=0$.

\item $[L_1(x_0),R_2(x_0)]=0$.

\item $[L_2(x_0),R_1(x_0)]=0$.
\end{enumerate}
\end{multicols}
\end{lem}
\begin{proof}
(i): Let $x\in \L(\Omega)$ be given. 
By Definition \ref{defn:L12&R12}(i), (ii) both of $(L_2(x_0)\circ L_1(x_0))(x)$ and $(L_1(x_0)\circ L_2(x_0))(x)$ are the linear combinations of 
$$
S=\{x'\in \L(\Omega)\,|\, 
x'+x_0/x_0\subsetdot x+x_0/x_0,\, 
x'\cap x_0\subsetdot x\cap x_0\}.
$$
Fix an element $x'\in S$. 
It is immediate from Definition \ref{defn:taux} that $\t(x')(u)\subseteq \t(x)(u)$ for all $u\in x'+x_0/x_0$.  
There is exactly one linear map $\sigma:x'+x_0/x_0\to x_0/x\cap x_0$ such that 
$\t(x')(u)\subseteq \sigma(u)\subseteq \t(x)(u)$ for all $u\in x'+x_0/x_0$. Hence the coefficient of $x'$ in $(L_2(x_0)\circ L_1(x_0))(x)$ is equal to one. 
There are exactly 
$$
|x\cap x_0/x'\cap x_0|=q^2
$$
linear maps $\sigma:x+x_0/x_0\to x_0/x'\cap x_0$ such that $\t(x')(u)\subseteq \sigma(u)$ for all $u\in x'+x_0/x_0$ and $\sigma(u)\subseteq \t(x)(u)$ for all $u\in x+x_0/x_0$. Hence the coefficient of $x'$ in $(L_1(x_0)\circ L_2(x_0))(x)$ is equal to $q^2$. 
The equation (i) follows.

(ii): Let $x\in \L(\Omega)$ be given. 
By Definition \ref{defn:L12&R12}(iii), (iv) both of $(R_1(x_0)\circ R_2(x_0))(x)$ and $(R_2(x_0)\circ R_1(x_0))(x)$ are the linear combinations of 
$$
S=\{x'\in \L(\Omega)\,|\, 
x+x_0/x_0\subsetdot x'+x_0/x_0,\, 
x\cap x_0\subsetdot x'\cap x_0\}.
$$
Fix an element $x'\in S$. 
It is immediate from Definition \ref{defn:taux} that $\t(x)(u)\subseteq \t(x')(u)$ for all $u\in x+x_0/x_0$. 
There is exactly one linear map $\sigma:x+x_0/x_0\to x_0/x'\cap x_0$ such that $\t(x)(u)\subseteq \sigma(u)\subseteq \t(x')(u)$ for all $u\in x+x_0/x_0$. Hence the coefficient of $x'$ in $(R_1(x_0)\circ R_2(x_0))(x)$ is equal to one. 
There are exactly 
$$
|x'\cap x_0/x\cap x_0|=q^2
$$ 
linear maps $\sigma:x'+x_0/x_0\to x_0/x\cap x_0$ such that $\t(x)(u)\subseteq \sigma(u)$ for all $u\in x+x_0/x_0$ and $\sigma(u) \subseteq \t(x')(u)$ for all $u\in x'+x_0/x_0$. Hence the coefficient of $x'$ in $(R_2(x_0)\circ R_1(x_0))(x)$ is equal to $q^2$.
The equation (ii) follows.

(iii): Let $x\in \L(\Omega)$ be given. 
Set $S$ to be the set consisting of all $(y,z,\tau)\in \L(\Omega)_{x_0}$ with $y\subsetdot x+x_0/x_0$, $x\cap x_0\subsetdot z$ and $\t(x)(u)\subseteq \tau(u)$ for all $u\in y$. 
By Definition \ref{defn:L12&R12}(i), (iv) both of $(L_1(x_0)\circ R_2(x_0))(x)$ and $(R_2(x_0)\circ L_1(x_0))(x)$ are the linear combinations of $S$.
Fix an element $(y,z,\tau)\in S$. 
There is exactly one linear map $\sigma:x+x_0/x_0\to x_0/z$ such that $\t(x)(u)\subseteq \sigma(u)$ for all $u\in x+x_0/x_0$ and $\tau(u)\subseteq \sigma(u)$ for all $u\in y$. Hence the coefficient of $(y,z,\tau)$ in $(L_1(x_0)\circ R_2(x_0))(x)$ is equal to one. 
There is exactly one linear map $\sigma:y\to x_0/x\cap x_0$ such that $\sigma(u)\subseteq \t(x)(u)$ and $\sigma(u)\subseteq \tau(u)$ for all $u\in y$. Hence the coefficient of $(y,z,\tau)$ in $(R_2(x_0)\circ L_1(x_0))(x)$ is equal to one. 
The equation (iii) follows.

(iv): Let $x\in \L(\Omega)$ be given. 
Set $S$ to be the set consisting of all $(y,z,\tau)\in \L(\Omega)_{x_0}$ with $x+x_0/x_0\subsetdot y$, $z\subsetdot x\cap x_0$ and $\tau(u)\subseteq \t(x)(u)$ for all $u\in x+x_0/x_0$. 
By Definition \ref{defn:L12&R12}(ii), (iii) both of $(L_2(x_0)\circ R_1(x_0))(x)$ and $(R_1(x_0)\circ L_2(x_0))(x)$ are the linear combinations of $S$. 
Fix an element $(y,z,\tau)\in S$. 
There is exactly one linear map $\sigma:x+x_0/x_0\to x_0/z$ such that $\sigma(u)\subseteq\t(x)(u)$ and $\sigma(u)\subseteq \tau(u)$ for all $u\in x+x_0/x_0$. Hence the coefficient of $(y,z,\tau)$ in $(R_1(x_0)\circ L_2(x_0))(x)$ is equal to one.
There is exactly one linear map $\sigma:y\to x_0/x\cap x_0$ such that $\t(x)(u)\subseteq \sigma(u)$ for all $u\in x+x_0/x_0$ and $\tau(u)\subseteq \sigma(u)$ for all $u\in y$.
Hence the coefficient of $(y,z,\tau)$ in $(L_2(x_0)\circ R_1(x_0))(x)$ is equal to one.
The equation (iv) follows.
\end{proof}

\section{The poset $(\L(\Omega)_{x_0},\subseteq)$ and the bilinear forms graphs}\label{s:bilinear}

In this section we study a natural connection between $(\L(\Omega)_{x_0},\subseteq)$ and the bilinear forms graphs. 
For any linear map $\tau$, the notations $\ker(\tau)$ and ${\rm Im}(\tau)$ stand for the kernel and image of $\tau$ respectively, and  ${\rm rk}(\tau)$ denotes the rank of $\tau$ provided that ${\rm Im}(\tau)$ is finite-dimensional.

\begin{lem}\label{lem1:bilinear}
Suppose that $x,x'\in \L(\Omega)$ with $x+x_0/x_0=x'+x_0/x_0$ and $x\cap x_0=x'\cap x_0$. 
Then 
$$
\ker\left(\t(x)-\t(x')\right)=(x\cap x')+x_0/x_0.
$$
\end{lem}
\begin{proof}
($\supseteq$): Suppose that $u+x_0\in (x\cap x')+x_0/x_0$ where $u\in x\cap x'$. Recall from Definition \ref{defn:taux} that $\t(x)(u+x_0)=u_0+(x\cap x_0)$ and $\t(x')(u+x_0)=u_0+(x'\cap x_0)$. 
Since $x\cap x_0=x'\cap x_0$ it follows that $\t(x)(u+x_0)=\t(x')(u+x_0)$. Therefore $u+x_0\in \ker\left(\t(x)-\t(x')\right)$. 

($\subseteq$): Suppose that $u+x_0\in \ker\left(\t(x)-\t(x')\right)$ where $u\in x$. Since $x+x_0/x_0=x'+x_0/x_0$ there exists a vector $u'\in x'$ such that $u+x_0=u'+x_0$. Hence $u-u'=u_0-u_0'$. Since $u+x_0=u'+x_0$ lies in $\ker\left(\t(x)-\t(x')\right)$ it follows that $u_0-u_0'$ lies in $x\cap x_0=x'\cap x_0$. Therefore $u\in x'$ and this implies that $u+x_0\in (x\cap x')+x_0/x_0$. 
\end{proof}

\begin{lem}\label{lem2:bilinear}
Suppose that $x,x'\in \L(\Omega)$ with $x+x_0/x_0=x'+x_0/x_0$ and $x\cap x_0=x'\cap x_0$. Then 
\begin{align*}
{\rm rk} \left( \t(x)-\t(x') \right)
&=\dim (x/x\cap x')=\dim (x'/x\cap x')
\\
&=\dim (x+x'/x)=\dim(x+x'/x').
\end{align*}
\end{lem}
\begin{proof}
Since $x+x_0/x_0=x'+x_0/x_0$ it follows from Lemma \ref{lem:x+x0/x0} that $\dim(x/x\cap x_0)=\dim(x'/x'\cap x_0)$.
Since $x\cap x_0=x'\cap x_0$ it follows that $\dim x=\dim x'$. Thus it suffices to show the first equality. 
Applying the rank-nullity theorem to the linear map $\t(x)-\t(x')$
 and by Lemma \ref{lem1:bilinear} it follows that ${\rm rk}(\t(x)-\t(x'))$ is equal to 
$$
\dim (x+x_0/x_0)-\dim ((x\cap x')+x_0/x_0).
$$ 
By Lemma \ref{lem:x+x0/x0} the above number is equal to 
$\dim (x/x\cap x_0)-\dim (x\cap x'/x\cap x'\cap x_0)$. 
Since $x\cap x'\cap x_0=x\cap x_0$ the first equality holds. The lemma follows.
\end{proof}

\begin{defn}
\label{defn:D3}
Define $D_3(x_0)$ to be the linear map $\C^{\L(\Omega)}\to \C^{\L(\Omega)}$ given by 
\begin{eqnarray*}
x &\mapsto & 
\!\!\!\!
\sum_{
\substack{
(x+x_0/x_0,x\cap x_0,\tau)\in \L(\Omega)_{x_0}
\\
{\rm rk}\left(\t(x)-\tau\right)=1}} 
\!\!\!\!
(x+x_0/x_0,x\cap x_0,\tau)
\qquad 
\hbox{for all $x\in \L(\Omega)$}.
\end{eqnarray*}
By Theorem \ref{thm:L(Omega)->Lx0(Omega)} and Lemma \ref{lem2:bilinear} the aforementioned linear map sends
\begin{eqnarray*}
x &\mapsto & 
\!\!\!\!
\sum_{
\substack{
x\cap x'\subsetdot x
\\
x+x_0/x_0=x'+x_0/x_0
\\
x\cap x_0=x'\cap x_0}} 
\!\!\!\!
x'
\qquad 
\hbox{for all $x\in \L(\Omega)$}.
\end{eqnarray*}
Therefore $D_3(x_0)$ is indeed independent of the choice of $x_1$.
\end{defn}

\begin{defn}
\label{defn:Cyz}
For any $y\in \L(\Omega/x_0)$ and $z\in \L(x_0)$ let $\C^{\L(\Omega)}_{y,z}$ denote the subspace of $\C^{\L(\Omega)}$ spanned by all $x\in \L(\Omega)$ with $x+x_0/x_0=y$ and $x\cap x_0=z$. Note that 
$$
\C^{\L(\Omega)}=
\bigoplus\limits_{y\in \L(\Omega/x_0)} 
\bigoplus\limits_{z\in \L(x_0)}
\C^{\L(\Omega)}_{y,z}.
$$
\end{defn}

\begin{lem}
\label{lem:Cyz}
Suppose that $y\in \L(\Omega/x_0)$ and $z\in \L(x_0)$.  
Then 
$$
\dim  \C^{\L(\Omega)}_{y,z}=q^{2 \dim y \dim (x_0/z)}.
$$ 
\end{lem}
\begin{proof}
There are exactly 
$
|x_0/z|^{\dim y}=q^{2\dim y \dim (x_0/z)}
$
linear maps from $y$ into $x_0/z$. Combined with Theorem \ref{thm:L(Omega)->Lx0(Omega)} the lemma follows.
\end{proof}

Let $y\in \L(\Omega/x_0)$ and $z\in \L(x_0)$ be given. 
Let ${\rm Bil}(y,x_0/z)$ denote the {\it bilinear forms graph} whose vertices are all linear maps $y\to x_0/z$ and two vertices $\sigma$ and $\tau$ are adjacent whenever ${\rm rk}(\sigma-\tau)=1$. 
Observe that $\C^{\L(\Omega)}_{y,z}$ is $D_3(x_0)$-invariant and $D_3(x_0)|_{\C^{\L(\Omega)}_{y,z}}$ is the adjacency operator of ${\rm Bil}(y,x_0/z)$.

\begin{lem}
\label{lem:D3}
Suppose that $y\in \L(\Omega/x_0)$ and $z\in \L(x_0)$.  
Then the eigenvalues of $D_3(x_0)$ in $\C^{\L(\Omega)}_{y,z}$ are 
\begin{gather}\label{D3:eign}
q^{\dim y-1}(q^{2\dim (x_0/z)-h}[\dim y-h]_q-[\dim y]_q)
\end{gather}
for all integers $h$ with $0\leq h\leq \min\{\dim y,\dim (x_0/z)\}$ and their multiplicities are 
\begin{gather}\label{D3:mult}
q^{h(\dim y-1)}
{\dim y \brack h}_q
\prod_{i=0}^{h-1}
(q^{2(\dim(x_0/z)-i)}-1)
\end{gather}
for all integers $h$ with $0\leq h\leq \min\{\dim y,\dim (x_0/z)\}$, respectively.
\end{lem}
\begin{proof}
According to the first eigenmatrix of ${\rm Bil}(y,x_0/z)$ given by \cite{BannaiIto1984,attenuated2013,attenuated2022}, the eigenvalues of $D_3(x_0)|_{\C^{\L(\Omega)}_{y,z}}$ are 
\begin{gather}\label{D3:eign'}
q^{2\dim(x_0/z)}|\L_1(\F^{\dim y-h})|
-|\L_1(\F^{\dim y})|
\end{gather}
for all integers $h$ with $0\leq h\leq \min\{\dim y,\dim (x_0/z)\}$.
By Lemma \ref{lem:|Lk|} the value (\ref{D3:eign'}) is equal to (\ref{D3:eign}).

The bilinear forms schemes are self-dual \cite{BannaiIto1984,attenuated2013}.
According to the first eigenmatrix of ${\rm Bil}(y,x_0/z)$ given by  \cite{BannaiIto1984,attenuated2013,attenuated2022}, the eigenvalue (\ref{D3:eign}) of $D_3(x_0)|_{\C^{\L(\Omega)}_{y,z}}$ is of multiplicity 
\begin{gather}\label{D3:mult'}
\sum_{i=0}^h
(-1)^{i-h}
q^{2i\dim(x_0/z)+(i-h)(i-h+1)}
|\L_{h-i}(\F^{\dim y-i})|
|\L_i(\F^{\dim y})|.
\end{gather}
Using (\ref{qbino}) and Lemma \ref{lem:|Lk|} yields that (\ref{D3:mult'}) is equal to $q^{h(\dim y-1)}{\dim y\brack h}_q$ times 
\begin{gather*}
\sum_{i=0}^h
(-1)^{i-h}
q^{i(2\dim(x_0/z)-h+1)}
{h\brack i}_q.
\end{gather*}
Evaluating the above sum by using the $q$-binomial theorem yields that (\ref{D3:mult'}) is equal to (\ref{D3:mult}). The lemma follows.
\end{proof}

In Lemmas \ref{lem:L1&D34}--\ref{lem:R2&D34} we will see how $D_3(x_0)$ is related to the following map:

\begin{defn}
\label{defn:D4}
Define $D_4(x_0)$ to be the linear map $\C^{\L(\Omega)}\to \C^{\L(\Omega)}$ given by 
\begin{eqnarray*}
x &\mapsto & 
\frac{|x\cup x_0|}{|x\cap x_0|} x
=\left(
\frac{|x|+|x_0|}{|x\cap x_0|}-1
\right)x
\qquad 
\hbox{for all $x\in \L(\Omega)$}.
\end{eqnarray*}
\end{defn}

\begin{lem}
\label{lem:D1234}
The linear maps $D_1(x_0)$, $D_2(x_0)$, $D_3(x_0)$, $D_4(x_0)$ mutually commute.
\end{lem}
\begin{proof}
It is routine to verify the lemma by using Definitions \ref{defn:D12}, \ref{defn:D3} and \ref{defn:D4}.
\end{proof}

\begin{lem}
\label{lem:D4}
Suppose that $y\in \L(\Omega/x_0)$ and $z\in \L(x_0)$.  
Then $D_4(x_0)$ acts on $\C^{\L(\Omega)}_{y,z}$ as scalar multiplication by 
$q^{2\dim y}+q^{2\dim (x_0/z)}-1$.
\end{lem}
\begin{proof}
Pick any $x\in \C^{\L(\Omega)}_{y,z}$. 
By Definition \ref{defn:D4} the linear map $D_4(x_0)$ maps $x$ to the scalar multiple of $x$ by $|x/z|+|x_0/z|-1$. By Lemma \ref{lem:x+x0/x0} the number $|x/z|$ is equal to $|y|$. The lemma follows.
\end{proof}

\begin{lem}
\label{lem:L1&D34}
The following equations hold:
\begin{enumerate}
\item $[D_3(x_0), L_1(x_0)]_q=q^{-1}(1-q^{\dim x_0} D_2(x_0))\circ L_1(x_0)$.

\item $[D_4(x_0), L_1(x_0)]_q=-(q-q^{-1})(1-q^{\dim x_0} D_2(x_0))\circ L_1(x_0)$.

\end{enumerate}
\end{lem}
\begin{proof}
(i): Let $x\in \L(\Omega)$ be given. Set 
\begin{align*}
S_0&=\{(y,x\cap x_0,\tau)\in \L(\Omega)_{x_0}\,|\, y\subsetdot x+x_0/x_0,{\rm rk}\left(\tau-\t(x)|_y\right)=0\},
\\
S_1&=\{(y,x\cap x_0,\tau)\in \L(\Omega)_{x_0}\,|\, y\subsetdot x+x_0/x_0,{\rm rk}\left(\tau-\t(x)|_y\right)=1\}.
\end{align*}
By Definitions \ref{defn:L12&R12}(i) and \ref{defn:D3} both of $(D_3(x_0)\circ L_1(x_0))(x)$ and $(L_1(x_0)\circ D_3(x_0))(x)$ are the linear combinations of $S_0\cup S_1$. 

Fix an element $(y,x\cap x_0,\tau)\in S_0$.
None of linear maps $\sigma:y\to x_0/x\cap x_0$ satisfies $\sigma(u)\subseteq \t(x)(u)$ for all $u\in y$ and ${\rm rk}(\tau-\sigma)=1$. Hence the coefficient of $(y,x\cap x_0,\tau)$ in $(D_3(x_0)\circ L_1(x_0))(x)$ is zero.
There are exactly 
$$
|x_0/x\cap x_0|-1=q^{2(\dim x_0/x\cap x_0)}-1
$$
linear maps $\sigma:x+x_0/x_0\to x_0/x\cap x_0$ such that ${\rm rk}(\t(x)-\sigma)=1$ and $\tau(u)\subseteq \sigma(u)$ for all $u\in y$. Hence the coefficient of $(y,x\cap x_0,\tau)$ in $(L_1(x_0)\circ D_3(x_0))(x)$ is $q^{2(\dim x_0/x\cap x_0)}-1$.

Fix an element $(y,x\cap x_0,\tau)\in S_1$. Observe that $\t(x)|_y$ is the unique linear map $\sigma:y\to x_0/x\cap x_0$ such that $\sigma(u)\subseteq \t(x)(u)$ for all $u\in y$ and ${\rm rk}(\tau-\sigma)=1$. 
Hence the coefficient of $(y,x\cap x_0,\tau)$ in $(D_3(x_0)\circ L_1(x_0))(x)$ is one.
There are exactly 
$$
\left|
{\rm Im}(\tau-\t(x)|_y)
\right|
=q^2
$$ 
linear maps $\sigma:x+x_0/x_0\to x_0/x\cap x_0$ such that ${\rm rk}(\t(x)-\sigma)=1$ and $\tau(u)\subseteq \sigma(u)$ for all $u\in y$. 
Hence the coefficient of $(y,x\cap x_0,\tau)$ in $(L_1(x_0)\circ D_3(x_0))(x)$ is $q^2$.

By the above comments the $q$-bracket $[D_3(x_0),L_1(x_0)]_q$ sends $x$ to 
$$
q^{-1}(1-q^{2(\dim x_0/x\cap x_0)})\sum_{(y,x\cap x_0,\tau)\in S_0} (y,x\cap x_0,\tau).
$$
By Definitions \ref{defn:L12&R12}(i) and \ref{defn:D12}(ii) the linear map $(1-q^{\dim x_0} D_2(x_0))\circ L_1(x_0)$ sends $x$ to 
\begin{gather}\label{D2L1}
(1-q^{2(\dim x_0/x\cap x_0)})\sum_{(y,x\cap x_0,\tau)\in S_0} (y,x\cap x_0,\tau).
\end{gather}
The equation (i) follows.

(ii): Let $x\in \L(\Omega)$ be given. By Definition \ref{defn:L12&R12}(i) and Lemma \ref{lem:D4} the $q$-bracket $[D_4(x_0), L_1(x_0)]_q$ sends $x$ to the scalar multiple of (\ref{D2L1}) by $-(q-q^{-1})$. The equation (ii) follows.
\end{proof}

For any $z,z'\in \L(x_0)$ with $z\subseteq z'$ let $1_z^{z'}:x_0/z\to x_0/z'$ denote the linear map given by $u+z\mapsto u+z'$ for all $u\in x_0$.

\begin{lem}
\label{lem:L2&D34}
The following equations hold:
\begin{enumerate}
\item $[L_2(x_0),D_3(x_0)]_q=q^{-1}(1-q^{D-\dim x_0} D_1(x_0)^{-1})\circ L_2(x_0)$.

\item $[L_2(x_0),D_4(x_0)]_q=-(q-q^{-1})(1-q^{D-\dim x_0} D_1(x_0)^{-1})\circ L_2(x_0)$.

\end{enumerate}
\end{lem}
\begin{proof}
(i): Let $x\in \L(\Omega)$ be given. 
Set 
\begin{align*}
S_0 &=\{(x+x_0/x_0,z,\tau)\in \L(\Omega)_{x_0}\,|\, z\subsetdot x\cap x_0, {\rm rk}\left(\t(x)-1_z^{x\cap x_0}\circ \tau\right)=0\},
\\
S_1 &=\{(x+x_0/x_0,z,\tau)\in \L(\Omega)_{x_0}\,|\, z\subsetdot x\cap x_0, {\rm rk}\left(\t(x)-1_z^{x\cap x_0}\circ \tau\right)=1\}.
\end{align*}
By Definitions \ref{defn:L12&R12}(ii) and \ref{defn:D3} both of $(L_2(x_0)\circ D_3(x_0))(x)$ and $(D_3(x_0)\circ L_2(x_0))(x)$ are the linear combinations of $S_0\cup S_1$. 

Fix an element $(x+x_0/x_0,z,\tau)\in S_0$. 
None of linear maps $\sigma:x+x_0/x_0\to x_0/x\cap x_0$ satisfies ${\rm rk}(\t(x)-\sigma)=1$ and $\tau(u)\subseteq \sigma(u)$ for all $u\in x+x_0/x_0$. 
Hence the coefficient of $(x+x_0/x_0,z,\tau)$ in $(L_2(x_0)\circ D_3(x_0))(x)$ is zero.
There are exactly 
$$
|x\cap x_0/z|^{\dim (x+x_0/x_0)}-1
=q^{2\dim (x+x_0/x_0)}-1
$$
linear maps $\sigma:x+x_0/x_0\to x_0/z$ such that $\sigma(u)\subseteq \t(x)(u)$ for all $u\in x+x_0/x_0$ and ${\rm rk}(\tau-\sigma)=1$. Hence the coefficient of $(x+x_0/x_0,z,\tau)$ in $(D_3(x_0)\circ L_2(x_0))(x)$ is $q^{2\dim (x+x_0/x_0)}-1$.

Fix an element $(x+x_0/x_0,z,\tau)\in S_1$.  Clearly $1_z^{x\cap x_0}\circ \tau$ is the unique linear map $\sigma:x+x_0/x_0\to x_0/x\cap x_0$ satisfying ${\rm rk}(\t(x)-\sigma)=1$ and $\tau(u)\subseteq \sigma(u)$ for all $u\in x+x_0/x_0$. Hence the coefficient of $(x+x_0/x_0,z,\tau)$ in $(L_2(x_0)\circ D_3(x_0))(x)$ is one. 
Since ${\rm rk}(\t(x)-1_z^{x\cap x_0}\circ\tau)=1$ 
there exists a vector $v\in x_0\setminus(x\cap x_0)$ such that 
$$
{\rm Im}(\t(x)-1_z^{x\cap x_0}\circ\tau)=\F v+(x\cap x_0).
$$ 
Suppose that $\sigma$ is a linear map $x+x_0/x_0\to x_0/z$ satisfying $\sigma(u)\subseteq \t(x)(u)$ for all $u\in x+x_0/x_0$ and ${\rm rk}(\sigma-\tau)=1$. 
Since $1_z^{x\cap x_0}\circ \sigma=\t(x)$ it follows that 
$1_z^{x\cap x_0}$ maps ${\rm Im}(\sigma-\tau)$ onto $\F v+(x\cap x_0)$.
Hence
\begin{align}
\label{e1:L2&D34}
{\rm Im}(\sigma-\tau)&=\F w+z
\end{align}
for some $w\in v+(x\cap x_0)$.
Let $u\in x+x_0/x_0$ be given. By (\ref{e1:L2&D34}) there is a scalar $c\in \F$ such that $\sigma(u)-\tau(u)=cw+z$. 
Applying $1_z^{x\cap x_0}$ to either side of the above equation yields that $c$ is uniquely determined by the following equation:
$$
\t(x)(u)-(1_z^{x\cap x_0}\circ \tau)(u)=c
v+(x\cap x_0).
$$
This shows that $\sigma$ is uniquely determined by the choice of $w+z\subseteq v+(x\cap x_0)$. There are exactly $|x\cap x_0/z|=q^2$ cosets of $z$ in $v+(x\cap x_0)$. 
Hence the coefficient of $(x+x_0/x_0,z,\tau)$ in $(D_3(x_0)\circ L_2(x_0))(x)$ is $q^2$.

By the above comments the $q$-bracket $[L_2(x_0),D_3(x_0)]_q$ sends $x$ to
$$
q^{-1}(1-q^{2\dim (x+x_0/x_0)})\sum_{(x+x_0/x_0,z,\tau)\in S_0} (x+x_0/x_0,z,\tau).
$$
By Definitions \ref{defn:L12&R12}(ii) and \ref{defn:D12}(i) the linear map $(1-q^{D-\dim x_0}D_1(x_0)^{-1})\circ L_2(x_0)$ sends $x$ to 
\begin{gather}\label{D1L2}
(1-q^{2\dim(x+x_0/x_0)})\sum_{(x+x_0/x_0,z,\tau)\in S_0} (x+x_0/x_0,z,\tau).
\end{gather}
The equation (i) follows.

(ii): Let $x\in \L(\Omega)$ be given. By Definition \ref{defn:L12&R12}(ii) and Lemma \ref{lem:D4} the $q$-bracket $[L_2(x_0),D_4(x_0)]_q$ sends $x$ to the scalar multiple of (\ref{D1L2}) by $-(q-q^{-1})$. The equation (ii) follows.
\end{proof}

\begin{lem}
\label{lem:R1&D34}
The following equations hold:
\begin{enumerate}
\item $[R_1(x_0),D_3(x_0)]_q=q^{-1} (1-q^{\dim x_0} D_2(x_0))\circ R_1(x_0)$.

\item $[R_1(x_0),D_4(x_0)]_q=-(q-q^{-1}) (1-q^{\dim x_0} D_2(x_0))\circ R_1(x_0)$.

\end{enumerate}
\end{lem}
\begin{proof}
(i): Let $x\in \L(\Omega)$ be given. Set 
\begin{align*}
S_0 &=\{(y,x\cap x_0,\tau)\in \L(\Omega)_{x_0}\,|\,x+x_0/x_0\subsetdot y, {\rm rk}\left(\t(x)-\tau|_{x+x_0/x_0}\right)=0\},
\\
S_1 &=\{(y,x\cap x_0,\tau)\in \L(\Omega)_{x_0}\,|\,x+x_0/x_0\subsetdot y, {\rm rk}\left(\t(x)-\tau|_{x+x_0/x_0}\right)=1\}.
\end{align*}
By Definitions \ref{defn:L12&R12}(iii) and \ref{defn:D3} both of $(R_1(x_0)\circ D_3(x_0))(x)$ and $(D_3(x_0)\circ R_1(x_0))(x)$ are the linear combinations of $S_0\cup S_1$.

Fix an element $(y,x\cap x_0,\tau)\in S_0$.
None of linear maps $\sigma:x+x_0/x_0\to x_0/x\cap x_0$ satisfies ${\rm rk}(\t(x)-\sigma)=1$ and $\sigma(u)\subseteq \tau(u)$ for all $u\in x+x_0/x_0$. Hence the coefficient of $(y,x\cap x_0,\tau)$ in $(R_1(x_0)\circ D_3(x_0))(x)$ is zero. 
There are exactly 
$$
|x_0/x\cap x_0|-1=q^{2(\dim x_0/x\cap x_0)}-1
$$
linear maps $\sigma:y\to x_0/x\cap x_0$ such that $\t(x)(u)\subseteq \sigma(u)$ for all $u\in x+x_0/x_0$ and ${\rm rk}(\tau-\sigma)=1$. Hence the coefficient of $(y,x\cap x_0,\tau)$ in $(D_3(x_0)\circ R_1(x_0))(x)$ is $q^{2(\dim x_0/x\cap x_0)}-1$.

Fix an element $(y,x\cap x_0,\tau)\in S_1$. Observe that $\tau|_{x+x_0/x_0}$ is the unique linear map $\sigma:x+x_0/x_0\to x_0/x\cap x_0$ such that ${\rm rk}(\t(x)-\sigma)=1$ and $\sigma(u)\subseteq \tau(u)$ for all $u\in x+x_0/x_0$. 
Hence the coefficient of $(y,x\cap x_0,\tau)$ in $(R_1(x_0)\circ D_3(x_0))(x)$ is one. 
There are exactly 
$$
\left|{\rm Im} (\t(x)-\tau|_{x+x_0/x_0})\right|=q^2
$$ 
linear maps $\sigma:y\to x_0/x\cap x_0$ such that $\t(x)(u)\subseteq \sigma(u)$ for all $u\in x+x_0/x_0$ and ${\rm rk}(\tau-\sigma)=1$. Hence the coefficient of $(y,x\cap x_0,\tau)$ in $(D_3(x_0)\circ R_1(x_0))(x)$ is $q^2$.

By the above comments the $q$-bracket $[R_1(x_0),D_3(x_0)]_q$ sends $x$ to 
$$
q^{-1}(1-q^{2(\dim x_0/x\cap x_0)})
\sum_{(y,x\cap x_0,\tau)\in S_0}(y,x\cap x_0,\tau). 
$$
By Definitions \ref{defn:L12&R12}(iii) and \ref{defn:D12}(ii) 
the linear map $(1-q^{\dim x_0} D_2(x_0))\circ R_1(x_0)$ sends $x$ to 
\begin{gather}\label{D2R1}
(1-q^{2(\dim x_0/x\cap x_0)})\sum_{(y,x\cap x_0,\tau)\in S_0}(y,x\cap x_0,\tau).
\end{gather}
The equation (i) follows.

(ii): Let $x\in \L(\Omega)$ be given. By Definition \ref{defn:L12&R12}(iii) and Lemma \ref{lem:D4} the $q$-bracket $[R_1(x_0),D_4(x_0)]_q$ sends $x$ to the scalar multiple of (\ref{D2R1}) by $-(q-q^{-1})$. The equation (ii) follows.
\end{proof}

\begin{lem}
\label{lem:R2&D34}
The following equations hold:
\begin{enumerate}
\item $[D_3(x_0),R_2(x_0)]_q=q^{-1}(1-q^{D-\dim x_0} D_1(x_0)^{-1})\circ R_2(x_0)$.

\item $[D_4(x_0),R_2(x_0)]_q=-(q-q^{-1})(1-q^{D-\dim x_0} D_1(x_0)^{-1})\circ R_2(x_0)$.

\end{enumerate}
\end{lem}
\begin{proof}
(i): Let $x\in \L(\Omega)$ be given. Set 
\begin{align*}
S_0 &=\{(x+x_0/x_0,z,\tau)\in \L(\Omega)_{x_0}\,|\, x\cap x_0\subsetdot z, {\rm rk}\left(1_{x\cap x_0}^z\circ \t(x)- \tau\right)=0\},
\\
S_1 &=\{(x+x_0/x_0,z,\tau)\in \L(\Omega)_{x_0}\,|\, x\cap x_0\subsetdot z, {\rm rk}\left(1_{x\cap x_0}^z\circ \t(x)- \tau\right)=1\}.
\end{align*}
By Definitions \ref{defn:L12&R12}(iv) and \ref{defn:D3} both of $(D_3(x_0)\circ R_2(x_0))(x)$ and $(R_2(x_0)\circ D_3(x_0))(x)$ are the linear combinations of $S_0\cup S_1$. 

Fix an element $(x+x_0/x_0,z,\tau)\in S_0$. None of linear maps $\sigma:x+x_0/x_0\to x_0/z$ satisfies $\t(x)(u)\subseteq \sigma(u)$ for all $u\in x+x_0/x_0$ and ${\rm rk}(\tau-\sigma)=1$. 
Hence the coefficient of $(x+x_0/x_0,z,\tau)$ in $(D_3(x_0)\circ R_2(x_0))(x)$ is zero. 
There are exactly
$$
|\F|^{\dim (x+x_0/x_0)}-1=q^{2\dim (x+x_0/x_0)}-1
$$
linear maps $\sigma:x+x_0/x_0\to x_0/x\cap x_0$ satisfying ${\rm rk}(\t(x)-\sigma)=1$ and $\sigma(u)\subseteq \tau(u)$ for all $u\in x+x_0/x_0$. Hence the coefficient of $(x+x_0/x_0,z,\tau)$ in $(R_2(x_0)\circ D_3(x_0))(x)$ is $q^{2\dim (x+x_0/x_0)}-1$.

Fix an element $(x+x_0/x_0,z,\tau)\in S_1$. Clearly $1_{x\cap x_0}^z\circ \t(x)$ is the unique linear map $\sigma:x+x_0/x_0\to x_0/z$ satisfying $\t(x)(u)\subseteq \sigma(u)$ for all $u\in x+x_0/x_0$ and ${\rm rk}(\sigma-\tau)=1$. 
Hence the coefficient of $(x+x_0/x_0,z,\tau)$ in $(D_3(x_0)\circ R_2(x_0))(x)$ is one. 
Since ${\rm rk}(1_{x\cap x_0}^z\circ \t(x)-\tau)=1$  
there exists a vector $v\in x_0\setminus z$ such that
$$
{\rm Im}(1_{x\cap x_0}^z\circ \t(x)-\tau)=\F v+z.
$$ 
Suppose that $\sigma$ is any linear map $x+x_0/x_0\to x_0/x\cap x_0$ satisfying ${\rm rk}(\t(x)-\sigma)=1$ and $\sigma(u)\subseteq \tau(u)$ for all $u\in x+x_0/x_0$. Since $1_{x\cap x_0}^z\circ \sigma=\tau$ it follows that $1_{x\cap x_0}^z$ maps ${\rm Im}(\t(x)-\sigma)$ onto $\F v+z$. 
Hence 
\begin{gather}
\label{e1:R2&D34}
{\rm Im}(\t(x)-\sigma)=\F w+(x\cap x_0)
\end{gather}
for some $w\in v+z$. Let $u\in x+x_0/x_0$ be given.  By (\ref{e1:R2&D34}) there is a scalar $c\in \F$ such that 
$\t(x)(u)-\sigma(u)=cw+(x\cap x_0)$. Applying $1_{x\cap x_0}^z$ to the above equation yields that $c$ is uniquely determined by the following equation: 
$$
(1_{x\cap x_0}^z\circ \t(x))(u)-\tau(u)=c 
v+z.
$$ 
This shows that $\sigma$ is uniquely determined by the choice of $w+(x\cap x_0)\subseteq v+z$.
There are exactly $|z/x\cap x_0|=q^2$ cosets of $x\cap x_0$ in $v+z$. Hence the coefficient of $(x+x_0/x_0,z,\tau)$ in $(R_2(x_0)\circ D_3(x_0))(x)$ is $q^2$.

By the above comments the $q$-bracket $[D_3(x_0),R_2(x_0)]_q$ sends $x$ to 
$$
q^{-1}(1-q^{2\dim (x+x_0/x_0)})\sum_{(x+x_0/x_0,z,\tau)\in S_0} (x+x_0/x_0,z,\tau).
$$
By Definitions \ref{defn:L12&R12}(iv) and \ref{defn:D12}(i) the linear map $(1-q^{D-\dim x_0} D_1(x_0)^{-1})\circ R_2(x_0)$ sends $x$ to 
\begin{gather}\label{D1R2}
(1-q^{2\dim (x+x_0/x_0)})\sum_{(x+x_0/x_0,z,\tau)\in S_0} (x+x_0/x_0,z,\tau).
\end{gather}
The equation (i) follows.

(ii): Let $x\in \L(\Omega)$ be given. By Definition \ref{defn:L12&R12}(iv) and Lemma \ref{lem:D4} the $q$-bracket $[D_4(x_0),R_2(x_0)]_q$ sends $x$ to the scalar multiple of (\ref{D1R2}) by $-(q-q^{-1})$. The equation (ii) follows.
\end{proof}

We finish this section with some results on the linear map $(q^2-1)D_3(x_0)+D_4(x_0)$.

\begin{lem}
\label{lem:K3_eigenvalue}
Suppose that $y\in \L(\Omega/x_0)$ and $z\in \L(x_0)$. Then the eigenvalues of $(q^2-1)D_3(x_0)+D_4(x_0)$ in $\C^{\L(\Omega)}_{y,z}$ are 
$$
q^{2(\dim y+\dim(x_0/z)-h)}
$$
for all integers $h$ with $0\leq h\leq \min\{\dim y,\dim(x_0/z)\}$ and their multiplicities 
\begin{gather*}
q^{h(\dim y-1)}
{\dim y \brack h}_q
\prod_{i=0}^{h-1}
(q^{2(\dim(x_0/z)-i)}-1)
\end{gather*}
for all integers $h$ with $0\leq h\leq \min\{\dim y,\dim(x_0/z)\}$, respectively.
\end{lem}
\begin{proof}
It is routine to verify the lemma by using Lemmas \ref{lem:D3} and \ref{lem:D4}.
\end{proof}

\begin{lem}
[Lemma 7.9, \cite{Watanabe:2017}]
\label{lem:K3}
The following equations hold:
\begin{enumerate}
\item $[(q^2-1)D_3(x_0)+D_4(x_0),L_1(x_0)]_q=0$.

\item $[L_2(x_0),(q^2-1)D_3(x_0)+D_4(x_0)]_q=0$.

\item $[R_1(x_0),(q^2-1)D_3(x_0)+D_4(x_0)]_q=0$.

\item $[(q^2-1)D_3(x_0)+D_4(x_0),R_2(x_0)]_q=0$.
\end{enumerate}
\end{lem}
\begin{proof}
Immediate from Lemmas \ref{lem:L1&D34}--\ref{lem:R2&D34}.
\end{proof}

Recall the $\U$-module $\C^{\L(\Omega)}(\lambda)$ for any nonzero $\lambda\in \C$ from Definition  \ref{defn:L(Omega)_Umodule} and the linear map $\iota(x_0)$ from Definition \ref{defn:iota}.

\begin{lem}
\label{lem:K3diagram}
Suppose that $\lambda$ and $\mu$ are two nonozero scalars in $\C$. Then the following diagram commutes:
\begin{table}[H]
\centering
\begin{tikzpicture}
\matrix(m)[matrix of math nodes,
row sep=4em, column sep=4em,
text height=1.5ex, text depth=0.25ex]
{
\C^{\L(\Omega)}
&\C^{\L(\Omega/ x_0)}(\lambda)\otimes \C^{\L(x_0)}(\mu)\\
\C^{\L(\Omega)}
&\C^{\L(\Omega/ x_0)}(\lambda)\otimes \C^{\L(x_0)}(\mu)\\
};
\path[->,font=\scriptsize,>=angle 90]
(m-1-1) edge node[left] {$(q^2-1)D_3(x_0)+D_4(x_0)$} (m-2-1)
(m-1-1) edge node[above] {$\iota(x_0)$} (m-1-2)
(m-2-1) edge node[below] {$\iota(x_0)$} (m-2-2)
(m-1-2) edge node[right] {$q^DK^{-1}\otimes K$} (m-2-2);
\end{tikzpicture}
\end{table}
\end{lem}
\begin{proof}
Let $x\in \L(\Omega)$ be given. There are exactly 
\begin{align*}
|\L_1(x_0/x\cap x_0)|(|\F|^{\dim (x+x_0/x_0)}-1)
\end{align*}
linear maps $\tau:x+x_0/x_0\to x_0/x\cap x_0$ satisfying ${\rm rk}(\tau-\t(x))=1$. 
By Lemmas \ref{lem:|Lk|} and \ref{lem:x+x0/x0} the above number is equal to $\frac{1}{q^2-1}$ times 
\begin{gather}\label{rank1}
\left(\frac{|x_0|}{|x\cap x_0|}-1\right)\left(\frac{|x|}{|x\cap x_0|}-1\right).
\end{gather}
By Definitions \ref{defn:iota} and \ref{defn:D3} the linear map $(q^2-1)\iota(x_0)\circ D_3(x_0)$ sends $x$ to the scalar  multiple of $x+x_0/x_0\otimes x\cap x_0$ by (\ref{rank1}). Combined with Definition \ref{defn:D4} the linear map $\iota(x_0)\circ ((q^2-1)D_3(x_0)+D_4(x_0))$ sends $x$ to the  scalar multiple of $x+x_0/x_0\otimes x\cap x_0$ by 
\begin{gather}\label{rank1'}
\left(\frac{|x_0|}{|x\cap x_0|}-1\right)\left(\frac{|x|}{|x\cap x_0|}-1\right)
+\frac{|x\cup x_0|}{|x\cap x_0|}
=
\frac{|x|\cdot |x_0|}{|x\cap x_0|^2}.
\end{gather}
Using Proposition \ref{prop:L(Omega)_Umodule} and Lemma \ref{lem:x+x0/x0} yields that $x+x_0/x_0\otimes x\cap x_0$ is an eigenvector of $q^D K^{-1}\otimes K$ in $\C^{\L(\Omega/x_0)}(\lambda)\otimes \C^{\L(x_0)}(\mu)$ corresponding to the eigenvalue (\ref{rank1'}). The lemma follows.
\end{proof}

\section{The $\W$-module $\C^{\L(\Omega)}(x_0)$}\label{s:WmoduleL(Omega)}

To develop a $\W$-module structure on $\C^{\L(\Omega)}$ we need two more equations  from \cite[Lemma 7.7]{Watanabe:2017}. Here we verify the two equations from the viewpoint of the poset $\L(\Omega)_{x_0}$.

\begin{lem}
\label{lem1:L1R1-R1L1}
Suppose that $x\in \L(\Omega)$. Then the following statements hold:
\begin{enumerate}
\item There are exactly 
\begin{gather}\label{e1:L1R1-R1L1}
q^{D-\dim x+\dim x_0-\dim(x\cap x_0)-1}
\left[
D-\dim (x+x_0)
\right]_q
\end{gather}
elements $x'\in \L(\Omega)$ such that $x\subsetdot x'$ and $x'\cap x_0=x\cap x_0$.

\item There are exactly 
\begin{gather}\label{e2:L1R1-R1L1}
q^{\dim x-\dim(x\cap x_0)-1}
\left[
\dim (x+x_0/x_0)
\right]_q
\end{gather}
elements $x'\in \L(\Omega)$ such that $x'\subsetdot x$ and $x'\cap x_0=x\cap x_0$.
\end{enumerate}
\end{lem}
\begin{proof}
(i): 
There are exactly 
\begin{gather}\label{e1':L1R1-R1L1}
|\L_1(\F^{\dim(\Omega/x_0)-\dim (x+x_0/x_0)})|\cdot 
|x_0/x\cap x_0|
=
|\L_1(\F^{D-\dim (x+x_0)})|\cdot 
q^{2(\dim x_0-\dim (x\cap x_0))}
\end{gather}
elements $(y,x\cap x_0,\tau)\in \L(\Omega)_{x_0}$ satisfying $x+x_0/x_0\subsetdot y$ and $\t(x)(u)\subseteq \tau(u)$ for all $u\in x+x_0/x_0$. 
Using Lemmas \ref{lem:|Lk|} and \ref{lem:x+x0/x0} yields that (\ref{e1':L1R1-R1L1}) is equal to (\ref{e1:L1R1-R1L1}). The statement (i) follows.

(ii): 
There are exactly 
\begin{gather}\label{e2':L1R1-R1L1}
|\L_{\dim (x+x_0/x_0)-1}(\F^{\dim (x+x_0/x_0)})|
\end{gather}
elements $(y,x\cap x_0,\tau)$ satisfying $y\subsetdot x+x_0/x_0$ and $\tau(u)\subseteq \t(x)(u)$ for all $u\in y$. 
Using Lemmas \ref{lem:|Lk|} and \ref{lem:x+x0/x0} yields that (\ref{e2':L1R1-R1L1}) is equal to (\ref{e2:L1R1-R1L1}). The statement (ii) follows.
\end{proof}

\begin{lem}
\label{lem2:L1R1-R1L1}
Suppose that $x,x'\in \L(\Omega)$ satisfy the following conditions: 
\begin{enumerate}
\item[{\rm (a)}] $(x+x_0/x_0)\cap (x'+x_0/x_0)\subsetdot x+x_0/x_0$.

\item[{\rm (b)}] $(x+x_0/x_0)\cap (x'+x_0/x_0)\subsetdot x'+x_0/x_0$.

\item[{\rm (c)}] $\t(x)(u)=\t(x')(u)$ for all $u\in (x+x_0/x_0)\cap (x'+x_0/x_0)$.
\end{enumerate}
Then the following statements hold:
\begin{enumerate}
\item There exists a unique $x''\in \L(\Omega)$ with $x\subsetdot x''$, $x'\subsetdot x''$ and $x''\cap x_0=x\cap x_0$.

\item There exists a unique $x''\in \L(\Omega)$ with $x''\subsetdot x$, $x''\subsetdot x'$ and $x''\cap x_0=x\cap x_0$.
\end{enumerate}
\end{lem}
\begin{proof}
(i): Suppose that $(y,x\cap x_0,\tau)\in \L(\Omega)_{x_0}$ satisfies the following conditions:
\begin{enumerate}
\item[(1)] $x+x_0/x_0\subsetdot y$. 

\item[(2)] $x'+x_0/x_0\subsetdot y$.

\item[(3)] $\t(x)(u)=\tau(u)$ for all $u\in x+x_0/x_0$.

\item[(4)] $\t(x')(u)=\tau(u)$ for all $u\in x'+x_0/x_0$.
\end{enumerate}
Since $y$ satisfies the conditions (1) and (2) it follows from (a) and (b) that $y=x+x'+x_0/x_0$. 
By the condition (c) there exists a unique map $y\to x_0/x\cap x_0$ given by
\begin{eqnarray*}
u+u' &\mapsto & \t(x)(u) +\t(x')(u')
\quad
\hbox{for all $u\in  x+x_0/x_0$ and $u'\in  x'+x_0/x_0$}.
\end{eqnarray*}
It is routine to verify that the above map is linear. 
Since $\tau$ satisfies the conditions (3) and (4) it follows that $\tau$ is the above linear map. 
The statement (i) follows.

(ii): Suppose that $(y,x\cap x_0,\tau)\in \L(\Omega)_{x_0}$ satisfies the following conditions:
\begin{enumerate}
\item[(1)] $y\subsetdot x+x_0/x_0$. 

\item[(2)] $y\subsetdot x'+x_0/x_0$.

\item[(3)] $\tau(u)=\t(x)(u)$ for all $u\in y$.

\item[(4)] $\tau(u)=\t(x')(u)$ for all $u\in y$.
\end{enumerate}
Since $y$ satisfies the conditions (1) and (2) it follows from (a) and (b) that $y=(x+x_0/x_0)\cap (x'+x_0/x_0)$. 
Since $\tau$ satisfies the conditions (3) and (4) it follows from (c) that $\tau$ is equal to $\t(x)|_{y}=\t(x')|_y$. The statement (ii) follows.
\end{proof}

\begin{lem}
\label{lem3:L1R1-R1L1}
Suppose that $x,x'\in \L(\Omega)$ satisfy the following conditions:
\begin{enumerate}
\item[{\rm (a)}] $x+x_0/x_0=x'+x_0/x_0$.

\item[{\rm (b)}] $x\cap x_0=x'\cap x_0$.

\item[{\rm (c)}] ${\rm rk}\left(\t(x)-\t(x')\right)=1$.
\end{enumerate}
Then the following statements hold:
\begin{enumerate}
\item None of the elements $x''\in \L(\Omega)$ satisfies $x\subsetdot x''$, $x'\subsetdot x''$ and $x''\cap x_0=x\cap x_0$.

\item There exists a unique $x''\in \L(\Omega)$ with $x''\subsetdot x$, $x''\subsetdot x'$ and $x''\cap x_0=x\cap x_0$.
\end{enumerate}
\end{lem}
\begin{proof}
(i): Suppose that there exists an element $(y,x\cap x_0,\tau)\in \L(\Omega)_{x_0}$ satisfying the following conditions:
\begin{enumerate}
\item[(1)] $x+x_0/x_0\subsetdot y$. 

\item[(2)] $x'+x_0/x_0\subsetdot y$. 

\item[(3)] $\t(x)(u)=\tau(u)$ for all $u\in x+x_0/x_0$.

\item[(4)] $\t(x')(u)=\tau(u)$ for all $u\in x'+x_0/x_0$.
\end{enumerate}
Applying (a) to (3) and (4) both of $\t(x)$ and $\t(x')$ are equal to $\tau|_{x+x_0/x_0}$, a contradiction to (c). The statement (i) follows.

(ii): Suppose that $(y,x\cap x_0,\tau)\in \L(\Omega)_{x_0}$ satisfies the following conditions:
\begin{enumerate}
\item[(1)] $y\subsetdot x+x_0/x_0$. 

\item[(2)] $y\subsetdot x'+x_0/x_0$.

\item[(3)] $\tau(u)=\t(x)(u)$ for all $u\in y$.

\item[(4)] $\tau(u)=\t(x')(u)$ for all $u\in y$.
\end{enumerate}
Applying the rank-nullity theorem to the condition (c) yields that the nullity of $\t(x)-\t(x')$ is equal to $\dim (x+x_0/x_0)-1=\dim (x'+x_0/x_0)-1$. 
Combined with the conditions (1)--(4) this implies that $y=\ker \left(\t(x)-\t(x')\right)$ 
and $\tau=\t(x)|_y=\t(x')|_y$. 
The statement (ii) follows.
\end{proof}

\begin{prop}
[Lemma 7.7(ii), \cite{Watanabe:2017}]
\label{prop:L1R1-R1L1}
The following equation holds:
$$
(q^2-1)[L_1(x_0), R_1(x_0)]=
q^D D_1(x_0)\circ D_2(x_0)
-(q^2-1) D_3(x_0)-D_4(x_0).
$$
\end{prop}
\begin{proof}
Let $x\in \L(\Omega)$ be given. Set 
\begin{align*}
S_0 &=\{x'\in \L(\Omega)\,|\, \hbox{$x'$ satisfies the conditions (a)--(c) of Lemma \ref{lem2:L1R1-R1L1}}\},
\\
S_1 &=\{x'\in \L(\Omega)\,|\, \hbox{$x'$ satisfies the conditions (a)--(c) of Lemma \ref{lem3:L1R1-R1L1}}\}.
\end{align*}
By Definition \ref{defn:L12&R12}(i), (iii) both of $(L_1(x_0)\circ R_1(x_0))(x)$ and  $(R_1(x_0)\circ L_1(x_0))(x)$ are the linear combinations of $\{x\}\cup S_0\cup S_1$.

By Lemma \ref{lem1:L1R1-R1L1}(i) the coefficient of $x$ in $(L_1(x_0)\circ R_1(x_0))(x)$ is equal to (\ref{e1:L1R1-R1L1}). By Lemma \ref{lem1:L1R1-R1L1}(ii) the coefficient of $x$ in $(R_1(x_0)\circ L_1(x_0))(x)$ is equal to (\ref{e2:L1R1-R1L1}). 
Fix an element $x'\in S_0$. 
By Lemma \ref{lem2:L1R1-R1L1}(i) the coefficient of $x'$ in $(L_1(x_0)\circ R_1(x_0))(x)$ is equal to one. By Lemma \ref{lem2:L1R1-R1L1}(ii) the coefficient of $x'$ in $(R_1(x_0)\circ L_1(x_0))(x)$ is equal to one. Fix an element $x'\in S_1$. 
By Lemma \ref{lem3:L1R1-R1L1}(i) the coefficient of $x'$ in $(L_1(x_0)\circ R_1(x_0))(x)$ is equal to zero. By Lemma \ref{lem3:L1R1-R1L1}(ii) the coefficient of $x'$ in $(R_1(x_0)\circ L_1(x_0))(x)$ is equal to one. 
By the above comments the bracket $[L_1(x_0), R_1(x_0)]$ sends $x$ to 
$$
c_0 x-\sum_{x'\in S_1} x'
$$
where the scalar $c_0$ is equal to the subtraction of (\ref{e2:L1R1-R1L1}) from (\ref{e1:L1R1-R1L1}). A routine calculation shows that 
\begin{align*}
(q^2-1) c_0
&=
\frac{|\Omega|}{|x|}-
\frac{|x\cup x_0|}{|x\cap x_0|}.
\end{align*}
By Definition \ref{defn:D3} the map $D_3(x_0)$ sends $x$ to 
$\sum_{x'\in S_1} x'$. 
Using Definition \ref{defn:D12} and Lemma \ref{lem:x+x0/x0} yields that $q^D D_1(x_0)\circ D_2(x_0)$ sends $x$ to $\frac{|\Omega|}{|x|} x$.
Combined with Definition \ref{defn:D4} the proposition follows.
\end{proof}

\begin{lem}
\label{lem1:L2R2-R2L2}
Suppose that $x\in \L(\Omega)$. Then the following statements hold:
\begin{enumerate}
\item There are exactly 
\begin{gather}\label{e1:L2R2-R2L2}
q^{\dim x_0-\dim (x\cap x_0)-1}
[\dim(x_0/x\cap x_0)]_q
\end{gather}
elements $x'\in \L(\Omega)$ such that $x\subsetdot x'$ and $x'+ x_0/x_0=x+x_0/x_0$.

\item There are exactly 
\begin{gather}\label{e2:L2R2-R2L2}
q^{2\dim x-\dim(x\cap x_0)-1}[\dim (x\cap x_0)]_q
\end{gather}
elements $x'\in \L(\Omega)$ such that $x'\subsetdot x$ and $x'+ x_0/x_0=x+x_0/x_0$.
\end{enumerate}
\end{lem}
\begin{proof}
(i): There are exactly
\begin{gather}
\label{e1':L2R2-R2L2}
|\L_1(\F^{\dim (x_0/x\cap x_0)})|
\end{gather}
elements $(x+x_0/x_0,z,\tau)\in \L(\Omega)_{x_0}$ satisfying $x\cap x_0\subsetdot z$ and $\t(x)(u)\subseteq \tau(u)$ for all $u\in x+x_0/x_0$. Using Lemma \ref{lem:|Lk|} yields that (\ref{e1':L2R2-R2L2}) is equal to (\ref{e1:L2R2-R2L2}). The statement (i) follows.

(ii): There are exactly
\begin{gather}
\label{e2':L2R2-R2L2}
|\L_{\dim (x\cap x_0)-1}(\F^{\dim (x\cap x_0)})|\cdot 
|\F|^{\dim (x+x_0/x_0)}
\end{gather}
elements $(x+x_0/x_0,z,\tau)\in \L(\Omega)_{x_0}$ satisfying $z\subsetdot x\cap x_0$ and $\tau(u)\subseteq \t(x)(u)$ for all $u\in x+x_0/x_0$. Using Lemmas \ref{lem:|Lk|} and \ref{lem:x+x0/x0} yields that  (\ref{e2':L2R2-R2L2}) is equal to (\ref{e2:L2R2-R2L2}). The statement (ii) follows.
\end{proof}

\begin{lem}
\label{lem2:L2R2-R2L2}
Suppose that $x,x'\in \L(\Omega)$ satisfy the following conditions: 
\begin{enumerate}
\item[{\rm (a)}] $x+x_0/x_0=x'+x_0/x_0$.
\item[{\rm (b)}] $x\cap x'\cap x_0\subsetdot x\cap x_0$.
\item[{\rm (c)}] $x\cap x'\cap x_0\subsetdot x'\cap x_0$.
\item[{\rm (d)}] $\t(x)(u)\cap \t(x')(u)\not=\emptyset$ for all $u\in x+x_0/x_0$.  
\end{enumerate}
Then the following statements hold:
\begin{enumerate}
\item There exists a unique $x''\in \L(\Omega)$ with $x\subsetdot x''$, $x'\subsetdot x''$ and $x''+x_0/x_0=x+x_0/x_0$.

\item There exists a unique $x''\in \L(\Omega)$ with
$x''\subsetdot x$, $x''\subsetdot x'$ and $x''+x_0/x_0=x+x_0/x_0$.
\end{enumerate}
\end{lem}
\begin{proof}
(i): Suppose that $(x+x_0/x_0,z,\tau)\in \L(\Omega)_{x_0}$ satisfies the following conditions:
\begin{enumerate}
\item[(1)] $x\cap x_0\subsetdot z$.

\item[(2)] $x'\cap x_0\subsetdot z$.

\item[(3)] $\t(x)(u)\subseteq \tau(u)$ for all $u\in x+x_0/x_0$.

\item[(4)] $\t(x')(u)\subseteq \tau(u)$ for all $u\in x+x_0/x_0$.
\end{enumerate}
Since $z$ satisfies the conditions (1) and (2) it follows from (b) and (c) that $z=(x\cap x_0)+(x'\cap x_0)$. 
By the condition (d) the sets $\t(x)(u)$ and $\t(x')(u)$ are in the same coset of $z$ in $x_0$ for all $u\in x+x_0/x_0$. 
Hence $1_{x\cap x_0}^z\circ \t(x)=1_{x'\cap x_0}^z\circ\t(x')$. 
Since $\tau$ satisfies the conditions (3) and (4) it follows that $\tau$ is the above map. The statement (i) follows.

(ii): Suppose that $(x+x_0/x_0,z,\tau)\in \L(\Omega)_{x_0}$ satisfies the following conditions:
\begin{enumerate}
\item[(1)] $z\subsetdot x\cap x_0$.

\item[(2)] $z\subsetdot x'\cap x_0$.

\item[(3)] $\tau(u)\subseteq \t(x)(u)$ for all $u\in x+x_0/x_0$.

\item[(4)] $\tau(u)\subseteq \t(x')(u)$ for all $u\in x+x_0/x_0$.
\end{enumerate}
Since $z$ satisfies the conditions (1) and (2) it follows from (b) and (c) that $z=x\cap x'\cap x_0$. 
By the condition (d) there exists a unique map $x+x_0/x_0\to x_0/z$ given by 
\begin{eqnarray*}
u &\mapsto & \t(x)(u)\cap \t(x')(u)
\qquad \hbox{for all $u\in x+x_0/x_0$}.
\end{eqnarray*}
It is routine to verify that the above map is linear. Since $\tau$ satisfies the conditions (3) and (4) it follows that $\tau$ is the aforementioned linear map. The statement (ii) follows. 
\end{proof}

\begin{lem}
\label{lem3:L2R2-R2L2}
Suppose that $x,x'\in \L(\Omega)$ satisfy the following conditions: 
\begin{enumerate}
\item[{\rm (a)}] $x+x_0/x_0=x'+x_0/x_0$.

\item[{\rm (b)}] $x\cap x_0=x'\cap x_0$.

\item[{\rm (c)}] ${\rm rk}(\t(x)-\t(x'))=1$.
\end{enumerate}
Then the following statements hold:
\begin{enumerate}
\item There exists a unique $x''\in \L(\Omega)$ with $x\subsetdot x''$, $x'\subsetdot x''$ and $x''+x_0/x_0=x+x_0/x_0$.

\item None of the elements $x''\in \L(\Omega)$ satisfies
$x''\subsetdot x$, $x''\subsetdot x'$ and $x''+x_0/x_0=x+x_0/x_0$.
\end{enumerate}
\end{lem}
\begin{proof}
(i): Suppose that $(x+x_0/x_0,z,\tau)\in \L(\Omega)_{x_0}$ satisfies the following conditions:
\begin{enumerate}
\item[(1)] $x\cap x_0\subsetdot z$.

\item[(2)] $x'\cap x_0\subsetdot z$.

\item[(3)] $\t(x)(u)\subseteq \tau(u)$ for all $u\in x+x_0/x_0$.

\item[(4)] $\t(x')(u)\subseteq \tau(u)$ for all $u\in x+x_0/x_0$.
\end{enumerate}
Using Lemma \ref{lem:x+x0/x0} yields that 
\begin{align}
\dim x &= \dim(x+x_0/x_0)+\dim(x\cap x_0),
\label{dimx}
\\
\dim(x+x') &= \dim(x+x'+x_0/x_0)+\dim((x+x')\cap x_0).
\label{dim(x+x')}
\end{align}
By the condition (a) the quotient space $x+x'+x_0/x_0=x+x_0/x_0$. 
Applying Lemma \ref{lem2:bilinear} to (c) yields that 
the left-hand side of (\ref{dim(x+x')}) is equal to 
$1+\dim x$. 
Subtracting (\ref{dim(x+x')}) from (\ref{dimx})
yields $\dim((x+x')\cap x_0)=1+\dim(x\cap x_0)$. Hence 
$x\cap x_0\subsetdot (x+x')\cap x_0$ and by a similar argument  $x'\cap x_0\subsetdot (x+x')\cap x_0$. Since $z$ satisfies the conditions (1) and (2) it follows that $z=(x+x')\cap x_0$. 
Since $\tau$ satisfies the conditions (3) and (4) it follows from Definition \ref{defn:taux} that $\tau=\t(x+x')$. The statement (i) follows.

(ii): Suppose that there exists an element $(x+x_0/x_0,z,\tau)\in \L(\Omega)_{x_0}$ satisfying the following conditions:
\begin{enumerate}
\item[(1)] $z\subsetdot x\cap x_0$.

\item[(2)] $z\subsetdot x'\cap x_0$.

\item[(3)] $\tau(u)\subseteq \t(x)(u)$ for all $u\in x+x_0/x_0$.

\item[(4)] $\tau(u)\subseteq \t(x')(u)$ for all $u\in x+x_0/x_0$.
\end{enumerate}
Applying (b) to (3) and (4) both of $\t(x)$ and $\t(x')$ are equal to $1_z^{x\cap x_0}\circ \tau$, a contradiction to (c). 
The statement (ii) follows.
\end{proof}

\begin{prop}
[Lemma 7.7(i), \cite{Watanabe:2017}]
\label{prop:L2R2-R2L2}
The following equation holds:
$$
(q^2-1)[L_2(x_0), R_2(x_0)]=(q^2-1)D_3(x_0)+D_4(x_0)-q^{D} D_1(x_0)^{-1}\circ D_2(x_0)^{-1}.
$$
\end{prop}
\begin{proof}
Let $x\in \L(\Omega)$ be given. Set 
\begin{align*}
S_0 &=\{x'\in \L(\Omega)\,|\, \hbox{$x'$ satisfies the conditions (a)--(d) of Lemma \ref{lem2:L2R2-R2L2}}\},
\\
S_1 &=\{x'\in \L(\Omega)\,|\, \hbox{$x'$ satisfies the conditions (a)--(c) of Lemma \ref{lem3:L2R2-R2L2}}\}.
\end{align*}
By Definition \ref{defn:L12&R12}(ii), (iv) both of $(L_2(x_0)\circ R_2(x_0))(x)$ and $(R_2(x_0)\circ L_2(x_0))(x)$ are the linear combinations of $\{x\}\cup S_0\cup S_1$. By Lemma \ref{lem1:L2R2-R2L2}(i) the coefficient of $x$ in $(L_2(x_0)\circ R_2(x_0))(x)$ is equal to (\ref{e1:L2R2-R2L2}). By Lemma \ref{lem1:L2R2-R2L2}(ii) the coefficient of $x$ in $(R_2(x_0)\circ L_2(x_0))(x)$ is equal to (\ref{e2:L2R2-R2L2}). Fix an element $x'\in S_0$. By Lemma \ref{lem2:L2R2-R2L2}(i) the coefficient of $x'$ in $(L_2(x_0)\circ R_2(x_0))(x)$ is equal to one. By Lemma \ref{lem2:L2R2-R2L2}(ii) the coefficient of $x'$ in $(R_2(x_0)\circ L_2(x_0))(x)$ is equal to one. Fix an element $x'\in S_1$. By Lemma \ref{lem3:L2R2-R2L2}(i) the coefficient of $x'$ in $(L_2(x_0)\circ R_2(x_0))(x)$ is equal to one. By Lemma \ref{lem3:L2R2-R2L2}(ii) the coefficient of $x'$ in $(R_2(x_0)\circ L_2(x_0))(x)$ is equal to zero. By the above comments the bracket $[L_2(x_0),R_2(x_0)]$ sends $x$ to 
$$
c_0 x+\sum_{x'\in S_1} x'
$$
where $c_0$ is equal to the subtraction of (\ref{e2:L2R2-R2L2}) from (\ref{e1:L2R2-R2L2}). A direct calculation shows that 
$$
(q^2-1)c_0=\frac{|x\cup x_0|}{|x\cap x_0|}-|x|.
$$
By Definition \ref{defn:D3} the map $D_3(x_0)$ sends $x$ to $\sum_{x'\in S_1} x'$.
Using Definition \ref{defn:D12} yields that $q^{D} D_1(x_0)^{-1}\circ D_2(x_0)^{-1}$ sends $x$ to the scalar multiple of $x$ by $|x|$. Combined with Definition \ref{defn:D4} the proposition follows.
\end{proof}

Lemmas \ref{lem:E1diagram}--\ref{lem:K2diagram} and \ref{lem:K3diagram} are the inspiration to create the following $\W$-module.
Note that the linear map $(q^2-1)D_3(x_0)+D_4(x_0)$ is invertible by Definition \ref{defn:Cyz} and Lemma \ref{lem:K3_eigenvalue}.

\begin{thm}
\label{thm:Wmodule_L(Omega)}
There exists a unique $\W$-module $\C^{\L(\Omega)}$ given by 
\begin{align*}
E_1 &= q^{\dim x_0-D} L_1(x_0),
\\ 
E_2 &= q^{\dim x_0-D} D_1(x_0)\circ L_2(x_0),
\\
F_1 &= q^{1-\dim x_0} R_1(x_0)\circ D_2(x_0)^{-1},
\\ 
F_2 &= q^{1-\dim x_0} R_2(x_0),
\\
K_1^{\pm 1} &= D_1(x_0)^{\pm 1},
\\
K_2^{\pm 1} &= D_2(x_0)^{\pm 1},
\\
I^{\pm 1} &= q^{\mp D} D_1(x_0)^{\pm 1}\circ D_2(x_0)^{\mp 1}\circ ((q^2-1)D_3(x_0)+D_4(x_0))^{\pm 1}.
\end{align*}
\end{thm}
\begin{proof}
For the moment set $E_1,E_2,F_1,F_2,K_1^{\pm 1},K_2^{\pm 1},I^{\pm 1}$ to be the corresponding linear maps on $\C^{\L(\Omega)}$ stated in Theorem \ref{thm:Wmodule_L(Omega)}.
By construction the relations (\ref{W:Ii})--(\ref{W:K2i}) hold. 
By Lemma \ref{lem:D1234} the element $I$ commutes with $K_1^{\pm 1}$ and $K_2^{\pm 1}$. By Lemmas \ref{lem:L12R12D12} and \ref{lem:K3} the element $I$ commutes with $E_1,E_2,F_1,F_2$. Therefore the relation (\ref{W:Icenter}) holds.
Applying Lemmas \ref{lem:L12R12D12}, \ref{lem:L1L2R1R2} and \ref{lem:D1234} it is routine to verify (\ref{W:K1E2})--(\ref{W:E1E2F1F2}).   
By Lemmas \ref{lem:L12R12D12}(i), \ref{lem:D1234} and Proposition \ref{prop:L1R1-R1L1} the relation (\ref{W:E1F1}) holds. 
By Lemmas \ref{lem:L12R12D12}(iv), \ref{lem:D1234} and Proposition \ref{prop:L2R2-R2L2} the relation (\ref{W:E2F2}) holds. 
The existence of the $\W$-module $\C^{\L(\Omega)}$ follows. 
Since the algebra $\W$ is generated by $E_1,E_2,F_1,F_2,K_1^{\pm 1},K_2^{\pm 1},I^{\pm 1}$ the uniqueness follows.
\end{proof}

\begin{defn}
\label{defn:Wmodule_L(Omega)(x0)}
We denote by $\C^{\L(\Omega)}(x_0)$ the $\W$-module shown in Theorem \ref{thm:Wmodule_L(Omega)}. By Theorem \ref{thm1:H->W} the $\W$-module $\C^{\L(\Omega)}(x_0)$ is also an $\H$-module.
\end{defn}

Recall the algebra homomorphism $\widetilde{\Delta}:\U\to \W$ from Theorem \ref{thm1:U->W}. In light of the following lemma the $\W$-module $\C^{\L(\Omega)}(x_0)$ can be considered as an extension of the $\U$-module $\C^{\L(\Omega)}(q^{\dim x_0})$.

\begin{thm}
\label{thm:X&tildeDeltaXq}
The following diagram commutes for any $X\in \U$:
\begin{table}[H]
\centering
\begin{tikzpicture}
\matrix(m)[matrix of math nodes,
row sep=4em, column sep=4em,
text height=1.5ex, text depth=0.25ex]
{
\C^{\L(\Omega)}(q^{\dim x_0})
&\C^{\L(\Omega)}(x_0)
\\
\C^{\L(\Omega)}(q^{\dim x_0})
&\C^{\L(\Omega)}(x_0)
\\
};
\path[->,font=\scriptsize,>=angle 90]
(m-1-1) edge node[left] {$X$} (m-2-1)
(m-1-1) edge node[above] {$1$} (m-1-2)
(m-2-1) edge node[below] {$1$} (m-2-2)
(m-1-2) edge node[right] {$\widetilde{\Delta}(X)$} (m-2-2);
\end{tikzpicture}
\end{table}
\noindent Here $1$ denotes the identity map on $\C^{\L(\Omega)}$.
\end{thm}
\begin{proof}
Recall the images of $E,F,K^{\pm 1}$ under $\widetilde{\Delta}$ from (\ref{E'})--(\ref{K'}). By Theorem \ref{thm:Wmodule_L(Omega)} 
the following equations hold on the $\W$-module $\C^{\L(\Omega)}(x_0)$: 
\begin{align*}
\widetilde{\Delta}(E)&=q^{\dim x_0-D} (L_1(x_0)+L_2(x_0)), 
\\
\widetilde{\Delta}(F)&=q^{1-\dim x_0}(R_1(x_0)+R_2(x_0)), 
\\
\widetilde{\Delta}(K^{\pm 1})&=D_1(x_0)^{\pm 1}\circ D_2(x_0)^{\pm 1}
=D_2(x_0)^{\pm 1}\circ D_1(x_0)^{\pm 1}.
\end{align*}
Combined with Lemma \ref{lem:L1+L2} the above diagram commutes for $X=E,F,K^{\pm 1}$. Since the algebra $\U$ is generated by $E,F,K^{\pm 1}$ the result follows.
\end{proof}

Recall the algebra homomorphism $\flat:\W\to \U\otimes \U$ from Theorem \ref{thm:W->U2}. In light of the following lemma the $\W$-module $\C^{\L(\Omega)}(x_0)$ can be considered as a lift of the $\U\otimes \U$-module $\C^{\L(\Omega/x_0)}(1)\otimes \C^{\L(x_0)}(q^{\dim x_0})$ across $\iota(x_0)$.

\begin{thm}
\label{thm:flatdiagram}
The following diagram commutes for any $X\in \W$:
\begin{table}[H]
\centering
\begin{tikzpicture}
\matrix(m)[matrix of math nodes,
row sep=4em, column sep=4em,
text height=1.5ex, text depth=0.25ex]
{
\C^{\L(\Omega)}(x_0)
&\C^{\L(\Omega/ x_0)}(1)\otimes \C^{\L(x_0)}(q^{\dim x_0})\\
\C^{\L(\Omega)}(x_0)
&\C^{\L(\Omega/ x_0)}(1)\otimes \C^{\L(x_0)}(q^{\dim x_0})\\
};
\path[->,font=\scriptsize,>=angle 90]
(m-1-1) edge node[left] {$X$} (m-2-1)
(m-1-1) edge node[above] {$\iota(x_0)$} (m-1-2)
(m-2-1) edge node[below] {$\iota(x_0)$} (m-2-2)
(m-1-2) edge node[right] {$\flat(X)$} (m-2-2);
\end{tikzpicture}
\end{table}
\end{thm}
\begin{proof}
Recall the images of $E_1,E_2,F_1,F_2,K_1^{\pm 1},K_2^{\pm 1},I^{\pm 1}$ under $\flat$ from Theorem \ref{thm:W->U2}. 
By Lemmas \ref{lem:E1diagram}--\ref{lem:K2diagram} and \ref{lem:K3diagram} the diagram commutes for $X=E_1,E_2,F_1,F_2,K_1^{\pm 1},K_2^{\pm 1},I^{\pm 1}$.  Since the algebra $\W$ is generated by $E_1,E_2,F_1,F_2,K_1^{\pm 1},K_2^{\pm 1},I^{\pm 1}$ the result follows.
\end{proof}

We remark that Theorem \ref{thm:X&DeltaXq} is an amalgamation of Theorems \ref{thm:X&tildeDeltaXq} and \ref{thm:flatdiagram}.

\section{The decomposition of the $\W$-module $\C^{\L(\Omega)}(x_0)$}\label{s:dec_WmoduleL(Omega)}

Let $V$ denote a vector space. For any integer $n \geq 1$ the notation $n\cdot V$ stands for 
$$
\underbrace{V\oplus V\oplus \cdots \oplus V}_{\hbox{{\tiny $n$ copies of $V$}}}.
$$ 
Recall the preliminaries on finite-dimensional irreducible $\W$-modules from Section \ref{s:Wmodule}.

\begin{lem}
\label{lem:L(Omega)_dec}
The $\W$-module $\C^{\L(\Omega)}(x_0)$ is equal to 
\begin{gather}\label{K}
\bigoplus_{i=0}^{D-\dim x_0}
\bigoplus_{j=0}^{\dim x_0}
\bigoplus_{h=0}^{\min\{i,\dim x_0-j\}}
V_{hij}
\end{gather}
where $V_{hij}$ is the subspace of $\C^{\L(\Omega)}(x_0)$ spanned by the weight vectors of $\C^{\L(\Omega)}(x_0)$ that have the weight $(q^{D-\dim x_0-2i}, q^{\dim x_0-2j},q^{-2h})$.  Moreover the dimension of $V_{hij}$ is equal to 
\begin{align*}
q^{h(i-1)+i(D-\dim x_0-i)+j(\dim x_0-j)}
{i\brack h}_q 
{D-\dim x_0\brack i}_q
{\dim x_0\brack j}_q
\prod_{\ell=0}^{h-1}(q^{2(\dim x_0-j-\ell)}-1).
\end{align*}
\end{lem}
\begin{proof}
Immediate from Definitions \ref{defn:D12}, \ref{defn:Cyz} and Lemmas \ref{lem:|Lk|}, \ref{lem:K3_eigenvalue}.
\end{proof}

\begin{defn}
\label{defn:mhij}
Define 
\begin{align*}
m_{hij}(k)&=
q^{(h+i)(D-k-i+1)+j(k-j+1)-2h}
\frac{[D-k-h-2i+1]_q[k-h-2j+1]_q}{[D-k-h-i+1]_q[k-h-j+1]_q}
\\
&\qquad \times \;
{h+i\brack h}_q
{D-k\brack h+i}_q
{k\brack j}_q
\prod_{\ell=0}^{h-1}(q^{2(k-j-\ell)}-1)
\end{align*}
for any $h,i,j,k\in \N$ with $i\not=D-k-h+1$ and $j\not=k-h+1$.
\end{defn}

\begin{thm}
\label{thm:L(Omega)_Wdec}
The $\W$-module $\C^{\L(\Omega)}(x_0)$ is isomorphic to 
\begin{gather*}
\bigoplus_{(h,i,j)\in \P}
m_{hij}(\dim x_0)
\cdot 
(L_{D-\dim x_0-h-2i}\otimes L_{\dim x_0-h-2j})^{1,q^{-h}}
\end{gather*}
where $\P$ is the set consisting of all triples $(h,i,j)\in \Z^3$ satisfying the following conditions:
\begin{align*}
&0\leq h\leq \min\{D-\dim x_0,
\dim x_0\},
\\
&0\leq i\leq \textstyle\frac{D-\dim x_0-h}{2},
\\
&0\leq j\leq \textstyle\frac{\dim x_0-h}{2}.
\end{align*}
\end{thm}
\begin{proof}
Let $\mathbf K$ denote the set of all indices $(h,i,j)$ in the direct sum (\ref{K}). It is routine to verify that $(h,h+i,j)\in \mathbf K$ for all $(h,i,j)\in \P$. 
Recall from Lemma \ref{lem:L(Omega)_dec} the subspaces $V_{hij}$ of $\C^{\L(\Omega)}(x_0)$ for all $(h,i,j)\in \mathbf K$.
Set 
$$
d_{hij}=\dim V_{h,h+i,j} 
\qquad 
\hbox{for all $(h,i,j)\in \P$}.
$$

Let $(h,i,j)\in \Z^3$ with $(h,h+i,j)\in \mathbf K$ be given. 
By the setting of $\mathbf K$ the triple $(h,i,j)$ satisfies
\begin{align}\label{scope_hij}
0\leq h\leq \min\{D-\dim x_0,\dim x_0\},
\quad 
i\geq 0,
\quad 
j\geq 0.
\end{align}
Now suppose that $v\in V_{h,h+i,j}$ is a highest weight vector of $\C^{\L(\Omega)}(x_0)$. 
By Theorem \ref{thm:irrWmodule} there is a unique quadruple $(m,n,\delta,\lambda)\in \I(\W)$ such that the $\W$-module $(L_m\otimes L_n)^{\delta,\lambda}$ is isomorphic to 
the $\W$-submodule of $\C^{\L(\Omega)}(x_0)$ generated by $v$. By Corollary \ref{cor:irrWmodule} it is routine to verify that 
$$
(m,n,\delta,\lambda)=(D-\dim x_0-h-2i,\dim x_0-h-2j,1,q^{-h}).
$$
Using $m\geq 0$ and $n\geq 0$ yields that 
$i\leq \textstyle \frac{D-\dim x_0-h}{2}$
and 
$j\leq \textstyle \frac{\dim x_0-h}{2}$. 
Combined with (\ref{scope_hij}) the triple $(h,i,j)\in \P$. 
By Theorem \ref{thm:W_reducible} the $\W$-module $\C^{\L(\Omega)}(x_0)$ is completely reducible. 
In view of the above comments,  
there are $m_{hij}\in \N$ for all $(h,i,j)\in \P$ such that 
the $\W$-module $\C^{\L(\Omega)}$ is isomorphic to 
\begin{gather}
\label{e1:L(Omega)_Wdec}
\bigoplus_{(h,i,j)\in \P} m_{hij}\cdot (L_{D-\dim x_0-h-2i}\otimes L_{\dim x_0-h-2j})^{1,q^{-h}}.
\end{gather}
If $m_{hij}=0$ then $m_{hij}\cdot (L_{D-\dim x_0-h-2i}\otimes L_{\dim x_0-h-2j})^{1,q^{-h}}$ is interpreted as the zero vector space.

Let $(h,i,j),(h',i',j')\in \P$ be given. 
By Proposition \ref{prop:Lmn} and Lemma \ref{lem:delta&lambda} the $\W$-module $(L_{D-\dim x_0-2i'-h'}\otimes L_{\dim x_0-2j'-h'})^{1,q^{-h'}}$ contains a weight vector of $\C^{\L(\Omega)}(x_0)$ with weight 
$$
(q^{D-\dim x_0-2(h+i)}, q^{\dim x_0-2j},q^{-2h})
$$ 
if and only if $h'=h$, $i'\leq i$ and $j'\leq j$.
Using the above comment to count $d_{hij}$ in the decomposition (\ref{e1:L(Omega)_Wdec}) of $\C^{\L(\Omega)}(x_0)$ yields that 
\begin{align}
\label{d_hij&m_hij}
d_{hij}=
\sum_{\substack{(h,i',j')\in \P
\\
i'\leq i,\,
j'\leq j}} 
m_{hi'j'}
\qquad 
\hbox{for all $(h,i,j)\in \P$}.
\end{align}
We rewrite the right-hand side of (\ref{d_hij&m_hij}) as follows:
\begin{align*}
d_{hij}=
m_{hij}
+
\sum_{\substack{(h,i',j')\in \P
\\
i'< i,\,
j'\leq j}} 
m_{hi'j'}
+
\sum_{\substack{(h,i',j')\in \P
\\
i'\leq i,\,
j'< j}} 
m_{hi'j'}
-
\sum_{\substack{(h,i',j')\in \P
\\
i'< i,\,
j'< j}} 
m_{hi'j'}.
\end{align*}
Using (\ref{d_hij&m_hij}) to simplify the above equation yields that 
\begin{align*}
m_{hij}=\left\{
\begin{array}{ll}
d_{hij}-d_{h,i-1,j}-d_{h,i,j-1}+d_{h,i-1,j-1}
\qquad
&\hbox{if $i,j>0$},
\\
d_{hij}-d_{h,i,j-1}
\qquad
&\hbox{if $i=0$ and $j>0$},
\\
d_{hij}-d_{h,i-1,j}
\qquad
&\hbox{if $i>0$ and $j=0$},
\\
d_{hij}
\qquad
&\hbox{if $i=j=0$}
\end{array}
\right.
\end{align*}
for all $(h,i,j)\in \P$. Substituting the formula for $d_{hij}$ given in Lemma \ref{lem:L(Omega)_dec} into the above equation, it is routine to derive that 
$$
m_{hij}=m_{hij}(\dim x_0)>0
\qquad
\hbox{for all $(h,i,j)\in \P$}.
$$
The result follows.
\end{proof}

We remark that Theorem \ref{thm:L(Omega)_Wdec} matches the results of  \cite[Sections 4 and 5]{Dunkl77}.

\section{The decomposition of the $\H$-module $\C^{\L_k(\Omega)}(x_0)$ for $x_0\in \L_k(\Omega)$}
\label{s:dec_HmoduleL(Omega)}

Recall from Definition \ref{defn:V(theta)} the notation $V(\theta)$ for any $\W$-module $V$ and $\theta\in \C$.

\begin{lem}
\label{lem:Lk(Omega)=L(Omega)(D-2k)}
$
\C^{\L_k(\Omega)}=\C^{\L(\Omega)}(x_0)(q^{D-2k})
$
for any integer $k$ with $0\leq k\leq D$.
\end{lem}
\begin{proof}
By Theorem \ref{thm:X&tildeDeltaXq} the action of $\widetilde{\Delta}(K)$ on the $\W$-module $\C^{\L(\Omega)}(x_0)$ is identical to the action of $K$ on the $\U$-module $\C^{\L(\Omega)}(q^{\dim x_0})$. 
Combined with Proposition \ref{prop:L(Omega)_Umodule} the lemma follows.
\end{proof}

\begin{defn}
\label{defn:H_Lk(Omega)(x0)}
Recall the $\H$-module $\C^{\L(\Omega)}(x_0)$ from Definition \ref{defn:Wmodule_L(Omega)(x0)}. 
Suppose that $k$ is an integer with $0\leq k\leq D$. 
It follows from Lemmas \ref{lem:widetildeK&Hmodule} and \ref{lem:Lk(Omega)=L(Omega)(D-2k)} that $\C^{\L_k(\Omega)}$ 
is an $\H$-submodule of $\C^{\L(\Omega)}(x_0)$. We denote this $\H$-module by $\C^{\L_k(\Omega)}(x_0)$.
\end{defn}

\begin{defn}
\label{defn:Pk}
For any integer $k$ with $0\leq k\leq D$ let $\mathbf P(k)$ denote the set of all triples $(h,i,j)\in \Z^3$ satisfying the following conditions:
\begin{align*}
&0\leq h\leq \min\{D-k,k\},
\\
&0\leq i\leq 
\textstyle \frac{D-k-h}{2},
\\
&0\leq j\leq 
\textstyle
\min\left\{
D-k-h-i,
k-h-i,
\frac{k-h}{2}
\right\}.
\end{align*}
\end{defn}

\begin{prop}
\label{prop:L(Omega)_Hdec}
Suppose that $x_0\in \L_k(\Omega)$ where $k$ is an integer with $0\leq k\leq D$. 
Then the following statements hold:
\begin{enumerate}
\item The $\H$-module $\C^{\L_k(\Omega)}(x_0)$ is isomorphic to 
$$
\bigoplus_{(h,i,j)\in \mathbf P(k)}
m_{hij}(k)\cdot 
(L_{D-k-h-2i}\otimes L_{k-h-2j})^{1,q^{-h}}
(q^{D-2k}).
$$

\item For any $(h,i,j)\in \mathbf P(k)$ the $\H$-module $(L_{D-k-h-2i}\otimes L_{k-h-2j})^{1,q^{-h}}(q^{D-2k})$ is isomorphic to the irreducible $\H$-module $V_d(a,b,c)$ where
\begin{align*}
a&=q^{k-h-\min\{D-k-i,k-j\}-\max\{i,j\}},
\\
b&=q^{D-h-i-j-\min\{D-k-i,k-j\}-\min\{i,j\}+1},
\\
c&=q^{k-h-D+\min\{D-k-i,k-j\}+\max\{i,j\}},
\\
d&=\min\{D-k-i,k-j\}-\max\{i,j\}-h.
\end{align*}
\end{enumerate}
\end{prop}
\begin{proof}
(i): By Theorem \ref{thm:L(Omega)_Wdec} and Lemma \ref{lem:Lk(Omega)=L(Omega)(D-2k)} the $\H$-module $\C^{\L_k(\Omega)}(x_0)$ is isomorphic to 
$$
\bigoplus_{h=0}^{\min\{D-k,k\}}
\bigoplus_{i=0}^{\floor{\frac{D-k-h}{2}}}
\bigoplus_{j=0}^{\floor{\frac{k-h}{2}}}
m_{hij}(k)\cdot 
(L_{D-k-h-2i}\otimes L_{k-h-2j})^{1,q^{-h}}
(q^{D-2k}).
$$
Applying Theorem \ref{thm:Lmn_Hmodule}(i) with $(m,n,\delta,\lambda)=(D-k-h-2i,k-h-2j,1,q^{-h})$ yields that $(L_{D-k-h-2i}\otimes L_{k-h-2j})^{1,q^{-h}}
(q^{D-2k})=\{0\}$ for any integer $j>\min\{D-k-h-i,k-h-i\}$. The statement (i) follows.

(ii): It is routine to verify (ii) by applying Theorem \ref{thm:Lmn_Hmodule}(ii).
\end{proof}

\begin{thm}
\label{thm:L(Omega)_Hdec}
Suppose that $x_0\in \L_k(\Omega)$ where $k$ is an integer with $0\leq k\leq D$. Then the following statements hold:
\begin{enumerate}
\item Suppose that $k\not=\frac{D}{2}$. Then the $\H$-module $\C^{\L_k(\Omega)}(x_0)$ is isomorphic to 
$$
\bigoplus_{(h,i,j)\in \mathbf P(k)}
m_{hij}(k)\cdot 
(L_{D-k-h-2i}\otimes L_{k-h-2j})^{1,q^{-h}}
(q^{D-2k}).
$$
Moreover the irreducible $\H$-modules $(L_{D-k-h-2i}\otimes L_{k-h-2j})^{1,q^{-h}}
(q^{D-2k})$ for all $(h,i,j)\in \mathbf P(k)$ are mutually non-isomorphic.

\item Suppose that $k=\frac{D}{2}$. Then the $\H$-module $\C^{\L_k(\Omega)}(x_0)$ is isomorphic to 
$$
\bigoplus_{h=0}^{\frac{D}{2}}
\bigoplus_{i=0}^{\floor{\frac{D-2h}{4}}}
\bigoplus_{j=i}^{\floor{\frac{D-2h}{4}}}
m_{hij}\cdot
(L_{\frac{D}{2}-h-2i}\otimes L_{\frac{D}{2}-h-2j})^{1,q^{-h}}
(1)
$$
where 
$$
m_{hij}
=\left\{
\begin{array}{ll}
2m_{hij}(\frac{D}{2}) \qquad 
&\hbox{if $i<j$},
\\
m_{hii}(\frac{D}{2}) \qquad 
&\hbox{if $i=j$}.
\end{array}
\right.
$$
Moreover the irreducible $\H$-modules $(L_{\frac{D}{2}-h-2i}\otimes L_{\frac{D}{2}-h-2j})^{1,q^{-h}}(1)$ for all integers $h,i,j$ with $0\leq h\leq \frac{D}{2}$ and $0\leq i\leq j\leq \frac{D-2h}{4}$ are mutually non-isomorphic.
\end{enumerate}
\end{thm}
\begin{proof}
(i): Let $(h,i,j), (h',i',j')\in \mathbf P(k)$ be given. Applying Theorem \ref{thm:Lmn_Hmodule_iso} the $\H$-module $(L_{D-k-h-2i}\otimes L_{k-h-2j})^{1,q^{-h}} 
(q^{D-2k})$ is isomorphic to $(L_{D-k-h'-2i'}\otimes L_{k-h'-2j'})^{1,q^{-h'}}
(q^{D-2k})$ if and only if $(h,i,j)=(h',i',j')$. Combined with Proposition \ref{prop:L(Omega)_Hdec} the statement (i) follows.

(ii): In this case 
$$
\textstyle\mathbf P(\frac{D}{2})=\{(h,i,j)\in \Z^3\,|\, 
0\leq h\leq \frac{D}{2},0\leq i,j\leq \frac{D-2h}{4}\}.
$$ 
Let $(h,i,j), (h',i',j')\in \mathbf P(\frac{D}{2})$ be given. Applying Theorem \ref{thm:Lmn_Hmodule_iso} the $\H$-module $(L_{D-k-h-2i}\otimes L_{k-h-2j})^{1,q^{-h}}
(1)$ is isomorphic to $(L_{D-k-h'-2i'}\otimes L_{k-h'-2j'})^{1,q^{-h'}}(1)$ if and only if $(h',i',j')\in\{(h,i,j),(h,j,i)\}$. 
Using Definition \ref{defn:mhij} it is routine to verify that 
$$
m_{hij}({\textstyle\frac{D}{2}})=m_{hji}({\textstyle\frac{D}{2}}) 
\qquad \hbox{for all $(h,i,j)\in \mathbf P({\textstyle\frac{D}{2}})$}.
$$ 
Combined with Proposition \ref{prop:L(Omega)_Hdec} the statement (ii) follows.
\end{proof}

\begin{defn}
\label{defn:tildeT}
Assume that $x_0\in \L_k(\Omega)$ where $k$ is an integer with $0\leq k\leq D$. 
Let 
\begin{gather*}
\widetilde{\T}(x_0)={\rm Im}
\left(\H\to {\rm End}(\C^{\L_k(\Omega)})\right).
\end{gather*}
Here $\H\to {\rm End}(\C^{\L_k(\Omega)})$ is the representation of $\H$ into ${\rm End}(\C^{\L_k(\Omega)})$ corresponding to the $\H$-module $\C^{\L_k(\Omega)}(x_0)$. 
\end{defn}

\begin{thm}
\label{thm:tildeT}
Suppose that $x_0\in \L_k(\Omega)$ where $k$ is an integer with $0\leq k\leq D$. Then the following statements hold:
\begin{enumerate}
\item Suppose that $k\not=\frac{D}{2}$. Then the algebra $\widetilde{\T}(x_0)$ is isomorphic to 
\begin{gather*}
\bigoplus_{(h,i,j)\in \P(k)}
{\rm End}(\C^{\min\{D-k-i,k-j\}-\max\{i,j\}-h+1}).
\end{gather*}

\item Suppose that $k=\frac{D}{2}$. Then the algebra $\widetilde{\T}(x_0)$ is isomorphic to 
\begin{gather*}
\bigoplus_{h=0}^{\frac{D}{2}}
\bigoplus_{i=0}^{\floor{\frac{D-2h}{4}}}
(i+1)\cdot {\rm End}(\C^{\frac{D}{2}-h-2i+1}).
\end{gather*}
\end{enumerate}
\end{thm}
\begin{proof}
By Proposition \ref{prop:L(Omega)_Hdec} the finite-dimensional algebra  $\widetilde{\T}(x_0)$ is semisimple. Hence (i) and (ii) are immediate from Theorem \ref{thm:L(Omega)_Hdec}(i), (ii) respectively.
\end{proof}

Assume that $k$ is an integer with $0\leq k\leq D$. Observe that there exists a bijection $\P(k)\to \P(D-k)$ given by 
\begin{eqnarray*}
(h,i,j) &\mapsto & (h,j,i)
\qquad 
\hbox{for all $(h,i,j)\in \P(k)$}.
\end{eqnarray*}
Combined with Theorem \ref{thm:tildeT} it suffices to study the algebra $\widetilde{\T}(x_0)$ under the assumption $x_0\in \L_k(\Omega)$ with $0\leq k\leq \frac{D}{2}$. 
Recall that the binomial coefficients are defined as  
$$
{n\choose \ell}=\prod_{i=1}^\ell \frac{n-i+1}{\ell-i+1}
\qquad 
\hbox{for all $n\in \Z$ and $\ell\in \N$}.
$$
Recall that for any real number $a$ the notation $\floor{a}$ denotes the greatest integer less than or equal to $a$; the notation $\ceil{a}$ denotes the least integer greater than or equal to $a$.  
We establish the following formula for $\dim \widetilde{\T}(x_0)$.

\begin{thm}
\label{thm:dim_tildeT}
Suppose that $x_0\in \L_k(\Omega)$ where $k$ is an integer with $0\leq k\leq \frac{D}{2}$. Then the following statements hold:
\begin{enumerate}
\item Suppose that $0\leq k<\frac{D}{3}$. Then the dimension of $\widetilde{\T}(x_0)$ is equal to 
\begin{gather*}
{k+4\choose 5}
+\frac{1}{2}{k+5\choose 5}
+\frac{1}{4}{k+5\choose 4}
-\frac{1}{8}{k+5\choose 3}
+\frac{1}{16}{k+5\choose 2}
-\frac{1}{16}\floor{\frac{k+5}{2}}.
\end{gather*}

\item Suppose that $\frac{D}{3}\leq k<\frac{D}{2}$. Then the dimension of $\widetilde{\T}(x_0)$ is equal to 
\begin{align*}
&{k+4\choose 5}
+\frac{1}{2}{k+5\choose 5}
+\frac{1}{4}{k+5\choose 4}
-\frac{1}{8}{k+5\choose 3}
+\frac{1}{16}{k+5\choose 2}
-\frac{1}{2}{3k-D+4\choose 5}
\\
&\quad -\;
\frac{1}{4}{3k-D+4\choose 4}
+\frac{1}{8}{3k-D+4\choose 3}
-\frac{1}{16}{3k-D+4\choose 2}
+\frac{1}{16}\left(
\floor{\frac{k}{2}}
-\ceil{\frac{D-k}{2}}
\right).
\end{align*}

\item Suppose that $k=\frac{D}{2}$. Then the dimension of $\widetilde{\T}(x_0)$ is equal to 
$$
\frac{1}{2}{\frac{D}{2}+5\choose 5}
+\frac{1}{4}{\frac{D}{2}+5\choose 4}
-\frac{1}{8}{\frac{D}{2}+5\choose 3}
+\frac{1}{16}{\frac{D}{2}+5\choose 2}
-\frac{1}{16}\floor{\frac{D+10}{4}}.
$$
\end{enumerate}
\end{thm}
\begin{proof}
In this proof we will use the notation 
$$
s_n=\frac{1}{2}{n+1\choose 4}
+\frac{1}{4}{n+1\choose 3}
-\frac{1}{8}{n+1\choose 2}
+\frac{1}{16}{n+1\choose 1}
-
\left\{
\begin{array}{ll}
\frac{1}{16}
\qquad &\hbox{if $n$ is even},
\\
0
\qquad &\hbox{if $n$ is odd}
\end{array}
\right.
$$
for all $n\in \N$ and the equality
\begin{gather}\label{bino}
\sum_{i=0}^n
{i\choose \ell}
={n+1\choose \ell+1}
\qquad 
\hbox{for all $\ell,n\in \N$}.
\end{gather}

(i), (ii): Suppose that $k<\frac{D}{2}$. For any integer $h$ with $0\leq h\leq k$ let 
\begin{align*}
d_h
&=
\sum_{
(h,i,j)\in \mathbf P(k)}
\left(
\min\{D-k-i,k-j\}-\max\{i,j\}-h+1
\right)^2.
\end{align*}
By Theorem \ref{thm:tildeT}(i) the dimension of $\widetilde{\T}(x_0)$ is equal to 
\begin{gather}\label{dimT}
\sum_{h=0}^{k} 
d_h.
\end{gather}
Note that $d_k=1$.
By Theorem \ref{thm:Lmn_Hmodule}(ii) and \cite[Theorems 1.6, 5.7(i) and 5.9]{Huang:CG&Johnson} the number $d_h$ is equal to the dimension of the Terwilliger algebra of the Johnson graph $J(D-2h,k-h)$ for any integer $h$ with $0\leq h\leq k-1$. 
By \cite[Propositions 6.4--6.6]{Huang:CG&Johnson} the number 
$$
d_h=\left\{
\begin{array}{ll}
{k-h+3\choose 4}+s_{k-h+3}
\qquad 
&\hbox{if $3k-D<h\leq k-1$},
\\
{k-h+3\choose 4}+s_{k-h+3}-s_{3k-h-D+2}
\qquad 
&\hbox{if $0\leq h\leq 3k-D$}.
\end{array}
\right.
$$
To obtain (i) and (ii) one may substitute the above equation into (\ref{dimT}) and apply (\ref{bino}) to evaluate the resulting equation.

(iii): Suppose that $k=\frac{D}{2}$. 
For any integer $h$ with $0\leq h\leq \frac{D}{2}$ let 
\begin{align*}
d_h
&=
\sum_{i=0}^{\floor{\frac{D-2h}{4}}}
\textstyle
(i+1)
(\frac{D}{2}-h-2i+1)^2.
\end{align*} 
By Theorem \ref{thm:tildeT}(ii) the dimension of $\widetilde{\T}(x_0)$ is equal to 
\begin{gather}\label{dimT2}
\sum_{h=0}^{\frac{D}{2}} 
d_h.
\end{gather}
Note that $d_{\frac{D}{2}}=1$. 
By Theorem \ref{thm:Lmn_Hmodule}(ii) and \cite[Theorems 1.6, 5.7(ii) and 5.9]{Huang:CG&Johnson} the number $d_h$ is equal to the dimension of the Terwilliger algebra of the Johnson graph $J(D-2h,\frac{D}{2}-h)$ for any integer $h$ with $0\leq h\leq \frac{D}{2}-1$. 
By \cite[Proposition 6.7]{Huang:CG&Johnson} the number $d_h=s_{\frac{D}{2}-h+3}$ for any integer $h$ with $0\leq h\leq \frac{D}{2}-1$. 
The statement (iii) now follows by applying (\ref{bino}) to evaluate the equation (\ref{dimT2}). 
\end{proof}

\section{An application to the Terwilliger algebras of Grassmann graphs}\label{s:H&Grassmann}

Recall the $\H$-module $\C^{\L(\Omega)}(x_0)$ from Definition \ref{defn:Wmodule_L(Omega)(x0)}. By Theorem \ref{thm1:H->W} the action of $A$ on the $\H$-module $\C^{\L(\Omega)}(x_0)$ is identical to the action of $K_2^{-1}$ on the $\W$-module $\C^{\L(\Omega)}(x_0)$. 
Combined with Definition \ref{defn:D12}(ii) and Theorem \ref{thm:Wmodule_L(Omega)} the element $A$ acts on the $\H$-module $\C^{\L(\Omega)}(x_0)$ as follows:
\begin{gather}
\label{A_L(Omega)}
A x = q^{2\dim(x\cap x_0)-\dim x_0} x
\qquad 
\hbox{for all $x\in \L(\Omega)$}.
\end{gather}
By Theorem \ref{thm1:H->W} the action of $B$ on the $\H$-module $\C^{\L(\Omega)}(x_0)$ is identical to the action of $\widetilde{\Delta}(\Lambda)$ on the $\W$-module $\C^{\L(\Omega)}(x_0)$. By Theorem \ref{thm:X&tildeDeltaXq} the latter is also identical to the action of $\Lambda$ on the $\U$-module $\C^{\L(\Omega)}(q^{\dim x_0})$. 
Combined with Lemma \ref{lem:Lambda_L(Omega)} the element $B$ acts on the $\H$-module $\C^{\L(\Omega)}(x_0)$ as follows:
\begin{align}
\label{B_L(Omega)}
B x &= 
(q^{D-2\dim x+1}+q^{2\dim x-D+1}+q^{-1-D}-q^{1-D}) x
\\
&\qquad +\;
q^{1-D}(q-q^{-1})^2\sum_{\substack{x'\in \L_{\dim x}(\Omega) 
\\
x\cap x'\subsetdot x}}
x'
\qquad 
\hbox{for all $x\in \L(\Omega)$}.
\notag
\end{align}

Assume that $k$ is an integer with $1\leq k\leq D-1$. Recall that the {\it Grassmann graph} of $k$-dimensional subspaces of $\Omega$, 
denoted by $J_q(D,k)$ in this paper, is a finite simple connected graph whose vertex set is $\L_k(\Omega)$ and two vertices $x,x'\in \L_k(\Omega)$ are adjacent whenever $x\cap x'\subsetdot x$. In other words, the adjacency operator $\A$ of $J_q(D,k)$ is a linear map $\C^{\L_k(\Omega)}\to \C^{\L_k(\Omega)}$ given by 
\begin{gather}\label{A}
\A x=\sum_{\substack{x'\in \L_k(\Omega)\\ x\cap x'\subsetdot x}} x'
\qquad 
\hbox{for all $x\in \L_k(\Omega)$}.
\end{gather}
Suppose that $x_0\in \L_k(\Omega)$. By \cite{BannaiIto1984,TerAlgebraIII} the dual adjacency operator $\A^*(x_0)$ of $J_q(D,k)$ with respect to $x_0$ is a linear map $\C^{\L_k(\Omega)}\to \C^{\L_k(\Omega)}$ given by 
\begin{gather}\label{A*}
\A^*(x_0) x
=
\frac{[D-1]_q}{q-q^{-1}}
\left(
\frac{[D]_q}{[k]_q[D-k]_q} q^{D+2(\dim x\cap x_0-k)}
-\frac{q^k}{[D-k]_q}
-\frac{q^{D-k}}{[k]_q}
\right) x
\end{gather}
for all $x\in \L_k(\Omega)$. 
The {\it Terwilliger algebra} $\T(x_0)$ of $J_q(D,k)$ with respect to $x_0$ is the subalgebra of ${\rm End}(\C^{\L_k(\Omega)})$ generated by $\A$ and $\A^*(x_0)$.
Please see \cite{TerAlgebraI, TerAlgebraII, TerAlgebraIII, BannaiIto1984} for details. 
Recall the $\H$-submodule $\C^{\L_k(\Omega)}(x_0)$ of $\C^{\L(\Omega)}(x_0)$ from Definition \ref{defn:H_Lk(Omega)(x0)}. By (\ref{A_L(Omega)})--(\ref{A*}) the following equations hold on the $\H$-module $\C^{\L_k(\Omega)}(x_0)$:
\begin{align}
 \A 
&=
\frac{q^{D-1}B-q^{2D-2k}-q^{2k}}{(q-q^{-1})^2}
+\frac{1}{q^2-1},
\label{A&B}
\\
\A^*(x_0)
&=
\frac{[D-1]_q}{q-q^{-1}}
\left(
\frac{q^D[D]_q}{[k]_q[D-k]_q} A
-\frac{q^k}{[D-k]_q}
-\frac{q^{D-k}}{[k]_q}
\right).
\label{A*&A}
\end{align} 
Therefore $\T(x_0)$ is a subalgebra of $\widetilde{\T}(x_0)$.

\begin{lem}
\label{lem:T<tildeT}
Suppose that $x_0\in \L_k(\Omega)$ where $k$ is an integer with $1\leq k\leq D-1$. Then the following statements hold:
\begin{enumerate}
\item Every irreducible $\H$-submodule of $\C^{\L_k(\Omega)}(x_0)$ is an irreducible $\T(x_0)$-submodule of $\C^{\L_k(\Omega)}$.

\item For any $(h,i,j)\in \P(k)$ the irreducible $\H$-module $(L_{D-k-h-2i}\otimes L_{k-h-2j})^{1,q^{-h}}(q^{D-2k})$ is an irreducible $\T(x_0)$-module on which the actions of $\A$ and $\A^*(x_0)$ are given by {\rm (\ref{A&B})} and {\rm (\ref{A*&A})} respectively.
\end{enumerate}
\end{lem}
\begin{proof}
(i): Immediate from the fact that $\T(x_0)$ is a subalgebra of $\widetilde{\T}(x_0)$. 

(ii): Immediate from Proposition \ref{prop:L(Omega)_Hdec} and Lemma \ref{lem:T<tildeT}(i).
\end{proof}

\begin{defn}
\label{defn:P0k}
Assume that $k$ is an integer with $0\leq k\leq D$. 
Let $\P^{(0)}(k)$ denote the set of all triples $(h,i,j)\in \P(k)$ satisfying the following condition:
\begin{gather}\label{P0k}
\min\{D-k-i,k-j\}=h+\max\{i,j\}.
\end{gather}
Equivalently, by Proposition \ref{prop:L(Omega)_Hdec}(ii) the set $\P^{(0)}(k)$ consists of every triple $(h,i,j)\in \P(k)$ for which the space $(L_{D-k-h-2i}\otimes L_{k-h-2j})^{1,q^{-h}}(q^{D-2k})$ is of dimension one.
\end{defn}

\begin{prop}
\label{prop:P-P0}
Suppose that $x_0\in \L_k(\Omega)$ where $k$ is an integer with $1\leq k\leq D-1$. 
For any $(h,i,j)$, $(h',i',j')\in \P(k)\setminus \P^{(0)}(k)$ the following conditions are equivalent:
\begin{enumerate}
\item The irreducible $\T(x_0)$-module $(L_{D-k-h-2i}\otimes L_{k-h-2j})^{1,q^{-h}}(q^{D-2k})$ is isomorphic to the irreducible $\T(x_0)$-module $(L_{D-k-h'-2i'}\otimes L_{k-h'-2j'})^{1,q^{-h'}}(q^{D-2k})$.

\item The irreducible $\H$-module $(L_{D-k-h-2i}\otimes L_{k-h-2j})^{1,q^{-h}}(q^{D-2k})$ is isomorphic to the irreducible $\H$-module $(L_{D-k-h'-2i'}\otimes L_{k-h'-2j'})^{1,q^{-h'}}(q^{D-2k})$.
\end{enumerate}
\end{prop}
\begin{proof}
(ii) $\Rightarrow$ (i): Immediate from Lemma \ref{lem:T<tildeT}(ii). 

(i) $\Rightarrow$ (ii): Immediate from Corollary \ref{cor':irrHmodule} and Lemma \ref{lem:T<tildeT}(ii). 
\end{proof}

\begin{defn}
\label{defn:P0k123}
Assume that $k$ is an integer with $0\leq k\leq \frac{D}{2}$. 
Let 
\begin{align*}
\P^{(0)}_{\rm I}(k)
&=\{(h,i,j)\in \P^{(0)}(k)\,|\,i< j\},
\\
\P^{(0)}_{\rm II}(k)
&=\{(h,i,j)\in \P^{(0)}(k)\,|\,0\leq i- j\leq D-2k\},
\\
\P^{(0)}_{\rm III}(k)
&=\{(h,i,j)\in \P^{(0)}(k)\,|\,D-2k <i- j\}.
\end{align*}
For any integer $j$ with $0\leq j\leq \floor{\frac{k}{2}}$ let 
$
\P^{(0)}_{\rm II}(k;j)
$ 
denote the set of all triples $(h,i,j)\in \P^{(0)}_{\rm II}(k)$. 
\end{defn}

\begin{lem}
\label{lem:P0k1}
Suppose that $k$ is an integer with $0\leq k\leq \frac{D}{2}$. 
Then 
$
k=h+2j
$
for any $(h,i,j)\in \P^{(0)}_{\rm I}(k)$.
\end{lem}
\begin{proof}
Let $(h,i,j)\in \P^{(0)}_{\rm I}(k)$ be given. Using Definition \ref{defn:P0k123} yields that $k-j< D-k-i$ and $i<j$. Combined with (\ref{P0k}) the lemma follows. 
\end{proof}

\begin{lem}
\label{lem:P0k2}
Suppose that $k$ is an integer with $0\leq k\leq \frac{D}{2}$. Then the following statements hold:
\begin{enumerate}
\item $\P^{(0)}_{\rm II}(k)$ is a disjoint union of the sets $\{\P^{(0)}_{\rm II}(k;j)\}_{j=0}^{\floor{\frac{k}{2}}}$.

\item For any integer $j$ with $0\leq j\leq \floor{\frac{k}{2}}$ the set $\P^{(0)}_{\rm II}(k;j)$ consists of 
the triples 
$$
(k-2j-p,j+p,j)
$$
for all integers $p$ with $0\leq p \leq \min\{k-2j,D-2k\}$.

\item $|\P_{\rm II}^{(0)}(k)|$ is equal to 
$$
\left\{
\begin{array}{ll}
(\floor{\frac{k}{2}}+1)(\ceil{\frac{k}{2}}+1)
\qquad 
&\hbox{if $0\leq k<\frac{D}{3}$},
\\
(D-2k+1)(\floor{\frac{3k-D}{2}}+1)+
(\ceil{\frac{D-k}{2}}-\floor{\frac{k}{2}})
(\ceil{\frac{D-k}{2}}-\ceil{\frac{k}{2}})
\qquad 
&\hbox{if $\frac{D}{3}\leq k\leq \frac{D}{2}$}.
\end{array}
\right.
$$
\end{enumerate}
\end{lem}
\begin{proof}
(i): By Definition \ref{defn:Pk}, if $(h,i,j)\in \P(k)$ then $j$ is an integer with $0\leq j\leq \floor{\frac{k}{2}}$.

(ii): Fix an integer $j$ with $0\leq j\leq \floor{\frac{k}{2}}$. Let $(h,i,j)\in \P^{(0)}_{\rm II}(k;j)$ be given.  By Definition \ref{defn:P0k123} there is $p\in \N$ such that $i=j+p$. Using (\ref{P0k}) yields that $h=k-2j-p$. 

Now let $p\in \N$ be given. By Definition \ref{defn:Pk} the triple $(k-2j-p,j+p,j)\in \P(k)$ if and only if $p\leq \min\{k-2j,D-2k\}$. Combined with Definition \ref{defn:P0k123}, if $(k-2j-p,j+p,j)\in \P(k)$ then $(k-2j-p,j+p,j)\in \P^{(0)}_{\rm II}(k;j)$. 
By the above comments the statement (ii) follows.

(iii): It follows from Lemma \ref{lem:P0k2}(i), (ii) that 
$$
|\P_{\rm II}^{(0)}(k)|=\sum_{j=0}^{\floor{\frac{k}{2}}}
(\min\{k-2j,D-2k\}+1).
$$
It is routine to verify (iii) by using the above equation.   
\end{proof}

\begin{lem}
\label{lem:P0k3}
Suppose that $k$ is an integer with $0\leq k\leq \frac{D}{2}$. 
Then the following statements hold:
\begin{enumerate}
\item  $\P^{(0)}_{\rm III}(k)$ consists of all triples $(h,i,j)\in \Z^3$ satisfying the following conditions:
\begin{align*}
&0\leq h< 3k-D,
\\ 
&i=\textstyle \frac{D-k-h}{2},
\\
&0\leq j< \textstyle  \frac{3k-h-D}{2}.
\end{align*}

\item $|\P^{(0)}_{\rm III}(k)|$ is equal to 
$$
\left\{
\begin{array}{ll}
0
\qquad &\hbox{if $0\leq k<\frac{D}{3}$},
\\
\displaystyle {\floor{\frac{3k-D}{2}}+1 \choose 2}
\qquad &\hbox{if $\frac{D}{3}\leq k\leq \frac{D}{2}$}.
\end{array}
\right.
$$

\item There exists a map $\E:\P^{(0)}_{\rm III}(k)\to \P^{(0)}_{\rm I}(k)$ given by 
\begin{eqnarray*}
(h,i,j) &\mapsto & 
\textstyle (h-2k+D,j,\frac{3k-h-D}{2})
\qquad 
\hbox{for all $(h,i,j)\in \P^{(0)}_{\rm III}(k)$}.
\end{eqnarray*}
Moreover $\E$ is injective.
\end{enumerate}
\end{lem}
\begin{proof}
(i): 
Let $(h,i,j)\in \P^{(0)}_{\rm III}(k)$ be given. Using Definition \ref{defn:P0k123} yields that $D-k-i<k-j$ and $j<i$. Combined with (\ref{P0k}) this implies  $i=\frac{D-k-h}{2}$. 

Now let $h,j\in \N$ be given. By Definition \ref{defn:Pk} the triple $(h,\frac{D-k-h}{2},j)\in \P(k)$ if and only if $h\leq 3k-D$ and $j\leq \frac{3k-h-D}{2}$. Suppose that $(h,\frac{D-k-h}{2},j)\in \P(k)$. By Definition \ref{defn:P0k123} the triple $(h,\frac{D-k-h}{2},j)\in \P^{(0)}_{\rm III}(k)$ if and only if $j<\frac{3k-h-D}{2}$. Note that if $h=3k-D$ then there is no $j\in \N$ satisfying $j<\frac{3k-h-D}{2}$. 
By the above comments the statement (i) follows.

(ii): It follows from Lemma \ref{lem:P0k3}(i) that 
$$
|\P^{(0)}_{\rm III}(k)|
=
\sum_{\substack{
h=0\\
h=D-k\bmod{2}
}}^{3k-D} \frac{3k-h-D}{2}.
$$
It is routine to verify (ii) by using the above equation.

(iii): Let $(h,i,j)\in \P^{(0)}_{\rm III}(k)$ be given. Using Lemma \ref{lem:P0k3}(i) it is routine to verify that $(h-2k+D,j,\frac{3k-h-D}{2})\in \P^{(0)}_{\rm I}(k)$. The existence of $\E$ follows. Suppose that $(h,i,j)$ and $(h',i',j')$ are in $\P^{(0)}_{\rm III}(k)$ with $\E(h,i,j)=\E(h',i',j')$. Clearly $(h,j)=(h',j')$. Combined with Lemma \ref{lem:P0k3}(i) this yields $i=i'$. Hence $\E$ is injective.
\end{proof}

\begin{prop}
\label{prop:P0}
Suppose that $x_0\in \L_k(\Omega)$ where $k$ is an integer with $1\leq k\leq \frac{D}{2}$. 
For any $(h,i,j),(h',i',j')\in \P^{(0)}(k)$ the following conditions are equivalent:
\begin{enumerate}
\item The irreducible $\T(x_0)$-module $(L_{D-k-h-2i}\otimes L_{k-h-2j})^{1,q^{-h}}(q^{D-2k})$ is isomorphic to the irreducible $\T(x_0)$-module $(L_{D-k-h'-2i'}\otimes L_{k-h'-2j'})^{1,q^{-h'}}(q^{D-2k})$.

\item The following equations hold:
\begin{align*}
h+\max\{i,j\}&=h'+\max\{i',j'\},
\\
h+i+j&=h'+i'+j'.
\end{align*}
\end{enumerate}
Moreover, if $(h,i,j)\not=(h',i',j')$ then {\rm (i)} and {\rm (ii)} are true if and only if one of the following conditions holds:
\begin{enumerate}
\item[{\rm (a)}] $(h,i,j), (h',i',j')\in \P^{(0)}_{\rm II}(k;j)$.

\item[{\rm (b)}] $(h,i,j)\in \P_{\rm III}^{(0)}(k)$ 
and $\E(h,i,j)=(h',i',j')$.

\item[{\rm (c)}] $(h',i',j')\in \P_{\rm III}^{(0)}(k)$ 
and $\E(h',i',j')=(h,i,j)$. 
\end{enumerate}
\end{prop}
\begin{proof}
Since $(h,i,j)\in \P^{(0)}(k)$ the $\H$-module $(L_{D-k-h-2i}\otimes L_{k-h-2j})^{1,q^{-h}}(q^{D-2k})$ is isomorphic to $V_0(a,b,c)$ where $a,b,c$ are given in Proposition \ref{prop:L(Omega)_Hdec}(ii). Using (\ref{P0k}) yields that 
\begin{gather}\label{(abc)}
(a,b,c)=
(
q^{k-2(h+\max\{i,j\})},q^{D-2(h+i+j)+1},q^{k-D+2\max\{i,j\}}
).
\end{gather}
Similarly the $\H$-module $(L_{D-k-h'-2i'}\otimes L_{k-h'-2j'})^{1,q^{-h'}}(q^{D-2k})$ is isomorphic to $V_0(a',b',c')$ where
\begin{gather}\label{(a'b'c')}
(a',b',c')=
(
q^{k-2(h'+\max\{i',j'\})},q^{D-2(h'+i'+j')+1},q^{k-D+2\max\{i',j'\}}
).
\end{gather}
The element $A$ acts on the $\H$-modules $V_0(a,b,c)$ and $V_0(a',b',c')$ as scalar multiplication by $a$ and $a'$ respectively. The element $B$ acts on the $\H$-modules $V_0(a,b,c)$ and $V_0(a',b',c')$ as scalar multiplication by $b+b^{-1}$ and $b'+b'^{-1}$ respectively. 
Recall from Definition \ref{defn:Pk} that $h+i+j\leq k$ and $h'+i'+j'\leq k$.
Since $k\leq \frac{D}{2}$ it follows that $D-2(h+i+j)\geq 0$ and $D-2(h'+i'+j')\geq 0$. 
Combined with Lemma \ref{lem:T<tildeT}(ii) 
the condition (i) is equivalent to $(a,b)=(a',b')$. By (\ref{(abc)}) and (\ref{(a'b'c')}) the latter is equivalent to the condition (ii). Therefore  (i) and (ii) are equivalent.

Assume that $(h,i,j),(h',i',j')\in \P^{(0)}(k)$ with $(h,i,j)\not=(h',i',j')$. 
By Definition \ref{defn:P0k123} along with Lemmas \ref{lem:P0k2}(ii) and \ref{lem:P0k3}, the condition (a), (b) or (c) implies (ii). Now we suppose that (ii) holds and show that (a), (b) or (c) holds. By symmetry we may divide the argument into the following six cases:

(1): Suppose that $(h,i,j)\in \P^{(0)}_{\rm I}(k)$ and $(h',i',j')\in \P^{(0)}_{\rm II}(k)$. Using Lemmas \ref{lem:P0k1} and \ref{lem:P0k2}(ii) yields that 
$$
k=h+2j=h'+i'+j'.
$$
Combined with the second equation in (ii) this implies $i=j$, a contradiction to $(h,i,j)\in \P^{(0)}_{\rm I}(k)$.

(2): Suppose that $(h,i,j)\in \P^{(0)}_{\rm II}(k)$ and $(h',i',j')\in \P^{(0)}_{\rm III}(k)$. 
Using Lemma \ref{lem:P0k2}(ii) yields that 
$k=h+i+j$. Combined with the second equation in (ii) the following equation holds:
$$
k=h'+i'+j'.
$$
Recall from Lemma \ref{lem:P0k3}(i) that $i'=\frac{D-k-h'}{2}$. Substituting the expression for $i'$ into the above equation yields 
$j'=\frac{3k-h'-D}{2}$, a contradiction to Lemma \ref{lem:P0k3}(i).

(3): Suppose that $(h,i,j)\in \P^{(0)}_{\rm I}(k)$ and $(h',i',j')\in \P^{(0)}_{\rm III}(k)$. By Definition \ref{defn:P0k123} the equations in (ii) are as follows:
\begin{gather*}
h+j=h'+i',
\\
h+i+j=h'+i'+j'.
\end{gather*}
Using Lemmas \ref{lem:P0k1} and \ref{lem:P0k3}(i) yields that 
$$
k=h+2j=D-h'-2i'.
$$
Solving the above three equations yields that $\E(h',i',j')=(h,i,j)$. 
Therefore the condition (c) holds. Note that the condition (b) holds by switching the above roles of $(h,i,j)$ and $(h',i',j')$.

(4): Suppose that $(h,i,j), (h',i',j')\in \P^{(0)}_{\rm I}(k)$. By Definition \ref{defn:P0k123} the equations in (ii) are as follows:
\begin{gather*}
h+j=h'+j',
\\
h+i+j=h'+i'+j'.
\end{gather*}
By Lemma \ref{lem:P0k1} the following equation holds:
$$
h+2j=h'+2j'.
$$ 
Solving the above three equations yields that $(h',i',j')=(h,i,j)$, a contradiction.

(5): Suppose that $(h,i,j),(h',i',j')\in \P^{(0)}_{\rm II}(k)$. 
By Lemma \ref{lem:P0k2}(ii) there are 
$p,p'\in \N$ such that 
\begin{align*}
h&=k-2j-p,
\qquad 
i=j+p,
\\
h'&=k-2j'-p',
\qquad 
i'=j'+p'.
\end{align*}
Substituting the above equations into the first equation given in (ii) yields $j=j'$. Therefore the condition (a) holds.

(6): Suppose that $(h,i,j),(h',i',j')\in \P^{(0)}_{\rm III}(k)$. 
By Definition \ref{defn:P0k123} the equations in (ii) are as follows:
\begin{gather*}
h+i=h'+i',
\\
h+i+j=h'+i'+j'.
\end{gather*}
By Lemma \ref{lem:P0k3}(i) the following equation holds:
$$
h+2i=h'+2i'.
$$
Solving the above three equations yields that $(h',i',j')=(h,i,j)$, a contradiction. 
The proposition follows.
\end{proof}

Mapping each $x\in \L_k(\Omega)$ to the orthogonal complement of $x$ in $\Omega$ with respect to a nondegenerate bilinear form, it induces a graph isomorphism from $J_q(D,k)$ to $J_q(D,D-k)$. By (\ref{A}) and (\ref{A*}) it suffices to study the algebra $\T(x_0)$ under the assumption $x_0\in \L_k(\Omega)$ with $1\leq k\leq \frac{D}{2}$.
We are all set to present the decomposition of $\T(x_0)$-module $\C^{\L_k(\Omega)}$ and the decomposition of $\T(x_0)$ well as the formula for $\dim \T(x_0)$.

\begin{thm}
\label{thm:T-module}
Suppose that $x_0\in \L_k(\Omega)$ where $k$ is an integer with $1\leq k\leq \frac{D}{2}$. Then the following statements hold:
\begin{enumerate}
\item Suppose that $k<\frac{D}{2}$. Then the $\T(x_0)$-module $\C^{\L_k(\Omega)}$ is isomorphic to 
\begin{align*}
&\bigoplus_{
\substack{
(h,i,j)\in \P(k)
\\
(h,i,j)\notin \E(\P^{(0)}_{\rm III}(k))
\cup \P^{(0)}_{\rm II}(k)
\cup \P^{(0)}_{\rm III}(k)
}
}
m_{hij}(k)\cdot 
(L_{D-k-h-2i}\otimes L_{k-h-2j})^{1,q^{-h}}
(q^{D-2k})
\\
&\oplus
\quad
\,\,\,
\bigoplus_{j=0}^{\floor{\frac{k}{2}}}
\bigg(
\sum_{(h,i,j)\in \P^{(0)}_{\rm II}(k;j)}
m_{hij}(k)
\bigg)
\cdot 
(L_{D-2k}\otimes L_0)^{1,q^{2j-k}}
(q^{D-2k})
\\
&\oplus
\bigoplus_{\substack{
(h,i,j)\in \P^{(0)}_{\rm III}(k)
}}
\left(
m_{hij}(k)+m_{\E(h,i,j)}(k)
\right)\cdot 
(L_0\otimes L_{k-h-2j})^{1,q^{-h}}
(q^{D-2k}).
\end{align*} 
Moreover 
the above irreducible $\T(x_0)$-modules 
are mutually non-isomorphic.

\item Suppose that $k=\frac{D}{2}$. Then the $\T(x_0)$-module $\C^{\L_k(\Omega)}$ is isomorphic to 
$$
\bigoplus_{h=0}^{\frac{D}{2}}
\bigoplus_{i=0}^{\floor{\frac{D-2h}{4}}}
\bigoplus_{j=i}^{\floor{\frac{D-2h}{4}}}
m_{hij}\cdot
(L_{\frac{D}{2}-h-2i}\otimes L_{\frac{D}{2}-h-2j})^{1,q^{-h}}
(1)
$$
where 
$$
m_{hij}
=\left\{
\begin{array}{ll}
2m_{hij}(\frac{D}{2}) \qquad 
&\hbox{if $i<j$},
\\
m_{hii}(\frac{D}{2}) \qquad 
&\hbox{if $i=j$}.
\end{array}
\right.
$$
Moreover 
the above irreducible $\T(x_0)$-modules 
are mutually non-isomorphic. 
\end{enumerate}
\end{thm}
\begin{proof}
Immediate from Theorem \ref{thm:L(Omega)_Hdec} and Propositions \ref{prop:P-P0} and \ref{prop:P0}.
\end{proof}

\begin{thm}
\label{thm:T}
Suppose that $x_0\in \L_k(\Omega)$ where $k$ is an integer with $1\leq k\leq \frac{D}{2}$. Then the following statements hold:
\begin{enumerate}
\item Suppose that $k<\frac{D}{2}$. Then the algebra $\T(x_0)$ is isomorphic to 
\begin{gather*}
({\textstyle\floor{\frac{k}{2}}}+1)\cdot \C
\oplus
\bigoplus_{(h,i,j)\in \P(k)\setminus(\P^{(0)}_{\rm II}(k)\cup \P^{(0)}_{\rm III}(k))}
{\rm End}(\C^{\min\{D-k-i,k-j\}-\max\{i,j\}-h+1}).
\end{gather*}

\item Suppose that $k=\frac{D}{2}$. The algebra $\T(x_0)$ is isomorphic to 
\begin{gather*}
\bigoplus_{h=0}^{\frac{D}{2}}
\bigoplus_{i=0}^{\floor{\frac{D-2h}{4}}}
(i+1)\cdot {\rm End}(\C^{\frac{D}{2}-h-2i+1}).
\end{gather*}
\end{enumerate}
\end{thm}
\begin{proof}
Immediate from Theorem \ref{thm:T-module}.
\end{proof}

Recall the closed-form formula for $\dim \widetilde{\T}(x_0)$ from Theorem \ref{thm:dim_tildeT}. Instead of a complete description for $\dim \T(x_0)$ we give the formula for $\dim \widetilde{\T}(x_0)-\dim \T(x_0)$.

\begin{thm}
\label{thm:dim_T}
Suppose that $x_0\in \L_k(\Omega)$ where $k$ is an integer with $1\leq k\leq \frac{D}{2}$. Then the following statements hold:
\begin{enumerate}
\item Suppose that $1\leq k<\frac{D}{3}$. Then $\dim \widetilde{\T}(x_0)-\dim \T(x_0)=
\ceil{\frac{k}{2}}
\left(
\floor{\frac{k}{2}}+1
\right)$.

\item Suppose that $\frac{D}{3}\leq k<\frac{D}{2}$. Then  $\dim \widetilde{\T}(x_0)-\dim \T(x_0)$ is equal to 
\begin{align*}
&
(D-2k)
\left(
\floor{\frac{3k-D}{2}}+1
\right)
+
\left(\ceil{\frac{D-k}{2}}-\floor{\frac{k}{2}}\right)
\left(\ceil{\frac{D-k}{2}}-\ceil{\frac{k}{2}}\right)
\\
&\qquad +\,
{\floor{\frac{3k-D}{2}}+2\choose 2}
-
\floor{\frac{k}{2}}-1.
\end{align*}

\item Suppose that $k=\frac{D}{2}$. Then $\dim \T(x_0)=\dim \widetilde{\T}(x_0)$.
\end{enumerate}
\end{thm}
\begin{proof}
Comparing Theorem \ref{thm:tildeT} with Theorem \ref{thm:T} yields that 
\begin{gather*}
\dim \widetilde{\T}(x_0)
=
\dim \T(x_0)
+
\left\{
\begin{array}{ll}
|\P^{(0)}_{\rm II}(k)|+|\P^{(0)}_{\rm III}(k)|
-{\textstyle\floor{\frac{k}{2}}}-1
\qquad 
&\hbox{if $1\leq k<\frac{D}{2}$},
\\
0
&\hbox{if $k=\frac{D}{2}$}.
\\
\end{array}
\right.
\end{gather*}
Combined with Lemmas \ref{lem:P0k2}(iii) and \ref{lem:P0k3}(ii) the result follows.
\end{proof}

\subsection*{Acknowledgements}
The research was supported by the National Science and Technology Council of Taiwan under the project NSTC 112-2115-M-008-009-MY2.

\bibliographystyle{amsplain}
\bibliography{MP}

\end{document}